\documentclass[11pt,a4paper,reqno]{amsart}
\usepackage{amsmath,amssymb,amsfonts,amsthm,mathtools}
\usepackage{setspace}
\usepackage{amsthm}
\usepackage{array,multirow,graphicx}
 \usepackage{float}
 \usepackage{epsfig}
 \usepackage{subcaption}
\usepackage{graphicx}
\usepackage{framed}
\usepackage[ruled,noline]{algorithm2e}
\usepackage[utf8]{inputenc}
\usepackage{multirow}
\usepackage[mathscr]{euscript}
\usepackage{tablefootnote}
\usepackage[margin=1in]{geometry}
\usepackage{enumitem}  
\usepackage[hidelinks,breaklinks]{hyperref}
\hypersetup{colorlinks = true, urlcolor = blue, linkcolor = blue,citecolor = black}
\usepackage{natbib}
\bibpunct[, ]{(}{)}{,}{a}{}{,}%
\usepackage{natbib}
\bibpunct[, ]{(}{)}{,}{a}{}{,}%
\renewcommand{\qed}{\hfill $\square$}

\newcommand{\xmath}[1]{\ensuremath{#1}\xspace}
\newcommand{\Lev}{\xmath{\textnormal{lev}}}
\newcommand{\RV}{\xmath{\mathscr{RV}}}

\newcommand{\R}{\xmath{\mathbb{R}}}

\newcommand{\mv}[1]{\boldsymbol{#1}}

\renewcommand{\Pr}{\mathbb{P}}

\newcommand{\CCPalpha}{{\tt{CCP}}(\alpha)}
\newcommand{\Dapprox}{{\tt{CCPapx}}}
\newcommand{\DROCCPalpha}{{\tt{DRO-CCP}}(\alpha)}

\newcommand{\xx}{\boldsymbol{x}}
\newcommand{\XX}{\boldsymbol{X}}
\newcommand{\YY}{\boldsymbol{Y}}

\newcommand{\yy}{\boldsymbol{y}}
\newcommand{\zz}{\boldsymbol{z}}
\newcommand{\xxi}{\boldsymbol{\xi}}

\newcommand{\pp}{\boldsymbol{p}}

\newcommand{\e}{\mathrm{e}}

\newcommand{\Real}{\mathbb{R}}
\newcommand{\Prob}{\mathbb{P}}
\newcommand{\E}{\mathbb{E}}

\newcommand{\pphi}{\boldsymbol{\phi}}

\renewcommand{\Pr}{\xmath{P}}

\renewcommand*{\Upsilon}{\mathcal{U}}

\newtheorem{assumption}{Assumption}
\newtheorem{theorem}{Theorem}[section]
\newtheorem{definition}{Definition}
\newtheorem{example}{Example}
\newtheorem{lemma}{Lemma}[section]
\newtheorem{proposition}{Proposition}[section]
\newtheorem{corollary}{Corollary}[section]

\newtheorem*{assumption*}{Assumption}
\newtheorem{num_example}{Numerical Illustration}

\newcounter{nappendix}

\title[Scaling Behaviors in Achieving High Reliability]{The Scaling Behaviors in Achieving High Reliability via Chance-Constrained Optimization}
\author{{\large A\MakeLowercase{nand}
    D\MakeLowercase{eo}}}
  \address{Indian Institute of Management, Bannerghatta Road, Bangalore 560076} \email{anand.deo@iimb.ac.in}
  \author{{\large K\MakeLowercase{arthyek}
    M\MakeLowercase{urthy}}}
  \address{Singapore University of Technology and Design, 8 Somapah Rd, Singapore 487372} \email{karthyek\_murthy@sutd.edu.sg}
\begin{document}



\begin{abstract}

{We study the problem of resource provisioning under stringent reliability or service level requirements, which arise in applications such as power distribution, emergency response, cloud server provisioning, and regulatory risk management.  With chance-constrained optimization serving as a natural starting point for modeling this class of problems, our primary   contribution is to characterize how the optimal costs and decisions scale for a generic joint chance-constrained model as the target probability of satisfying the service/reliability constraints  approaches its maximal level. Beyond providing insights into the behavior of optimal solutions, our scaling framework has three key algorithmic implications. First, in distributionally robust optimization (DRO) modeling of chance constraints, we show that widely used approaches based on KL-divergences, Wasserstein distances, and moments heavily distort the scaling properties of optimal decisions, leading to exponentially higher costs. In contrast, incorporating marginal distributions or using appropriately chosen 
$f$-divergence balls preserves the correct scaling, ensuring decisions remain conservative by at most a constant or logarithmic factor. Second, we leverage the scaling framework to quantify the conservativeness of common inner approximations and propose a simple line search to refine their solutions, yielding near-optimal decisions. Finally, given 
$N$ data samples, we demonstrate how the scaling framework enables the estimation of approximately Pareto-optimal decisions with constraint violation probabilities significantly smaller than the $\Omega(1/N)$-barrier that arises in the absence of parametric assumptions.
}\\

\end{abstract}
\maketitle
\vspace{-20pt}

\section{Introduction}
Consider a service provider striving to meet a target level of service in the face of uncertainty. For instance, a distributor of a commodity like electricity needs to fulfill uncertain demands at different nodes of a distribution network. The distribution firm must allocate supply capacities to various nodes in such a way that excessive demand shedding occurs in \textsc{no} more than an $\alpha$ fraction of the service instances, where $1 - \alpha \in (0,1)$ is a pre-specified service level agreement, typically close to 1. Similarly, an emergency response service provider might seek to minimize the costs of positioning and dispatching ambulances while ensuring that the uncertain spatially distributed  demand for emergency services is met with probability $1 - \alpha$. What is the minimum cost the service provider will incur in meeting such high service level agreements? Specifically, how large must this minimum cost and the optimal resource allocations be if the provider aims to ensure very high levels of service availability, as required in contexts like electricity distribution, cloud computing, or emergency medical response? This paper focuses on analytically treating this question, as  gaining a more explicit understanding of how various cost parameters and the uncertainty influence this minimum cost is a crucial first step towards gaining insights into the economics of high reliability.

\subsection{Chance constrained optimization} \label{sec:ccp-intro}
For several decades, chance constrained optimization   has served as a typical starting point for modeling the above class of problems; see  \cite{charnes1959chance, prekopa1970probabilistic,shapiro2021lectures}. A generic chance-constrained optimization formulation can be stated abstractly as, 
\begin{equation}
    \label{eqn:CCP}
\boldsymbol{\CCPalpha:} \qquad  \min_{\xx \in \mathcal X} \ c(\xx) \quad \text{ s.t. } \quad  \mathbb P\left\{ g_k(\xx, \xxi) \leq 0, \ k = 1,\ldots, K  \right\} \geq 1-\alpha,
\end{equation}
where the goal is to find a decision $\xx$ from  a set $\mathcal X\subset \mathbb R^m$ that minimizes a cost  function $c(\xx)$  while ensuring that the service/reliability constraints, modeled via $\{g_k(\xx,\xxi) \leq 0, k = 1,\ldots,K\},$ are satisfied with high probability despite the parameters $\xxi$ affecting the constraints being random. 
In this paper, we will be primarily interested in the case where $\xxi$ is a continuous $\R^d$-valued random vector admitting a probability density function. For any decision choice $\xx \in \mathcal{X},$ constraint violation $\{g_k(\xx,\xxi) > 0\}$ models an undesirable disruption event such as  excessive demand shedding or a failure to dispatch an ambulance within a 15-minute window. 

To gain a comprehensive view into how the formulation \eqref{eqn:CCP} serves as a powerful  vehicle for modeling high service availability or high reliability requirements in various  contexts,  refer sample applications in power systems and electricity markets  \citep{bienstock2014chance,6714594,pena2020dc}, cloud computing \citep{cohen_cloud,kwon_cloud},  emergency medical service \citep{ReVelle_Hogan, beraldi2004designing}, portfolio selection \citep{agnew_1969,el_ghaoui_robust_var,bonami2009exact}, healthcare management \citep{deng2016decomposition,want_li_mehrotra}, project management \citep{shen_smith_ahmed}, telecommunication networks \citep{5437169}, supply chain and logistics  \citep{jiamin_supply,multisourcing}, humanitarian relief operations \citep{ELCI201855}, and staffing call centres \citep{gurvich_luedtke_tezcan}.  

While the formulation $\CCPalpha$ in  \eqref{eqn:CCP} is conceptually attractive and conducive for quantitatively modeling contractual service-level agreements, it does not lend itself readily to tractable solution procedures. To begin with, observe that 
the potential non-convexity of the collection of decisions $\xx \in \mathcal{X}$ satisfying the   probability constraint renders \eqref{eqn:CCP} computationally intractable, generally speaking.  Beyond a narrow collection of problems (see, eg., \citealt{prekopa1998programming,dentcheva2000concavity,lagoa2005probabilistically}), it is well-known that instances of \eqref{eqn:CCP} where both (i) the efficient computation of the probability $\mathbb{P}[\xi \notin \Upsilon(\xx)]$ and (ii) the convexity of the feasible set hold simultaneously are rare.
Therefore,  clever algorithmic inventions have been necessary to computationally handle the $\CCPalpha$ formulation \eqref{eqn:CCP} in general. These include the use of convex inner approximations (see, eg., \citealt{nemirovski2007convex}), scenario approximation \citep{calafiore2006scenario}, strengthening mixed-integer program formulations (see, eg., \citealt{luedtke2010integer}), among many notable algorithmic developments 
for tackling  \eqref{eqn:CCP}. 
Despite the broader computational challenges, chance-constrained optimization remains a widely used model for decision-making under uncertainty. 


\subsection{Research questions and our contributions}
Diverging from the above discussed methodological research thrusts,  this paper aims to develop a qualitative understanding of the properties of the optimal value and optimal solutions to $\CCPalpha$  when aiming for a high degree of reliability, specified via a target reliability level $1-\alpha$ close to one in \eqref{eqn:CCP}. To illustrate this pursuit concretely, 
consider the question: ``By how much should an electricity distribution firm increase the generator/supply capacities at different nodes of a network if it aims to halve the likelihood of excess demand shedding, say, from $\alpha = 1/1000$ to $1/2000$?" 
Alternatively, how steeply does the capital expenditure 
increase as a function of the target level $1 - \alpha$ with which we wish to avoid network failures due to excess demand shedding? 
Currently, aside from  methods—which become computationally demanding and less accurate under stringent reliability requirements—we lack a means to qualitatively understand how different cost parameters and uncertainties influence the answers to such questions. However, gaining a qualitative understanding into these questions is equally crucial from an operations and risk management perspective.

In this paper, we aim to alleviate this challenge by analytically examining the  chance-constrained formulation $\CCPalpha$ in the  
high reliability regime where $1-\alpha \rightarrow 1.$ 
With large deviations theory (see, eg., \citealt{dembo2009large}) providing a systematic framework for studying how significant deviations in the behaviour of the random vector $\xxi$ leads to  atypical events,  we propose to treat the  probability constraint in $\CCPalpha$ under the lens of large deviations approximations. 
When the target service level $1-\alpha$ approaches one,  this approximation allows us to view $\CCPalpha$ as a perturbed and scaled version of a limiting optimization problem that we can explicitly write. Leveraging the rich literature on perturbation analysis of optimization models \citep{rockafellar2009variational, bonnans2013perturbation}, 
we uncover a remarkable regularity in the behavior of optimal values and  solutions to $\CCPalpha$ when the  distribution of $\xxi$ and the constraint functions $g_k$ in \eqref{eqn:CCP} admit sufficient regularity. Besides offering qualitative insights, this regularity  has several algorithmic implications. The rest of this section is dedicated to describing these contributions and discussing related literature.

\subsubsection{Scaling of optimal costs and  decisions under high reliability requirements.}
Let $v_\alpha^\ast$ denote the optimal value of $\CCPalpha$ and $\xx_{\alpha}^\ast$ denote an optimal solution to $\CCPalpha.$ 
As our first main contribution, we explicitly identify a suitable scaling function $\alpha \mapsto s_\alpha$ which increases to infinity 
and  a constant $r$ such that,   
\begin{align}
\label{eqn:CCP_conv}
v_{\alpha}^* \sim v^* s_{\alpha}^r  \quad \text{ and } \quad {\xx^*_\alpha} \sim \xx^* s_{\alpha}^r,
\end{align}
as the permissible probability of service disruption, denoted by $\alpha$ in $\CCPalpha$, is decreased to zero. Here the notation $\sim$ is used to indicate concisely that $v_\alpha^\ast/s_\alpha^r \rightarrow v^\ast$ and $\xx_\alpha^\ast/s_\alpha^r \rightarrow \xx^\ast,$ as $\alpha$ decreases to zero.  The precise limiting notion, the choice of the scaling function $s_\alpha,$ constant $r,$ and the non-zero limits $v^\ast \in  \R$ and $\xx^\ast \in \R^n$ are identified in Theorem \ref{thm:solution_convergence}. 
A key observation here is that  it is possible to explicitly  characterize the rate $s_\alpha^r$ at which the optimal value  and the optimal solution of $\CCPalpha$ scale as a function of the target  level $1-\alpha$ approach one, and it depends 
 on the probability density of $\xxi$ primarily via the marginal distribution of its components $\xi_1,\ldots,\xi_d.$ 

The observation that the rate $s_\alpha^r,$ at which the optimal cost of meeting high reliability scales, is free of the dependence structure  across the components  $\xi_1,\ldots,\xi_d$ might offer  a degree of relief for decision-makers who may have to estimate the joint distribution of   $\xxi = (\xi_1,\ldots,\xi_d)$ from data. Estimating the joint distribution well in the tail regions where  constraint violations happen   typically is a formidable statistical challenge. The results reveal that although misspecifying the copula may yet affect the constraint satisfaction, its impact on the prescribed decisions and their costs remains bounded when viewed as a function of $\alpha$. 

\subsubsection{Algorithmic implications of the scaling phenomenon \eqref{eqn:CCP_conv}.}
Similar to how analyzing the scaling of computational effort with problem size aids in discriminating and developing efficient algorithms in computer science,   the novel approach of analyzing the optimization instances 
$\CCPalpha$ as a function of the target level 
$1-\alpha$ carries the following novel algorithmic implications. 

\textbf{Application 1: Delineating chance-constrained DRO models with sharp characterizations of  their conservativeness.}  Given a collection of plausible probability distributions $\mathcal{P}$ for the random vector $\xxi,$ the distributionally robust optimization (DRO) approach towards tackling chance constraints involves replacing the probability constraint in \eqref{eqn:CCP} with the uniform requirement $\inf_{\mathbb Q \in \mathcal P} \,  \mathbb{Q} \big\{ g_k(\xx,\xxi) \leq 0, k \in [K] \big\} \geq 1-\alpha.$  
Motivated by considerations of tractability and finite-sample guarantees, 
the distributional ambiguity set $\mathcal{P}$ is  formulated typically via moment constraints (see, eg.,  \citealt{el_ghaoui_robust_var, natarajan_robust_var, grani_ccp}), $f$-divergence balls (see,  eg., \citealt{jiang2016data}), or Wasserstein balls (see \citealt{xie2021distributionally,ho2022distributionally,chen2024data}); see also the survey article \cite{KUCUKYAVUZ2022100030}. While the uniform requirement $\inf_{\mathbb Q \in \mathcal P} \,  \mathbb{Q} \big\{ \cdot \big\} \geq 1-\alpha$ is conceptually a conservative approach, some choices of $\mathcal{P}$ 
are intuitively considered more conservative that others. With the performance of a DRO model relying crucially  on the choice of the ambiguity set $\mathcal{P},$ can we characterize the cost scaling of the optimal decisions prescribed by the DRO models for different choices of $\mathcal{P},$ so that we can  precisely delineate them based on their conservativeness? 

Employing the aforementioned large-deviations machinery, another key contribution in this paper is to  characterize the scaling of the optimal decisions prescribed by some prominent DRO formulations together with their costs. The results reveal that  commonly used DRO approaches based on KL-divergences, Wasserstein distances, and moments heavily distort the scaling properties of optimal decisions and result in exponentially higher costs. On the other hand, incorporating marginal distributions (or) employing suitably chosen $f$-divergence balls preserves the correct scaling and ensures that their solutions remain conservative by at most a constant or logarithmic  factor. 


\textbf{Application 2: Quantifying and reducing the conservativeness of inner approximations based on CVaR and Bonferroni inequality.} 
In computationally handling the joint probability constraints in \eqref{eqn:CCP}, solving inner approximations based on conditional value-at-risk (CVaR) and Bonferroni inequality have served as two widely used approaches (eg., \citealt{nemirovski2007convex}).  While  solutions  from  the inner  approximations remains feasible for $\CCPalpha,$ the extent of optimality loss is less understood—specifically, how much more expensive are their solutions? Using our large deviations-based scaling framework, we develop sharp  characterizations akin to \eqref{eqn:CCP_conv}, providing key insights into their conservativeness. The results reveal that both CVaR- and Bonferroni-based inner approximations exhibit negligible relative loss in optimality when the underlying random variables are light-tailed. In the presence of heavy-tailed random variables, both methods yield decisions that are more expensive by a constant factor. Notably, CVaR-based approximations score better by producing decisions proportionate to those of $\CCPalpha$: Specifically, we establish  
\begin{align}
\label{eqn:CVaR_conv}
v_{\alpha}^{\tt cvar} \sim cv^* s_{\alpha}^r  \quad \text{ and } \quad {\xx^{\tt cvar}_\alpha} \sim c\xx^* s_{\alpha}^r,
\end{align}
where constant $c \geq 1,$ $v_{\alpha}^{\tt cvar}$ and $\xx_{\alpha}^{\tt cvar}$ denote the  optimal  value and solution of the CVaR-constrained  approximation of $\CCPalpha,$ and $v^\ast,\xx^\ast, s_\alpha^r$ are the same as in \eqref{eqn:CCP_conv}.  Although the  conservativeness, quantified by  $c,$ increases with heavier-tailed distributions, we propose an elementary  line-search (see Algorithm \ref{algo:Vanish_Regret})  that can strictly improve the solution $\xx_\alpha^{\text{cvar}}$  into near-optimal decisions for $\CCPalpha$ irrespective of  tail heaviness.  

\begin{figure}[h]  \includegraphics[width=0.95\textwidth]{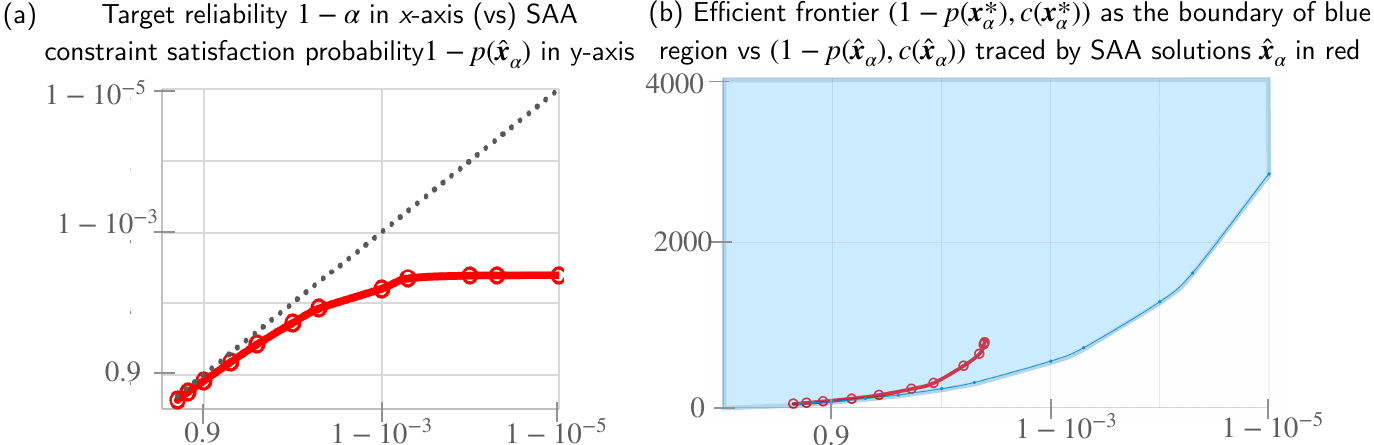}
  \caption{Performance of SAA optimal solutions $\hat{\xx}_\alpha$ obtained at target reliability levels $1-\alpha \in [0.8, 1-10^{-5}]$  using $N = 1000$ samples
    in a transportation example considered in Numerical Illustration \ref{num-eg:DD-Pareto-Front} (Section \ref{sec:Algo})}
\label{fig:dd_underestimate}
\end{figure}

\textbf{Application 3: Estimating  Pareto efficient  decisions from limited data.} In data-driven applications, the choice of the  distribution $\Prob$ in the formulation $\CCPalpha$ is informed typically by a limited number of independent observations of $\xxi.$ Given $N$ such samples and a target constraint violation probability $\alpha < 1/N$, observe that it is statistically impossible to nonparametrically  identify  a decision $\xx$ whose constraint violation probability $p(\xx) = \Prob\{g_k(\xx,\xxi) > 0 \text { for some } k = 1,\ldots,K\}$ is smaller than $1/N$ unless additional  assumptions are made about the  distribution of $\xxi.$ This is manifestly observed in any sample average approximation (SAA) of the chance constraint $\CCPalpha$ in which target $\alpha <  1/N,$ as illustrated via a numerical example in Figure \ref{fig:dd_underestimate}  below with $N = 1000$ observations. In Fig \ref{fig:dd_underestimate}(a), we see that an SAA optimal solution $\hat{\xx}_\alpha$ falls significantly short of satisfying the target reliability level $1-\alpha$ for all $\alpha < 10/N = 0.01.$ Next, we also observe from both the panels in Fig \ref{fig:dd_underestimate} that SAA is not able to arrive at a solution whose constraint violation probability is smaller than the level $1/N$ regardless of how large we set the target $1-\alpha$ to be. 

The underestimation reported in Figure \ref{fig:dd_underestimate} happens  because for any patently infeasible decision $\xx$ whose constraint violation probability $p(\xx)$ is in the range $(\alpha, 1/N),$ we have $N \times p(\xx) < 1$ samples falling, on an average, in the constraint violation region.
As a result, SAA is prone to dangerously underestimating the constraint violation probability to be zero and declare such an $\xx$ to be feasible, at least 50\% of the times. This is a fundamental statistical bottleneck which limits the use of  nonparameteric estimation approaches in applications requiring high reliability levels where $\alpha \ll 1/N.$ In an attempt to overcome this bottleneck, at least partially, we ask, ``can we develop a nonparametric estimation procedure that, under minimal assumptions, can yield decisions $\xx$ that are Pareto efficient in balancing the cost $c(\xx)$ and the constraint violation probability $p(\xx) = \Prob\{g_k(\xx,\xxi) > 0 \text{ for some } k = 1,\ldots,K\},$ even when $p(\xx) \ll 1/N$?''


As our final illustration of the utility of the large-deviations scaling framework, we demonstrate how \eqref{eqn:CCP_conv} can serve as a basis for addressing the above ambitious question, much like how the field of extreme value theory in statistics provides a rigorous framework for estimating quantiles at levels far beyond what is feasible with finite data \citep[see, eg.,][Chap.~1,~4]{deHaan}. Specifically, under a minimal nonparametric assumption on the distribution of $\xxi,$ we show in Section \ref{sec:Algo}  that a suitably  extrapolated trajectory of solutions, grounded in \eqref{eqn:CCP_conv}, is nearly optimal in minimizing the constraint violation probability $p(\xx)$ for any given cost target, even when the cost target is sufficiently large to allow $p(\xx) \ll 1/N.$

\subsection{Related literature utilizing large deviations theory in  optimization modeling} 
While large deviations theory is frequently used to analyze the quality of solutions from sampled approximations (see, eg., \citealt{shapiro_homemdemello}), its direct application in formulating or studying  optimization models is relatively limited. 
\cite{vanparys_kuhn,sutter_kuhn, sutter_mengmeng_kuhn} 
apply Sanov's theorem--a key result in large deviations theory--to identify data-driven formulations that optimally balance conservativeness of solutions with out-of-sample performance.
More closely  related to our pursuit  are \cite{Mainik}, \cite{8735917,tong2022optimization}, and \cite{blanchet2024optimization}. 
For linear portfolios comprising assets with heavy-tailed losses, \cite{Mainik} seek to minimize their extremal risk index, a notion 
arising from  large deviations approximation of the excess losses probabilities. 
\cite{8735917,tong2022optimization} approximate  chance-constraints using large deviations heuristics for suitably light-tailed random vectors, focusing on the computational aspects of solving the resulting bi-level problems. 

The recent independent study by \cite{blanchet2024optimization}, made public in arXiv about a month before the first version of this paper, shares our objectives of (i) characterizing the scaling of the optimal value in  chance-constrained  models and (ii) assessing the conservativeness of CVaR approximations. Their analysis focuses on a specific case of \eqref{eqn:CCP},  where $g_k(\xx,\xxi) = \xx^\intercal A_k \xxi - 1,$ for $k =1,\ldots,K,$ and additionally quantifies the quality of solutions provided by scenario approximation (\citealt{calafiore2006scenario}). 
We now outline the key differences: 
Fundamentally, in the model studied by  \cite{blanchet2024optimization}, both the optimal value $v_\alpha^\ast$ and the optimal resource allocation decisions $\xx_\alpha^\ast$ shrink to zero as the target reliability level $1-\alpha$ is raised to 1, implying that the scaling in 
\eqref{eqn:CCP_conv} must satisfy  $s_\alpha^r \rightarrow 0$. This qualitative phenomenon contrasts sharply with those in common resource provisioning tasks such as in  commodity distribution networks, emergency response,  or cloud servers, etc., where resource allocations must scale up to meet stricter service level requirements.  One of the   contributions of our paper lies  in identifying minimal structural assumptions on the constraint functions $g_k$ under which both  these qualitative phenomena can manifest. In particular,  our abstraction provides the flexibility to study  both the setting considered by \cite{blanchet2024optimization} and the more common setting where the optimal value of $\CCPalpha$  must  scale up  as  the target  level $1-\alpha$ is raised. Our framework is versatile enough to accommodate a wider 
variety of chance constraints, including those with nonlinearities and uncertainties on either the  left- or right-hand side 
quadratic constraints, and more.  
Additionally, the results characterizing the conservativeness of DRO models,  the algorithmic breakthrough of data-driven estimation of Pareto-optimal decisions even when $\alpha \ll 1/N,$ and the line-search procedure for reducing the conservativeness of CVaR approximation, are all novel and unique to this paper.


\textbf{Organization of the paper.} Section~\ref{sec:Modelling Assumptions} introduces the  precise chance-constrained model assumptions under which we derive our results. Section~\ref{sec:Asymptotics} is devoted to discussing the first main result on the scaling of optimal values and decisions. 
Section~\ref{sec:Robust_Asymototics} delves into the implications of the scaling framework for DRO modeling. 
Sections~\ref{sec:conv_rel}~-~\ref{sec:Algo} are devoted, respectively, to the applications relating to quantifying the conservativeness of inner approximations and data-driven estimation of Pareto optimal solutions. Numerical illustrations are provided shortly after the key results to quantitatively complement their understanding.  Proofs are furnished in the appendix.

\section{Model description and assumptions}
\label{sec:Modelling Assumptions}

\textbf{Notational Conventions.} Vectors are written in boldface to enable differentiation from
scalars. For any $\mv{a} = (a_1,\ldots,a_d) \in \Real^d$ and 
$\mv{b} = (b_1,\ldots,b_d) \in \Real^d,$  let
$\mv{a}\mv{b} = (a_1b_1,\ldots,a_d b_d)$ and 
$\mv{a}/\mv{b} = (a_1/b_1,\ldots,a_d/b_d)$
denote the respective component-wise operations. Let  $\bar{\R} = \R \cup \{+\infty, -\infty\}$ be the extended real-line,  
$\Real^d_+ = \{\xx \in \Real^d: \xx \geq \mv{0}\}$ denote the positive
orthant, and $\Real^d_{++}$ 
denote its interior. 
For $\xx \in\R^d$ and $A\subset \R^d$, let $d(\xx, A) =\inf\{\Vert \xx - \yy\Vert_\infty:\yy\in A\}$ denote the  distance between $\xx$ and the set $A$. 
For real-valued sequences $\{a_n\}_{n \geq 1}$ and $\{b_n\}_{n \geq 1}$, we write as $n \rightarrow \infty,$  $a_n\sim b_n$ if $\lim_n(a_n/b_n) = 1,$ $a_n = O(b_n)$ if $\limsup_n \vert a_n \vert /\vert b_n \vert < \infty,$ and $a_n = \Omega(b_n)$ if $\liminf_{n} \vert a_n \vert /\vert b_n \vert > 0.$ 
For any Borel measurable set $\mathcal{E} \subseteq \R^d$,  denote the set of all Borel probability measures on $\mathcal{E}$ as $\mathcal{P}(\mathcal{E})$. For any positive integer $K,$ we use $[K]$ to denote $[K] = \{1,\ldots,K\}.$

\subsection{Assumptions on the constraints and illustrative examples} 
\label{sec:unsafe-set} 
Recall the chance constrained optimization model \eqref{eqn:CCP},
which we  labeled as $\CCPalpha$ in the introduction. A key ingredient in the $\CCPalpha$ model is the collection of $K$ constraints $\{ g_k(\xx,\xxi) \leq 0, k \in [K]\},$ where $K$ is a positive integer,  and   
$g_k:\mathbb{R}^m \times \R^d \rightarrow \bar{\R} , k \in [K],$ are lower semicontinuous functions specified suitably for a problem at hand. With $\xxi$ denoting an $\mathbb{R}^d$-valued  random vector modeling all the uncertain factors affecting the decision problem, 
the   constraints $\{ g_k(\xx,\xxi) \leq 0, k \in [K]\}$ typically model critical requirements  such as meeting demands in  commodity supply networks or disaster relief networks (see Examples \ref{eg:transportation}-\ref{eg:nw-design} below), or  regulatory capital requirements in insurance-reinsurance networks (eg., \citealt{blanchet2023efficient}). 
In applications where it is either infeasible or unduly expensive to ensure that the requirements $\{ g_k(\xx,\xxi) \leq 0, k \in [K]\}$ are \textit{always} satisfied, the decision-maker strives to ensure that they are met, at least, with a pre-agreed target probability level $1-\alpha \in (0,1).$ 
In line with this goal, the chance-constrained optimization model $\CCPalpha$ seeks to find a  decision from the set $\mathcal{X} \subseteq \R^m$ that minimizes the cost function $c:\mathcal{X} \rightarrow \R$ while meeting the service level agreement $P\{g_k(\xx,\xxi) \leq 0, k \in [K]\} \geq 1-\alpha.$ 


Throughout the paper, we shall assume that the random vector $\xxi$ has a probability density supported on a closed cone $\mathcal{E} \subseteq \R^d.$ In applications, we typically have $\mathcal{E}$ as either the positive orthant $\R^d_+$ or the euclidean space $\R^d,$ the cost function $c(\xx)$ is  linear in the decision $\xx$, and the service requirements $\{ g_k(\xx,\xxi) \leq 0, k \in [K]\}$ are often specified via functions $g_k(\xx,\xxi)$ which are linear or bilinear in $\xx,\xxi.$ In more sophisticated instances, the functions $\{c,g_i:i \in [K]\}$ may be non-linear (or) may get specified by means of the value of an optimization problem; see,  for eg., \cite{yang2016distributionally,blanchet2023efficient, pena2020dc}. Without restricting to a specific functional form, our assumptions below specify requirements on the constraint functions $\{g_k:k \in [K]\}$ that bestow sufficient regularity 
while allowing broader use. 
We first introduce  the notion of ``safe set" $\mathcal{S}(\xx)$ for any decision $\xx$ which is helpful  towards this end.  For any $\xx \in \R^m,$ let 
\begin{align}
    \mathcal{S}(\xx) = \left\{ \zz \in \mathcal{E}: \max_{k \in [K]} g_k(\xx,\zz) \leq 0 \right\},
    \label{eq:safe-sets}
\end{align}
denote the set of scenarios $\zz$ for which the constraints $g_k(\xx,\zz) \leq 0$ hold for all $k \in [K].$ Observe that $\mathcal{S}: \R^m \rightrightarrows \R^d$ is a set-valued map. Let $\mathcal{X}^\prime := \{\xx \in \mathcal{X}: \text{ interior of } \mathcal{S}(\xx) \neq \emptyset \}$ 

\begin{definition}[\textbf{Non-vacuous set-valued mapping}]
\em 
We call a set-valued map $S: \R^m \rightrightarrows \R^d$ to be \textit{non-vacuous} for the support $\mathcal{E}$ and the decision set $\mathcal{X}$ if for every $\xx \in \mathcal{X},$ the ``unsafe set" 
$\mathcal{E} \setminus {S}(\xx)$ is (i) non-empty, and (ii) bounded away from $\{ \mv{0} \}$ for every $\xx \in \mathcal{X}^\prime.$
    \label{defn:nonvacous-sets}
\end{definition}

\noindent Equipped with these notions, we are now ready to introduce the regularity required for our analysis.

\subsubsection{Homogeneous safe-set model.}
\label{sec:hom-safeset}
We first present the simpler model based on homogeneity here, before moving to the more general assumption in Section \ref{sec:nonhom-safeset}. 
\begin{assumption}
\em 
 There exists a constant $r \neq 0$ such that  for any $\xx \in \mathcal{X}$ and $t > 1,$ we have $t^r\xx \in \mathcal{X}$ and $\mathcal{S}(t^r \xx) = t\mathcal{S}(\xx).$ Further,  the map $\mathcal{S},$ defined in \eqref{eq:safe-sets}, is continuous on $\mathcal{X}$ and 
        non-vacuous  for the support $\mathcal{E}$ and decision set $\mathcal{X}.$
    \label{assume:hom-safe-set}
\end{assumption}

Assumption \ref{assume:hom-safe-set} allows for an easy interpretation  as follows: Given any $\xx \in \R^m$ and $t > 1,$ suppose that a decision-maker wishes to identify a decision $\xx^\prime$ that makes $t\mathcal{S}(\xx)$ scenarios to be safe (that is, $\mathcal{S}(\xx^\prime) = t\mathcal{S}(\xx)$). Satisfaction of Assumption \ref{assume:hom-safe-set} means that they can achieve this by choosing $\xx^\prime = t^r\xx.$ In other words, $\mathcal{S}(\xx^\prime) = \mathcal{S}(t^r\xx) = t\mathcal{S}(\xx).$ As it will become evident from the first main result (Theorem \ref{thm:solution_convergence}) in Section \ref{sec:Asymptotics},  
non-vacuousness of the
safe-sets $\mathcal{S}(\xx)$ will ensure that the value of $\CCPalpha$ in \eqref{eqn:CCP} is not trivially zero or infinite. 

\begin{example}[Joint capacity sizing and probabilistic transportation]\label{eg:routing}
\em 
Consider the  formulation, 
    \begin{align*}
        \min_{\xx,\yy \geq \mv{0}} \sum_{i=1}^M c_i  x_i  + \sum_{i=1}^M\sum_{j=1}^N d_{ij} y_{ij}   \quad \text{s.t.} \ \ \Prob\left\{ \sum_{i:(i,j) \in E} y_{ij}  \geq \xi_j, \  \forall j \in [N] \right\} \geq 1-\alpha,\ \ \sum_{j:(i,j) \in E} y_{ij}  \leq x_i \ \forall i \in [M],
    \end{align*}
which includes the classical transportation problem as an instance with applications in commodity distribution and emergency response; see, eg., \cite{luedtke2010integer,beraldi2004designing}. In this joint capacity sizing and transportation problem,  we have $M$ factories producing a single commodity, $N$ distribution centers, and a distribution center $j \in [N]$ connected with a factory $i \in [M]$ only if $(i,j)$ lies in the edge set $E \subseteq [M] \times [N].$ The goal is to identify  factory supply capacities $\xx = (x_1,\ldots,x_M)$ and a transportation plan $\yy = (y_{ij}: (i,j) \in E)$ jointly such that the cumulative capacity allocation and transportation costs 
is minimized while ensuring  the factories are resourced sufficiently to meet the demand  $\xxi$ at the distribution centers with  probability $1-\alpha.$ 

To verify Assumption \ref{assume:hom-safe-set},
observe that $\mathcal{X} := \{(\xx,\yy) \in \R^M_+ \times \R^{\vert E \vert}_+: \sum_{j: (i,j) \in E} y_{ij} \leq x_i \ \forall i \in [M]\}$ is a cone.  Further, $t \mathcal{S}(\xx,\yy) = \{ t\zz \geq \mv{0} : \sum_{i: (i,j) \in E} y_{ij} 
 \geq z_j \  \forall j \in [N]\} = \{ \zz^\prime \geq \mv{0}: \sum_{i: (i,j) \in E} t y_{ij} 
 \geq z_j^\prime \  \forall j \in [N]\} = \mathcal{S}(t\xx,t\yy).$ Since the graph of the map $\mathcal{S}$ given by 
\[ \left\{(\xx,\yy,\zz) \in \R^M_+ \times \R^{\vert E \vert}_+ \times \R^d_+ : \sum_{i: (i,j) \in E} y_{ij} 
 \geq z_j \  \forall j \in [N], \sum_{j: (i,j) \in E} y_{ij} \leq x_i \ \forall i \in [M] \right\} \] is polyhedral,  the map $\mathcal{S}(\cdot)$ is continuous (\citealt[Eg. 9.35]{rockafellar2009variational}). Further, note that whenever $\mv{0} < \zz^\prime = (z_1^\prime, \ldots,z_d^\prime) \in \mathcal{S}(\xx,\yy)$ for  $(\xx,\yy) \in \mathcal{X},$ we have $(0,\ldots, 0, z_i^\prime, 0, \ldots, 0) \in \mathcal{S}(\xx,\yy),$ for any $i \leq d.$ As a result, $\{\zz \in \R^d_+: \Vert \zz \Vert_\infty \leq \min_{i \leq d} z_i^\prime\} \subseteq \mathcal{S}(\xx,\yy),$ and due to the convexity of $\mathcal{S}(\cdot),$ we have that the map  $\mathcal{S}(\cdot)$ is non-vacuous for the support $\mathcal{E} = \R^d_+$ and decision set $\mathcal{X}.$ \hfill$\Box$ 


\label{eg:transportation}
\end{example}

\noindent 
Note that the reasoning in Eg. \ref{eg:transportation} does not rely  on the network structure, and Assumption \ref{assume:hom-safe-set} can be verified to hold with $r=1$ more broadly for linear chance constraints of the form $\Prob\left\{ A\mv{x} \geq \xxi \right\} \geq 1-\alpha$ featuring  right-hand uncertainty.


\begin{example}[Network design]
\label{eg:nw-design}
\em 
    Consider the problem 
    \begin{align}
        \min_{\xx} \mv{c}_{_V}^\intercal \xx_{_V} + \mv{c}_{_E}^\intercal \xx_{_E}  \quad \text{s.t.} \quad & P\left( \exists \, \mv{f}  : \xx_{_V} - A\mv{f} \geq \xxi, \mv{0} \leq \mv{f} \leq \mv{x}_{_E} \right) \geq 1-\alpha, \label{eq:nw-flow}\\
        &\sum_{i \in V} \mathbf{1}(x_{i,V} > 0) = N, \ \xx_{_V} \geq \mv{0}, \xx_{_E} \geq \mv{0}, \nonumber
    \end{align}
    in which a network designer aims to locate $N$ outposts at suitable vertices of a given network, equip those vertices with supply capacities $\xx_{_V}$ and arcs with flow capacities $\xx_{_E}$ such that the following requirement is met: The stochastic demands arising in different nodes of the network, captured by the random vector $\xxi,$ should be met with a feasible flow with probability at least $1-\alpha.$ Here  $A$ is the node-arc incidence matrix. Network design formulations of this nature arise in disaster relief   \citep{disaster_response} and emergency medical response \citep{boutilier_chan}; see also \cite{atamturk}. 
   If we take $L(\xx,\zz) = \inf\{s \geq 0, \mv{f} : \xx_{_V} - A\mv{f} \geq \zz, s\mv{1} \geq \mv{f}, \mv{0} \leq \mv{f} \leq \mv{x}_{_E}\},$ then note that $L(\xx,\zz)$ is finite if and only if there is a feasible flow. Further, $L(\xx,\xxi) \leq \mv{1}^\intercal \xx_E$ if $L(\xx,\xxi)$ is finite. Therefore, the chance constraint in \eqref{eq:nw-flow} can be equivalently written as $P\{L(\xx,\xxi) - \mv{1}^\intercal \xx_{_E} \leq 0\} \geq 1-\alpha.$ Note that for any $\xx = (\xx_{_V},\xx_{_E}) \in \mathcal{X}= \{(\xx_{V},\xx_{_E}) \geq \mv{0}: \sum_{i \in V} \mathbf{1}(x_{i,V} > 0)= N\}$  and $t > 1,$ we have $t\xx \in \mathcal{X}$ and  $L(t\xx,t\zz)=tL(\xx,\zz).$  An immediate implication of this homogeneity is $\mathcal{S}(t\xx) = t\mathcal{S}(\xx),$ as
   \begin{align*}
   \mathcal{S}(t\xx) &= \{\zz \in \mathcal{E}: L(t\xx,\zz) - t\mathbf{1}^\intercal \xx_{_E} \leq 0\} = \{t\zz^\prime \in \mathcal{E}: L(t\xx,t\zz^\prime) - t\mathbf{1}^\intercal \xx_{_E} \leq 0\}\\
   &= t\{ \zz^\prime \in \mathcal{E}: L(\xx,\zz^\prime) 
 - \mathbf{1}^\intercal \xx_{_E} \leq 0\} = t\mathcal{S}(\xx).
 \end{align*}
   Since the graph of the map $\mathcal{S}(\cdot)$ given by $\{(\xx,\zz) \geq \mv{0}: \xx = (\xx_{_V}, \xx_{_E}), \xx_{_V} - A\mv{f} \geq \zz, \mv{0} \leq \mv{f} \leq \mv{x}_{_E}\}$ is polyhedral, the map $\mathcal{S}(\cdot)$ is continuous (see \citealt[Eg. 9.35]{rockafellar2009variational}). 
   Exactly following the same reason in  Example \ref{eg:transportation}, we have $\mathcal{S}(\cdot)$ to be non-vacuous  for the support $\mathcal{E}$ and the  decision set $\mathcal{X}= \{(\xx_{V},\xx_{_E}) \geq \mv{0}: \sum_{i \in V} \mathbf{1}(x_{i,V} > 0)= N\}.$ Therefore  
   Assumption \ref{assume:hom-safe-set} is satisfied with $r = 1.$   
   \hfill$\Box$
\end{example}

\subsubsection{Non-homogeneous safe-set model.} 
\label{sec:nonhom-safeset} 
We next present an alternative requirement on the constraint functions that holds more broadly than the homogeneous case in Section \ref{sec:hom-safeset}. 

\begin{assumption}
    \label{assume:non-hom-safeset}
    There exist constants $r \neq 0,\, \rho \geq 0,$ and  a  function $g^\ast:\R^m \times \R^d \rightarrow \bar{\R}$ such that the following are satisfied: 
\begin{itemize}
\item[i)] for any $\xx \in \mathcal{X}$ and $t > 1,$ we have $t^r\xx \in \mathcal{X};$
 \item[ii)] for any sequence $(\xx_n,\zz_n) \rightarrow (\xx,\zz),$ 
 \begin{align}
 \liminf_n \frac{\max_{k \in [K]} g_k(n^r\xx_n,n\zz_n)}{n^\rho} \geq g^\ast(\xx,\zz); 
     \label{eq:epi-liminf}
 \end{align}
 \item[iii)] given any sequence $\xx_n \rightarrow \xx$  and $ \zz \in \R^d,$ we have  \begin{align}
     \limsup_n \frac{ \max_{k \in [K]} g_k(n^r\xx_n,n\zz_n)}{n^\rho} \leq g^\ast(\xx,\zz)
     \label{eq:epi-limsup}
 \end{align}  for some sequence $\zz_n \rightarrow \zz;$  
 \item[(iv)] the set-valued  map $\mathcal{S}^\ast: \R^m \rightrightarrows \R^d$ defined by  $\mathcal{S}^\ast(\xx) = \{\zz \in \mathcal{E}: g^\ast(\xx,\zz) \leq 0\}$  is non-vacuous  for the support $\mathcal{E}$ and the decision set $\mathcal{X}.$
\end{itemize}
\end{assumption}

Observe that \eqref{eq:epi-liminf}-\eqref{eq:epi-limsup} are readily satisfied if $g_k(n^r\xx_n, n\zz_n)/n^{\rho} \rightarrow g_k^\ast(\xx,\zz)$
for any sequence $(\xx_n,\zz_n) \rightarrow (\xx,\zz).$ 
The weaker notion of convergence in \eqref{eq:epi-liminf} - \eqref{eq:epi-limsup} is related to the well-known notion of epi-convergence in the optimization literature \cite[Prop. 7.2]{rockafellar2009variational}. Specifically, if we let $g_{k,n}(\xx,\zz) = n^{-\rho} g_k(n^r\xx,n\zz)$ as suitably scaled versions of the constraint functions, then \eqref{eq:epi-liminf} - \eqref{eq:epi-limsup} are equivalent to saying that the sequence of  functions $\max_k g_{k,n}(\xx_n,\cdot),$ are epi-converging to $g^\ast(\xx,\cdot),$ whenever $\xx_n\to \xx$.
Another sufficient condition for \eqref{eq:epi-liminf} - 
 \eqref{eq:epi-limsup} is the epi-convergence of $g_{k,n}(\xx_n, \cdot)$ when  
$g_k(\xx,\cdot)$ is convex for every $\xx$ and $k$ \cite[Prop. 7.48]{rockafellar2009variational}.


\begin{example}[Linear portfolio selection]\label{eg:portfolio_opt}
\em 
Consider the problem, 
\begin{equation*}
       \min_{\xx \in \R^d: \xx \neq \mv{0} } \ c(\xx) \ \quad  \text{ s.t } \quad \ \Prob\left\{ \xx^\intercal \mv{\xi} + \mu_0 (1-\mv{1}^\intercal \xx)  \geq R \right\} \geq 1-\alpha, \ \mv{0} \leq \xx \leq \mv{1},
    \end{equation*}
which aims to select a minimum cost portfolio weight vector $\xx$ whose resulting portfolio return is above a prescribed minimal level  $R > 0$ with probability at least $1-\alpha;$ see, example, \cite{bonami2009exact,pagnoncelli2009computational}. Here $\mu_0$ is the risk-free return and $\xxi$ is the random vector modeling returns of $d$ risky assets. One may take $-c(\xx) = \xx^\intercal E[\xxi] -\mu_0(1-\mv{1}^\intercal \xx)$ signifying average portfolio return, as taken in \cite{pagnoncelli2009computational}, (or) combine it with a  risk measure $c(\xx) = (\xx^\intercal \text{Cov}[\xxi] \xx)^{1/2}$ signifying portfolio variance as in \cite{bonami2009exact}. Note that $g(\xx,\zz) = R - \xx^\intercal \xxi - \mu_0 (1-\mv{1}^\intercal \xx)$ satisfies 
\begin{align*}
    g(n^{-1}\xx_n,n\zz_n) = R - \left(\frac{\xx_n}{n}\right)^\intercal n\zz_n - \mu_0\left( 1- \mv{1}^\intercal \left(\frac{\xx_n}{n}\right)\right) \rightarrow R -\xx^\intercal \zz - \mu_0.
\end{align*}
whenever $(\xx_n,\zz_n) \rightarrow (\xx,\zz).$ As a result, Assumption \ref{assume:non-hom-safeset} holds with $r = -1, \rho = 0,$ and $g^\ast(\xx,\zz) = R - \xx^\intercal \zz - \mu_0.$ The resulting $\mathcal{S}^\ast(\xx) = \{\zz \in \R^d: R - \xx^\intercal \zz - \mu_0 \leq 0\}$ is non-vacuous for the support $\mathcal{E} = \R^d$ and decision set $\mathcal{X} = \{ \xx \in \R^d: \mv{0} \leq \xx \leq \mv{1}, \xx \neq \mv{0}\}$ when  $R < \mu_0.$ \hfill $\Box$
\end{example}

Note that if we had alternatively taken the scaling constants $r,\rho$ in \eqref{eq:epi-liminf} - \eqref{eq:epi-limsup} to be $r = -1,\rho > 0$ in Example \ref{eg:portfolio_opt}, we would obtain $g(n^{-1}\xx_n, n\zz_n)/n^\rho \rightarrow g^\ast(\xx,\zz) = 0$ which would lead to  vacuous $\mathcal{S}^\ast(\xx) = \R^d.$  On the other hand, if we had  taken $r = 1,\rho=2,$ we would obtain $g(n\xx_n, n\zz_n)/n^2 \rightarrow g^\ast(\xx,\zz) = -\xx^\intercal \zz$ which would again lead to a vacuous $\mathcal{S}^\ast(\xx) = \{\zz \in \R^d: -\xx^\intercal \zz \leq 0\}.$ Proposition \ref{prop:unique_r_rho} below asserts that an appropriate choice of scaling constants $r,\rho$ in Assumption \ref{assume:non-hom-safeset} that renders the resulting $\mathcal{S}^\ast(\xx)$ to be non-vacuous is unique. 

\begin{proposition}
    \label{prop:unique_r_rho}
    \em 
    If the collection of constants $(r,\rho)$ for which 
    Conditions (ii) - (iv) of Assumption \ref{assume:non-hom-safeset} are satisfied is non-empty, then it is unique.
\end{proposition}

\subsection{Assumptions on the cost  function}
We next introduce a mild structure  on the cost  $c(\cdot),$ which can be readily verified in the examples in Section \ref{sec:unsafe-set}. To state the assumption, recall the definition 
$\mathcal{X}^\prime := \{\xx \in \mathcal{X}:\,\text{interior of } \mathcal{S}(\xx) \neq \emptyset \}$ and the constant $r$ for which the constraint functions $\{g_k: k \in [K]\}$ satisfy Assumption  \ref{assume:hom-safe-set} or \ref{assume:non-hom-safeset}.

\begin{assumption}
\label{assume:cost}
\em 
The cost function $c:\mathcal{X} \rightarrow \R$ is positively homogeneous, that is, $c(t\xx) = tc(\xx),$ for every $t > 0$ and $\xx \in \mathcal{X}.$ Further, the mapping  $t \mapsto c(t^r\xx)$ is strictly increasing for every $\xx \in \mathcal{X}^\prime.$
\end{assumption}

Note that if $c$ is positively  homogeneous with degree $s > 0,$ then one can simply take $c^{1/s}$ as the cost function satisfying Assumption \ref{assume:cost}. The results also extend to the case of approximately homogeneous costs  where $ c(t^r \xx) \sim  t^r c(\xx),$ as $t \rightarrow \infty,$ uniformly over $\xx$ in compact subsets of $\R^m$ not containing the origin. 
For the ease of exposition, we limit the treatment to Assumption \ref{assume:cost}.

\subsection{Assumptions on the probability distribution of $\xxi$}
\label{sec:Xi_Assumptions}
We assume that the random vector $\xxi$ satisfies either Assumption $(\mathscr  L)$ or Assumption $(\mathscr  H)$ below, with the former corresponding to multivariate light-tailed distributions and the latter capturing multivariate heavy-tailed distributions. For $i \in [d],$ let $f_{\xxi}:\mathcal{E} \rightarrow \R$ denote the probability density of $\xxi,$  $\bar F_i(z) := P(\xi_i >z)$ denote the complementary cumulative distribution function (CDF) of $\xi_i.$  Let $\bar{F}(z) := \max_{i \in [d]} \bar{F}_i(z)$ capture the marginal distribution with the heaviest tail. 
The following definition is useful for introducing Assumptions $(\mathscr{L})$ and $(\mathscr{H}).$

\begin{definition}
\em 
A function $f:\R_+ \rightarrow \R_+$ is
said to be \textit{regularly varying} with index $\rho \in \mathbb{R}$ if $f(tx)/f(t) \rightarrow x^\rho$ for every $x > 0.$
\label{defn:regularly-varying}
\end{definition}
If $f$ is regularly varying with index $\rho,$ we abbreviate this as $f \in \RV(\rho).$ Regularly varying functions offer a systematic and general approach for studying functions with polynomial growth/decay rates, and have served as a fundamental tool for modeling and studying distribution tails: see, example, \cite{feller2008introduction,deHaan, Resnick}.

\begin{assumption*}[$\boldsymbol{\mathscr  L}$] 
\em 
\textbf{[Lighter-tailed distributions]} The  joint probability density $f_{\xxi}$ admits the following representation uniformly over $\zz$ in compact subsets of $\mathcal{E}$ not containing the origin:
\begin{align*}
    -\log f_{\xxi}(t\zz) \sim \varphi^\ast(\zz) {\lambda}(t), \quad \text{ as } t \rightarrow \infty,
\end{align*}
for a positive function $\varphi^\ast: \mathcal{E} \rightarrow \R_{++},$ $\lambda \in \RV(\gamma),$ and $\gamma > 0.$
\end{assumption*}

\begin{assumption*}[$\boldsymbol{\mathscr  H}$] 
  \em
  \textbf{[Heavier-tailed distributions]} 
The  joint probability density $f_{\xxi}$ admits the following representation uniformly over $\zz$ in compact subsets of $\mathcal{E}$ not containing the origin:
\begin{align*}
      f_{\xxi}(t\zz) \sim \varphi^\ast(\zz)t^{-d}\lambda(t), \quad \text{ as } t \rightarrow \infty,
\end{align*}
for a positive function $\varphi^\ast: \mathcal{E} \rightarrow \R_{++},$  $\lambda \in \RV(-\gamma),$ and $\gamma > 0.$
\end{assumption*} 
\noindent 
Lemma \ref{lem:marginals-light}-\ref{lem:marginals-heavy} below characterize  the marginal distributions of $\xxi$ under Assumptions ($\mathscr{L}$) and ($\mathscr{H}$).

\begin{lemma}[\textbf{Marginal distributions under Assumption $\boldsymbol{(\mathscr{L})}$}]
\label{lem:marginals-light}
Under Assumption $(\mathscr L),$  
\begin{itemize}
    \item[a)] the complementary CDFs possess an exponential decay rate as in $\bar{F}_i(tz_i) = \exp\{-c_i z_i^\gamma \lambda(t)[1+o(1)]\}$ as $t \rightarrow \infty,$ where the positive constants $c_i = \inf\{\varphi^\ast(\zz): \zz \in \mathcal{E}, z_i \geq 1\}$ for  $i \in [d];$  
    \item[b)] the heaviest tail $\bar{F}(tz) = \exp\{-\min_{i \in [d]} c_i \, z^\gamma \lambda(t)[1+o(1)]\};$ and  
    \item[b)]  consequently, one could have taken $\lambda(t) = -\log \bar{F}(t)$ in Assumption ($\mathscr L$)  without  loss of generality; in this case we will have $\min_{i \in [d]} c_i = 1.$
\end{itemize}
\end{lemma}

\begin{lemma}[\textbf{Marginal distributions under Assumption $\boldsymbol{(\mathscr{H})}$}]
\label{lem:marginals-heavy}
Under Assumption $(\mathscr H),$ 
\begin{itemize}
    \item[a)] the complementary CDFs possess a polynomial decay rate as in $\bar{F}_i(tz_i) = c_i z_i^{-\gamma} \lambda (t)[1+o(1)]$ as $t \rightarrow \infty,$ where the positive constants $c_i = \int_{z_i > 1} \varphi^\ast(\zz)d\zz$ for  $i \in [d];$   
    \item[b)] the heaviest tail $\bar{F}(tz) = \max_{i \in [d]} c_i \, z^{-\gamma} \lambda(t) [1 + o(1)]; $ and 
    \item[b)]  consequently, one could have taken $\lambda(t) = \bar{F}(t)$ in Assumption ($\mathscr H$)  without  loss of generality; in this case we will have $\max_{i \in [d]} c_i = 1.$ 
\end{itemize}

\end{lemma}



Several  commonly used multivariate distributions, including the Gaussian distributions, multivariate $t$, several classes of elliptical distributions, Archimedean copulas, log-concave distributions and exponential families, satisfy Assumptions $(\mathscr L)$ or $(\mathscr H)$; please refer \cite{deo2023achieving}, Section E.C.2 for a more comprehensive list, their verification, and related properties. Note that the nonparametric nature of $\varphi^*$  allows a great degree of flexibility in modeling various copula and tail dependence structures.
Example \ref{eg:Elliptical_densities} below serves as a pointer towards understanding how Assumptions $(\mathscr L)$ and $(\mathscr H)$ are natural for capturing light and heavy-tailed phenomena respectively. Further examples are available in \cite{Resnick} and references therein.

\begin{example}[Elliptical distributions]\label{eg:Elliptical_densities}
\em 
If $\xxi$ is elliptically distributed on $\R^d,$ then its pdf is given by $f_{\xxi}(\zz) \propto   h\big((\zz - \mv{\mu})^\intercal \Sigma^{-1}(\zz - \mv{\mu})\big),$ for some positive definite matrix $\Sigma,$ mean vector $\mv{\mu} \in \R^d,$ and a suitable generator function $h:\R_+ \rightarrow \R_+$  (see, eg., \cite{frahm2004generalized}, Corollary 4). As special examples, we have the generator $h(t) \propto \exp(-t/2)$ for $\xxi$ to be multivariate normal distributed, and $h(t) \propto (1 + t/\nu)^{-(d+\nu)/2}$ for $\xxi$ to possess multivariate $t$-distribution with $\nu > 0$ degrees of freedom. 

In the multivariate normal case, note that Assumption $(\mathscr L)$ is readily satisfied with $\lambda(t) = t^2$ as, 
\begin{align*}
    -\log f_{\xxi}(t\zz) &= c + (t\zz - \mv{\mu})^\intercal \Sigma^{-1}(t\zz - \mv{\mu})/2 \sim t^2 \zz^\intercal \Sigma^{-1}\zz/2, 
\end{align*}
as $t \rightarrow \infty$, uniformly  in compact sets, for a suitable constant $c.$ 
Similarly, any generator function $h$ satisfying $h(t) = \exp(- c^\prime t^{\gamma/2}(1+o(1)) $ leads to the resulting $-\log f_{\xxi}(t\zz) \sim c^\prime t^{\gamma} (\zz^\intercal \Sigma^{-1}\zz)^{\gamma/2},$ 
 as $t \rightarrow \infty.$ 
 The parameter choice $\gamma < 2$ leads to distributions with tails heavier than the normal distribution and $\gamma < 1$ leads to Weibullian tails that are  heavier than the exponential distribution. 

For multivariate $t$-distributions, the heavy-tailed Assumption $(\mathscr H)$ holds with $\lambda(t) = t^{-\nu}$ as, 
\begin{align*}
    f_{\xxi}(t\zz) \propto \left[1 + \nu^{-1}\left((t\zz - \mv{\mu})^\intercal \Sigma^{-1}(t\zz - \mv{\mu})\right) \right]^{-\frac{(d + \nu)}{2}} \sim \nu^{-1} t^{-(d+\nu)}\left(\zz^\intercal \Sigma^{-1}\zz \right)^{-\frac{(d + \nu)}{2}}, \text{ with } t \rightarrow \infty. \quad \Box 
\end{align*} 

\end{example}

\section{Scaling of  optimal value and solution in the high reliability regime}\label{sec:Asymptotics} 
\subsection{Large deviations characterizations for the probability of constraint violation}
In order to understand the behavior of optimal value and  solutions of $\CCPalpha,$ we first  derive novel characterizations of the probability of constraint violation, $P(\max_{k \in [K]}g_k(\xx,\xxi) > 0),$ under the assumptions introduced in Section \ref{sec:Modelling Assumptions}. To state the results governed by the cases in Assumption \ref{assume:hom-safe-set}-\ref{assume:non-hom-safeset} in a unified manner,   we take throughout the paper,
\begin{align*}
    g^\ast(\xx,\zz) = \max_{k \in [K]} g_k(\xx,\zz)  \quad \text{ if Assumption \ref{assume:hom-safe-set} is satisfied.}
\end{align*}


\begin{proposition}
\label{prop:prob-convergence_LT}
Under Assumption $(\mathscr L)$ and either of Assumptions \ref{assume:hom-safe-set} \underline{or}  \ref{assume:non-hom-safeset}, we have the following convergence,   as $t \rightarrow \infty,$ uniformly over $\xx$ in compact subsets of $\mathcal{X} \setminus \{\mv{0}\}$:
 \begin{align*}
     \log \Prob\left(\max_{k \in [K]} g_k(t^r\xx,\xxi) > 0 \right) \sim -G_0(\xx)\lambda(t),  
 \end{align*}
where $G_0(\xx) = \inf\{\varphi^\ast(\zz): g^\ast(\xx,\zz) \geq 0\}.$
\end{proposition}

\begin{proposition}
\label{prop:prob-convergence_HT}
Under Assumption $(\mathscr H)$ and either of Assumptions \ref{assume:hom-safe-set} \underline{or}  \ref{assume:non-hom-safeset}, we have the following convergence,   as $t \rightarrow \infty,$ uniformly over $\xx$ in compact subsets of $\mathcal{X} \setminus \{\mv{0}\}$:
 \begin{align*}
       \Prob\left(\max_{k \in [K]} g_k(t^r\xx,\xxi) > 0 \right) \sim G_0(\xx) \lambda(t) ,  
 \end{align*}
where $G_0(\xx) =\int_{g^*(\xx,\zz)\geq 0 } \varphi^*(\zz) d\zz.$
\end{proposition}

A well-known example of  explicit characterizations of the constraint violation probability, similar to those in Propositions \ref{prop:prob-convergence_LT} - \ref{prop:prob-convergence_HT}, is from the setting of linear  chance constraint, $g(\xx,\xxi) = v - \xxi^\intercal \xx \leq 0,$ under multivariate normal distribution $\mathcal{N}(\mv{\mu},\Sigma)$ for $\xxi.$ Here take $v$ to be a given constant, $\mv{\mu} \in \R^d,$ and $\Sigma$ is a positive definite covariance matrix.  Since $g(\xx,\xxi) \sim \mathcal{N}(\mv{\mu}^\intercal \xx, \xx^\intercal \Sigma \xx)$ in this instance, we readily have $\Prob(g(\xx,\xxi) >  0)) = \bar{\Phi}\left( (v - \mv{\mu}^\intercal \xx))/\sqrt{\xx^\intercal \Sigma \xx}\right),$ where $\bar{\Phi}(\cdot)$ is the complementary CDF of the standard normal distribution. Due to this explicit expression, the chance constraint $\{\xx: \Prob(g(\xx,\xxi) \leq 0) \geq 1-\alpha \}$ reduces to the deterministic equivalent $\{\xx: v - \mv{\mu}^\intercal \xx + \bar{\Phi}^{-1}(\alpha)\sqrt{\xx^\intercal \Sigma \xx} \leq 0 \}$ in this case. Such deterministic equivalents, while pivotal in the earlier literature, are typically difficult to obtain beyond stylized examples. The probability characterizations in Propositions \ref{prop:prob-convergence_LT} - \ref{prop:prob-convergence_HT} above can be viewed as offering a pathway towards such deterministic characterizations much more broadly, with the caveat that the resulting deterministic constraints set we introduce next becoming equivalent to the chance constrained set only asymptotically.

\subsection{Main result 1: Scaling of the optimal cost and solutions of $\CCPalpha$}
To ease the notational burden, define
\begin{equation}\label{eqn:scaling}
 s_{\alpha} := 
    \bar F^{-1}(\alpha)  
\end{equation}
where $\bar{F}^{-1}$ is the inverse function of $\bar F := \max_{i \in [d]} \bar{F}_i.$ In line with the notation in Lemma \ref{lem:marginals-light}-\ref{lem:marginals-heavy} and  Propositions \ref{prop:prob-convergence_LT} - \ref{prop:prob-convergence_HT}, assign for any $\xx \in \mathcal{X},$ 
\begin{align}
\label{eqn:I_star_limit}
 I(\xx) := 
\begin{cases}
\left[\inf\{\varphi^*(\zz) : g^\ast(\xx,\zz) \geq 0\}\right]^{-1} \min_{i \in [d]} c_i \quad \quad \quad  \ \ & \text{ if } \xxi \text{ satisfies Assumption~(${\mathscr L}$)}\\
[\max_{i \in [d]} c_i]^{-1}\int_{g^\ast(\xx,\zz) \geq 0}  \varphi^*(\zz) d\zz  \quad & \text{ if } \xxi \text{ satisfies Assumption~(${\mathscr H}$)}.  
\end{cases}
 \end{align}

\begin{theorem}
\label{thm:solution_convergence}
   Suppose that the constraint functions satisfy either Assumption \ref{assume:hom-safe-set} \underline{\tt or} Assumption \ref{assume:non-hom-safeset}, the cost function satisfies Assumption \ref{assume:cost}, and $\xxi$ satisfies Assumption ($\mathscr{L}$)  or ($\mathscr{H}$). Then, as $\alpha \rightarrow 0,$
\begin{itemize}
    \item[a)] the optimal value of 
$\CCPalpha$, call it  $v_\alpha^*,$ satisfies
     $v_\alpha^* \sim v^* s_\alpha^{r},$
     where $v^\ast$ is the optimal value of the deterministic optimization problem $\Dapprox$
    below in \eqref{eq:d-approx}:
    \begin{equation}
    \label{eq:d-approx}
    \Dapprox: \qquad\qquad \min_{\xx \in \mathcal{X}} \ c(\xx) \quad \textnormal{ s.t. } \quad   I(\xx) \leq 1. \qquad\qquad\qquad\qquad   
    \end{equation}
    \item[b)] any optimal solution of $\CCPalpha,$ call it $\xx_\alpha^\ast,$ satisfies  $ d(s_\alpha^{-r}\xx_{\alpha}^
    \ast, \mathcal{X}^\ast) \rightarrow 0,$ where 
    $\mathcal{X}^\ast$ denotes the set of optimal solutions for $\Dapprox.$  In particular, if $\Dapprox$ has a unique solution $\xx^\ast,$  then
    \begin{equation}
    \label{eqn:solution_convergence}
         \xx_\alpha^* \sim  \xx^* s_\alpha^{r}, \quad \text{ as } \alpha \to 0.
    \end{equation}
\end{itemize}   
\end{theorem}

Example \ref{eg:Multiple_linear_CCs} below contextualizes the application of Theorem \ref{thm:solution_convergence} to a specific setting. 
\begin{example}
    \label{eg:Multiple_linear_CCs}
    \em 
    Consider the linear chance constraint $P\{ T_k(\xx)\xxi \leq u_k, k \in [K] \} \geq 1 - \alpha,$ with $\xxi$ being elliptically distributed as in Example \ref{eg:Elliptical_densities}, $T_k(\xx) = \xx^\intercal B_k,$ matrices $B_k \in \R^{m \times d},$ and $u_k \in \R,$ for $k \in [K].$ Here the constraint functions $g_k(\xx,\zz) = \xx^\intercal B_k \zz - u_k,$ for $ k \in [K],$  readily satisfy Assumption \ref{assume:hom-safe-set} with $r = -1.$ If the generator  $h(\cdot)$ for the elliptical distribution satisfies $h(t) = \exp(-c^\prime t^{\gamma/2}(1+ o(1)),$ we have from Example \ref{eg:Elliptical_densities} that Assumption $(\mathscr L)$ is satisfied with $\lambda(t) = t^{\gamma}$ and $\varphi^\ast(\zz) = c^\prime (\zz^\intercal \Sigma^{-1} \zz)^{\gamma/2}.$ Let $\sigma_i = \Sigma_{ii}^{1/2}$ for $i \in [d].$ Then from Lemma \ref{lem:marginals-light}(b), we have $s_\alpha \sim \sigma_{\max} [(\ln 1/\alpha)/c^\prime]^{1/\gamma},$ as $\alpha \rightarrow 0,$   
    where $\sigma_{\max} := \max_{i \in [d]} \sigma_i$ and $\alpha \rightarrow 0.$ Further, from the definition of $I(\xx)$ in \eqref{eqn:I_star_limit}, we have  
    $1/I(\xx) = c^\prime\inf\{ (\zz^\intercal \Sigma^{-1} \zz)^{\gamma/2}: \cup_{k \in [K]} \{\xx^\intercal B_k \zz > u_k\}  \} = c^\prime\min_{k \in [K]} \inf \{\zz^\intercal \Sigma^{-1} \zz: \xx^\intercal B_k \zz > u_k \}^{\gamma/2}.$ Since the constraints $\xx^\intercal B_k \zz$ are linear in the variable $\zz,$ a typical application of Lagrange duality leads to $1/I(\xx) = c^\prime \min_{k \in [K]} (u_k^2/\xx^\intercal B_k \Sigma B_k^\intercal  \xx)^{\gamma/2}.$ Then the resulting $\Dapprox$ in Theorem \ref{thm:solution_convergence} is given by, 
    \begin{align}
        \min_{\xx \in \mathcal{X}} c(\xx) \quad \text{s.t.} \quad  \xx^\intercal B_k \Sigma B_k^\intercal \xx \leq c^\prime u_k^2 \quad \forall k \in [K].   
        \label{eq:ccp-apx-eg}
    \end{align}
    Therefore the characterization in Theorem \ref{thm:solution_convergence} translates to  
    \begin{align*}
        v_\alpha^\ast \sim \frac{v^\ast}{\sigma_{\max}}\left(\frac{c^\prime}{\ln 1/\alpha}\right)^{1/\gamma} \quad \text{ and } \quad \xx_{\alpha}^\ast \sim \frac{\xx^\ast}{\sigma_{\max}} \left(\frac{c^\prime}{\ln 1/\alpha}\right)^{1/\gamma},
    \end{align*}
    where  $v^\ast$ and $\xx^\ast$ denote the optimal value and solution of \eqref{eq:ccp-apx-eg} and the target reliability level $1-\alpha \rightarrow 1.$ If the optimal solution set $\mathcal{X}^\ast$ of \eqref{eq:ccp-apx-eg} is not a singleton, then $d(s_\alpha \xx_\alpha^\ast, \mathcal{X}^\ast) \rightarrow 0.$ \hfill$\Box$
\end{example}

\subsubsection{A discussion on the scaling rate.}
First, note that Theorem \ref{thm:solution_convergence} precisely characterizes the rate at which an optimal solution $\xx_\alpha^\ast$ and its respective cost $v_\alpha^\ast = c(\xx_\alpha^\ast)$  vary with respect to the target reliability level $1-\alpha.$ Interestingly, this scaling rate $\xx_\alpha^\ast = v_\alpha^\ast = O(s_\alpha^r)$ is determinable entirely from (i) the constant $r$ in Assumptions \ref{assume:hom-safe-set} or \ref{assume:non-hom-safeset} capturing a certain rate of growth of  the constraint functions $\{g_k: k \in [K]\}$, and (ii) the function $\bar{F}(\cdot)$ which captures the rate at which the complementary CDF of the marginal distributions of $\xxi$ decay to zero. 
In the cases where $r > 0,$ as in Examples \ref{eg:transportation}-\ref{eg:nw-design}, greater reliability is attained by increasing the capacities of supplies (optimal $\xx_\alpha^\ast \rightarrow \infty$ ). On the other hand, in the case where  $r < 0$ as in the portfolio optimization setting, greater reliability is attained by decreasing the weights of the risky assets (optimal $\xx_\alpha^\ast \rightarrow \mv{0}$ ). Further, the dependence on $s_\alpha = \bar{F}^{-1}(\alpha)$ highlights that as the tails of the marginal distribution of $\xxi$
become heavier, the corresponding value of $s_\alpha$ increases. A more precise characterization of the growth rate of $s_\alpha$ under Assumptions $(\mathscr L)$ and $(\mathscr H)$ stipulating light-tailed and heavy-tailed distributions, respectively, are given in Lemma \ref{lem:salpha-growth-rate} below. 
\begin{lemma}
    \label{lem:salpha-growth-rate}
    As $\alpha \rightarrow 0,$ we have
    \begin{itemize}
        \item[(i)] $\log s_\alpha \sim \frac{1}{\gamma} \log\log \frac{1}{\alpha}$ under Assumption $(\mathscr L);$ and 
        \item[(ii)] $\log s_\alpha \sim    \frac{1}{\gamma} \log \frac{1}{\alpha}$ under Assumption $(\mathscr H).$ 
    \end{itemize}
\end{lemma}



Suppose $r > 0$ and consider a reliability level $1-\alpha$ sufficiently close to 1. Then it follows from the scaling rate  characterization  in  Lemma \ref{lem:salpha-growth-rate} that that  a decision-maker endowed with, for instance, twice as much budget as the optimal cost $v_\alpha^\ast$ can reduce the probability of constraint violations at most by a factor $2^{-\gamma/r}$ if $\xxi$ is heavy-tailed. Indeed, this is because   (a) $v^\ast_{\alpha^\prime} = 2v^\ast_\alpha$ if and only if $s_{\alpha^\prime}^r = 2s_\alpha^r$ due to Theorem \ref{thm:solution_convergence}, and (b) $(s_{\alpha^\prime}/s_\alpha)^r = 2$ is solved at $\alpha^\prime = 2^{-\gamma/r}\alpha$ when $s_\alpha \sim \gamma^{-1}\log(1/\alpha).$ The heavier the tail,  smaller is  the index $\gamma,$ resulting in a less significant reduction in constraint violation probability.  If $\xxi$ is light-tailed on the other hand, Lemma \ref{lem:salpha-growth-rate} implies that the probability of constraint violations can be reduced to an exponentially smaller level $\alpha^\prime = \alpha^{2^{\gamma/r}}$ when the decision-maker is allowed to select a decision whose cost is $2v_\alpha^\ast.$  This  can be verified similarly by solving for $\alpha^\prime$ in $(s_{\alpha^\prime}/s_\alpha)^r = 2$  for $s_\alpha \sim \gamma^{-1}\log\log(1/\alpha).$

\begin{num_example}[Scaling of optimal costs]
\label{num-eg:asymptotics}
\em 
To quantitatively illustrate  the scaling in Theorem \ref{thm:solution_convergence}, we consider the joint capacity sizing and transportation problem from  Example \ref{eg:transportation} in a network with $M = 5$ factories, $N = 50$ distribution centers (DCs), and 100 edges. Each DC $j \in [N]$ in the network is connected to a Factory $i \in [M]$ only if $\lceil j/10 \rceil \text{ mod } M = i$ or $\lfloor j/10 \rfloor = i \text{ mod } M;$ the  cost $d_{ij}$ incurred for  transporting unit quantity in the respective edges are set by  sampling uniformly at random from the intervals $[0,1]$ and $[0,0.8].$ The production cost parameters  $(c_i: i \in [M])$  are drawn uniformly from the interval $[0,2].$ Demands $(\xi_j: j \in [N])$ for the commodity at the DCs are independent, with their expected values drawn uniformly from the interval $[1,2].$ 

As in Example \ref{eg:transportation}, we first consider the joint chance-constrained case where the decision-maker 
ensures that all DC demands  are met with probability at least $1-\alpha.$ Figure \ref{fig:asymp_char} illustrates the optimal cost $v_\alpha^\ast$ for different demand distributions, along with  (i)  the respective best-fitting functions with growth rate $s_\alpha^r = (\bar{F}^{-1}(\alpha))^r$ and (ii) the  asymptotic characterization $v^\ast s_\alpha^r$ identified in  Theorem \ref{thm:solution_convergence}.

\begin{figure}[h]
\caption{Target probability level $1-\alpha$  (x-axis, log scale) vs  optimal cost (y-axis). Panels  (a)-(d) correspond to one joint  constraint over 50 DCs. Panels (e)-(h) correspond to 50 individual chance constraints, one per DC}
\vspace{10pt} 
\includegraphics[width=1.02\textwidth]{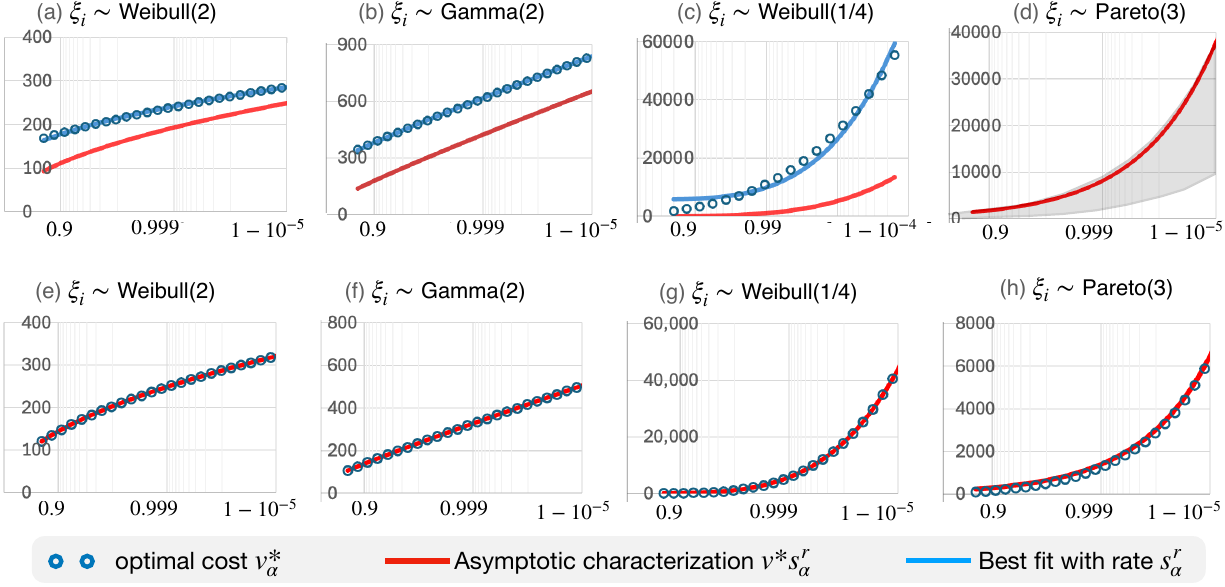}
\label{fig:asymp_char}
\end{figure}
\end{num_example}

The demand distributions in Panels (a)-(d) exhibit increasingly heavier tails from left to right and, correspondingly,  the respective optimal costs are  much larger in magnitude in the Panels (c)-(d). 
For example, if $(\xi_j: j \in [N])$ follows a Weibull distribution with shape parameter $\gamma$, then $s_\alpha \propto \log^{1/\gamma}(1/\alpha)$. Thus, the optimal cost $v_\alpha^\ast$ grows more steeply for smaller $\gamma$. This trend is evident in Panels (a) and (c), where the observed cost growth aligns precisely with the identified rates and Panel (c) exhibiting a significantly faster growth of optimal costs compared to Panel (a). For Gamma-distributed demands with shape parameter 2 (Panel b), we have $s_\alpha \propto \inf\{x \geq 0: (1+x)\exp(-x) = \alpha\} \sim \log(1/\alpha)$ as $\alpha \to 0$. The optimal cost $v_\alpha^\ast$ in Panel (b) reflects this, exhibiting linear growth when plotted on a log-scale. In Panel (d), the demands follow a Pareto distribution with $P(\xi_j > x) = (1+x/s_j)^{-3}$ for scale parameters $s_j > 0$. The corresponding scaling rate is $s_\alpha = \bar{F}^{-1}(\alpha) \propto (1/\alpha)^{1/3} - 1$, which matches the observed cost growth. Since solving the joint chance constraint problem exactly turns out to be intractable for this instance with common solvers, Panel (d) depicts a shaded region between upper and lower bounding individual chance constrained optimal values. Overall, we find that the scaling rate $s_\alpha^r$ and the cost characterization $v^\ast s_\alpha^r$ effectively capture both the magnitude and growth trends under various distributions considered.

Furthermore, the gap between the optimal cost $v_\alpha^\ast$ and its asymptotic counterpart $v^\ast s_\alpha^r$ narrows significantly  when the decision maker is  interested instead in individual chance constraints of the form $P(\sum_{i: (i,j) \in E}y_{i,j} \geq \xi_j) \geq 1-\alpha$ for each DC $j \in [N].$ This suggests that the observed gap in the joint chance constrained setting arises primarily from lower-order terms in the Bonferroni inequality, which are not captured in asymptotic analysis. While this gap may increase with larger $K$ in $\CCPalpha$, the scaling rate $s_\alpha^r$ remains accurate in characterizing the growth of both optimal decisions and costs, regardless of $K$.

\subsubsection{A discussion on the limiting constants $v^\ast,\xx^\ast.$} From Theorem \ref{thm:solution_convergence}, we find that the multiplicative constants $v^\ast,\xx^\ast$ in the relationships $v_\alpha^\ast \sim v^\ast s_\alpha^r$ and $\xx_\alpha^\ast \sim \xx^\ast s_\alpha^r$ are determined in $\Dapprox$ by the joint distribution and the constraint functions $\{g_k: k \in [K] \}$, as informed via the limiting counterparts $\varphi^\ast$ and $g^\ast$ in the assumptions.   
Although explicitly identifying the constraint function $I(\cdot)$ in $\Dapprox$ and solving it to arrive at the limiting constants $v^\ast, \xx^\ast$ 
is not our focus, Proposition \ref{prop:eq-dapprox-LT} below provides a characterization of $\Dapprox$ which can be useful towards that end. 

\begin{proposition}
    \label{prop:eq-dapprox-LT}
    Under Assumption $(\mathscr L)$ and either of Assumptions \ref{assume:hom-safe-set} \underline{or}  \ref{assume:non-hom-safeset}, the limiting deterministic formulation $\Dapprox$ is equivalent to solving,
        \begin{align}
        \min_{\xx \in \mathcal{X}, \, i \in [d]} \ \max_{s \geq 0, \, \zz \in \mathcal{E}} s^{-r/\gamma} c(\xx) \quad \textnormal{ s.t. } \quad \varphi^\ast(\zz) \leq s, \ \  g^\ast(\xx,\zz) \geq 0, \ \Vert \xx \Vert_\infty \leq 1, \ x_i = 1.
        \label{eq:LT-Dapprox-Equiv}
    \end{align}
 \end{proposition}
 The constraint set in \eqref{eq:LT-Dapprox-Equiv} is convex if $g^\ast(\cdot)$ is quasi-concave and $\varphi^\ast$ is convex. The convexity of $\varphi^\ast$ holds, for example, if the probability density $f_{\xxi}$ satisfying the assumption $(\mathscr L)$ is log-concave. Likewise, the constraint set in \eqref{eq:LT-Dapprox-Equiv} can be written as a union of convex sets in the case of joint chance constraints with $g^\ast(\xx,\zz) = \max_{k \in [K]} g_k^\ast(\xx,\zz),$ for quasi-concave functions $g_k^\ast.$  Though not the focus of this paper, Proposition \ref{prop:eq-dapprox-LT} may offer a new window into understanding the eventual convexity properties of chance constraints, which as a standalone topic, has been of great interest in the literature (see, e.g. \cite{van2015eventual,van2019eventual} and references therein).

\subsubsection{Impact of distribution mis-specification.}
Regardless of how challenging it is to precisely compute the constants $v^\ast,\xx^\ast$ for a given setting,  it is worthwhile to note that the asymptotic characterizations of the optimal value $v_\alpha^\ast$ and solution $\xx^\ast$ depend on the target reliability level $1-\alpha$ only via the scaling  $s_\alpha.$ This highlights how the marginal distributions and the copula governing joint distributions decouple in influencing the scaling rate $s_\alpha^r$ and the limiting constants $v^\ast, \xx_\alpha^\ast$ respectively.  
As noted in Corollary \ref{cor:poc} below, this decoupling has an interesting consequence on how mis-specifying the copula of $\xxi$ leads to decisions which are expensive, at most,  by a constant factor, even as the target reliability level $1-\alpha$ is raised to 1. Mis-specifying the rate at which the marginal distributions decay, on the other hand, will distort the rate $s_\alpha = \bar{F}^{-1}(\alpha)$ at which the resulting decisions scale.  
To formally state the former observation, let $\mathcal{P}(F_1,\ldots,F_d)$ denote the collection of all joint distributions of $\xxi = (\xi_1,\ldots,\xi_d)$ for which the  marginal CDF of $\xi_i$ is $F_i$ and the complementary CDF is $\bar{F}_i(x) = 1-F_i(x),$ for $i = 1,\ldots,d.$  
Additionally, let $v_\alpha^\ast(\mathbb{Q})$ denote the optimal value of $\CCPalpha$ obtained by solving  $\CCPalpha$ with a probability measure $\mathbb{Q}$ in place of $\Prob.$
   
\begin{corollary}[\textbf{Impact of copula misspecification}]
    \label{cor:poc}
    Under the assumptions  in Theorem \ref{thm:solution_convergence}, 
    \begin{align*}
       \varlimsup_{\alpha \rightarrow 0} \ \sup_{\mathbb{Q} \in \mathcal{P}(F_1,\ldots,F_d)}  \frac{v_\alpha^\ast(\mathbb{Q})}{v_\alpha^\ast} < \infty. 
    \end{align*}    
\end{corollary}



\section{Application I: Characterizing the costs of distributional robustness}
\label{sec:Robust_Asymototics}
We devote this section towards examining the implications of Theorem \ref{thm:solution_convergence} for  the  well-known Distributionally Robust Optimization (DRO) variants of the  chance-constrained formulation $\CCPalpha.$ 

Given a collection of probability distributions $\mathcal{P}$ defined on $\R^d,$ we consider
the DRO formulation, 
\begin{equation}\label{eqn:DRO_CCP}
    \DROCCPalpha: \qquad   \min_{\xx \in \mathcal X} c(\xx) \quad \text{ s.t. } \quad \inf_{\mathbb Q \in \mathcal P} \  \mathbb{Q} \big\{ g_k(\xx,\xxi) \leq 0, k \in [K] \big\} \geq 1-\alpha,
\end{equation}  
which seeks to identify an optimal decision $\xx$ whose  probability of constraint violation $\{g_k(\xx,\xxi) > 0, \text{ for some } k \in [K]\}$ continues to be smaller than $\alpha$ when evaluated with any  probability distribution $\mathbb{Q}$ in $\mathcal{P}.$ While chance-constrained DRO models of the form $\DROCCPalpha$ date back to \cite{el_ghaoui_robust_var,Erdogan_2005,calafiore2006distributionally}, recent literature has witnessed a surge in their study primarily due to their ability to model and hedge against uncertainty or shifts in the underlying operational environment. 
Commonly used models for the distributional ambiguity set $\mathcal{P}$ in chance-constrained setting include those  specified using moment constraints (see, eg.,  \citealt{el_ghaoui_robust_var, natarajan_robust_var, grani_ccp}), $f$-divergence balls (see,  eg., 
 \citealt{jiang2016data}), or Wasserstein balls (see \citealt{xie2021distributionally,ho2022distributionally,chen2024data}). Please refer the survey article \cite{KUCUKYAVUZ2022100030} for a comprehensive account.


 \subsection{Scaling of costs under $f$-divergence DRO} 
 \label{sec:f-div-DRO}
 First we consider the well-known $f$-divergence based distributional ambiguity set, 
 \begin{equation}
 \label{eqn:f-divergence-ambiguity}
    \mathcal P = \{\mathbb Q: D_{f}(\mathbb Q \Vert \mathbb P) \leq \eta\}, \quad  \text{ where } D_f(\mathbb Q \Vert \mathbb P) = E_{\mathbb  P}\left[f\left(\frac{d\mathbb Q}{d\mathbb P}\right)\right]
\end{equation}
denotes the $f$-divergence between $\mathbb Q$ and $\mathbb P$ and the radius parameter $\eta>0$ captures the extent of distributional ambiguity. 
\begin{assumption}
    \label{assume:f-div}
    The function $f:\R_+ \rightarrow \R_+$ specifying the  $f$-divergence $D_f$ in \eqref{eqn:f-divergence-ambiguity} is continuous, strictly convex, satisfies $\varliminf_{x \rightarrow \infty} f(x)/x = \infty,$ with 
    its minimum $\min_{x \geq 0}f(x) = 0$ attained at $x = 1.$
\end{assumption}
Then, as a consequence of the characterizations in Theorem \ref{thm:solution_convergence} above and the worst-case probability characterization in \cite[Theorem 1]{jiang2016data}, we obtain Theorem \ref{thm:f-div-DRO} below on the scaling of optimal costs and decisions of $\DROCCPalpha.$ Recall $\bar{F}(z) = \max_{i \in [d]} \Prob(\xi_i > z)$ and define 
\begin{align}
  t_\alpha := \bar{F}^{-1}\left(\frac{\alpha}{g(\eta/\alpha)} \right) \text{ where } g(u): = \inf\{x \geq 1: f(x)/x \geq u\} \text{ and } \alpha \in (0,1).
  \label{eq:talpha-defn} 
\end{align}

\begin{theorem}[\textbf{Scaling of  optimal value and solution for $f$-divergence DRO}]
\label{thm:f-div-DRO}
        Suppose that the constraint functions satisfy either Assumption \ref{assume:hom-safe-set} \underline{or} \ref{assume:non-hom-safeset} and Assumptions \ref{assume:cost}-\ref{assume:f-div} hold. Further assume that the  ambiguity set $\mathcal{P}$ in $\DROCCPalpha$  is defined via the $f$-divergence ball in \eqref{eqn:f-divergence-ambiguity}, for some $\eta \in (0,\infty)$ and a distribution  $\Prob$  satisfying either Assumption ($\mathscr{L}$) or ($\mathscr{H}$).   Then, as $\alpha \rightarrow 0,$
\begin{itemize}
    \item[a)] the optimal value of 
$\DROCCPalpha$, call it  $v_\alpha^{\tt f},$ satisfies   
     $v_\alpha^{\tt f} \sim v^* t_\alpha^{r},$
     where $v^\ast$ is the optimal value of the deterministic optimization problem $\Dapprox$ in  \eqref{eq:d-approx}; and 
    \item[b)] any optimal solution of $\DROCCPalpha,$ call it $\xx_\alpha^{\tt f},$ satisfies  $ d(t_\alpha^{-r}\xx_{\alpha}^{\tt f}, \mathcal{X}^\ast) \rightarrow 0,$ where 
    $\mathcal{X}^\ast$ denotes the set of optimal solutions for $\Dapprox.$ If $\Dapprox$ has a unique solution $\xx^\ast,$  then
    \[\xx_\alpha^{\tt f} \sim  \xx^* t_\alpha^r.\] %
\end{itemize}   
\end{theorem}
 \begin{table}[h!]
    \centering
    \caption{Cost $c(\xx_\alpha^{f})$ incurred  by $f$-divergence DRO optimal decisions expressed in terms of baseline scaling  $s_\alpha^r$}
    \begin{tabular}{|c|c|c|c|c|}
    \hline 
        $f$-divergence  &  $f(x)$ & scaling rate $t_\alpha$ & \multicolumn{2}{c|}{asymptotic for $ c(\xx_\alpha^f)$ under}\\[0.5ex]
        & & &  Assump. $(\mathscr{L})$ &  Assump. $(\mathscr{H})$\\[1.0ex]
        \hline \hline 
        KL-divergence & $x\log x$ & $\bar{F}^{-1} \left( \alpha \exp(-\eta/\alpha)\right)$ & $v^\ast \eta^{\frac{r}{\gamma}} \exp[rs_\alpha(1+o(1))]$ &   $v^\ast \exp\left[\frac{r\eta}{\gamma} s_\alpha^\gamma(1+o(1))\right] $\\      [2.0ex]
        $\chi^2$-divergence &  $\frac{1}{2} \vert x - 1 \vert^2$ & $\bar{F}^{-1}\left(\frac{\alpha^2}{2\eta} [1 + o(1)]\right)$ & $2^{\frac{r}{\gamma}} v^\ast  s_\alpha^r$ & $(2\eta)^{\frac{r}{\gamma}}v^\ast s_{\alpha}^{2r + o(1)} $ \\[2.0ex]  
        Polynomial   & $\frac{x^p - p(x-1) - 1}{p(p-1)},$ & $\bar{F}^{-1}\left( \frac{\alpha^{p/(p-1)}[1 + o(1)]}{ (p(p-1)\eta)^{1/(p-1)}} \right)$ & $\left( \frac{p}{p-1}\right)^{\frac{r}{\gamma}} v^\ast s_\alpha^r$ & $c_p v^\ast s_{\alpha}^{\frac{rp[1+o(1)]}{p-1}}$ \\ [-0.5ex]
         divergences & $p > 1$ & & & \\ 
         \hline 
    \end{tabular}
    \label{tab:f-div-scalingrates}
\end{table}

Theorem \ref{thm:f-div-DRO} reveals that while the limiting multiplicative constants $v_\ast$ and $\xx_\ast$ in Theorems \ref{thm:solution_convergence}-\ref{thm:f-div-DRO} are  identical, the scaling rate $t_\alpha$ arising with $f$-divergence DRO can be substantially different from that of the baseline distribution $\mathbb{P}.$ Utilizing the definition of $g(\cdot)$ in \eqref{eq:talpha-defn}, Table \ref{tab:f-div-scalingrates} below  furnishes (i) the scaling rate $t_\alpha$ arising in Theorem \ref{thm:f-div-DRO}, and (ii) an asymptotic for the cost incurred by deploying $f$-divergence DRO optimal decisions. To facilitate comparison, the latter is expressed in terms of the baseline scaling rate $s_\alpha^r$ witnessed in Theorem  \ref{thm:solution_convergence}. The 
constant $c_p = \left[p(p-1)\eta\right]^{\frac{r}{\gamma (p-1)}}$ in Table \ref{tab:f-div-scalingrates}. 

Suppose the parameter $r$ in Assumptions \ref{assume:hom-safe-set} or \ref{assume:non-hom-safeset} is positive. 
Then the first striking observation from Table \ref{tab:f-div-scalingrates} is that the well-known KL-divergence DRO yields decisions which are  exponentially more expensive than those of $\CCPalpha,$ regardless of the radius $\eta.$ To see this, recall that optimal cost $v_\alpha^\ast$ for $\CCPalpha$  grows only at the rate $v_\alpha^\ast = v^\ast s_\alpha^r, $ as $1-\alpha \rightarrow 1$ (see Thm.  \ref{thm:solution_convergence}).  As we shall see in the Numerical Illustration \ref{num-eg:DRO} below, the extreme conservativeness of KL-divergence DRO is observed even at not-so-stringent reliability levels, say, when $1-\alpha$ is  $0.95$ and the radius $\eta$ is taken as small as $0.1.$
In contrast to DRO employing KL-divergences, we learn from Table \ref{tab:f-div-scalingrates} that DRO employing alternative divergence measures such as $\chi^2$ and polynomial divergences are conservative only by a constant factor in light-tailed settings. In the presence of heavy-tailed random variables however,  they yield decisions which are more expensive by a factor  growing polynomially in $s_\alpha$.  

\begin{num_example}[Scaling of DRO decisions]
    \label{num-eg:DRO}
    \em 
    Considering the same transportation problem setting in Numerical Illustration \ref{num-eg:asymptotics}, we empirically explore in Figure \ref{fig:DRO} the impact of $f$-divergence choice on the optimal decisions prescribed by the respective  DRO formulations.
\begin{figure}[h]
\caption{Target probability level $1-\alpha$  (x-axis, log scale) vs  cost of DRO decisions (y-axis). Panels  (a)-(c) correspond to one joint DRO chance constraint over 50 DCs under light-tailed Gamma distributed demands. Panels (d)-(f) correspond to 50 individual DRO chance constraints under heavy-tailed  demands}
\vspace{10pt} 
\includegraphics[width=1.02\textwidth]{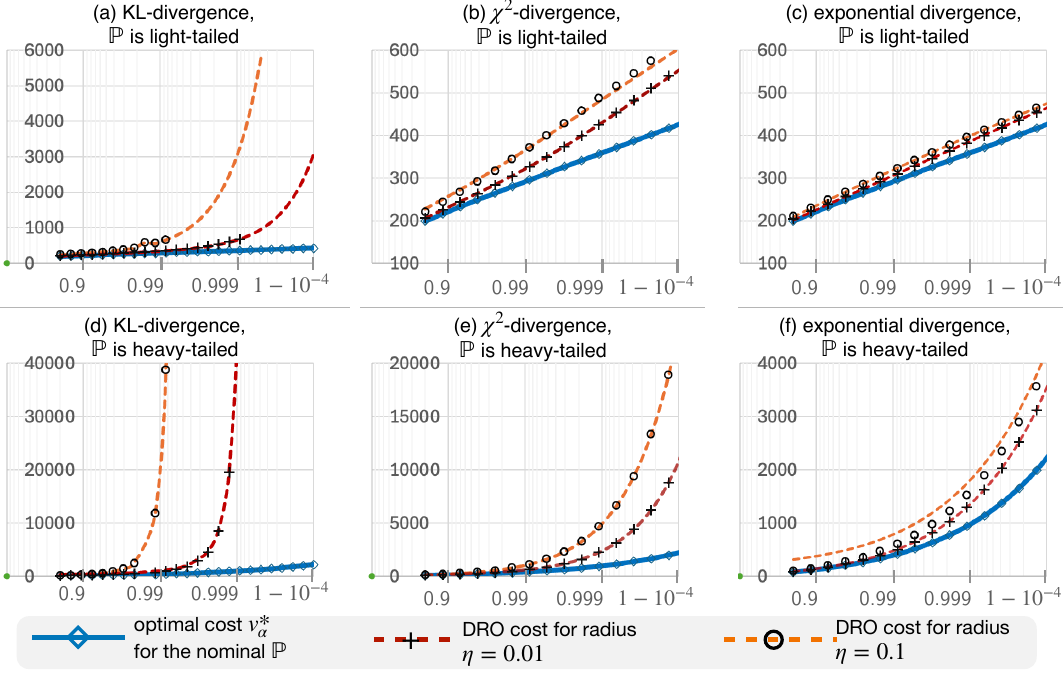}
\label{fig:DRO}
\end{figure}
The markers in Fig \ref{fig:DRO} indicate the actual costs and the dashed lines illustrate the best fitting functions with rate $t_\alpha$ predicted by Theorem \ref{thm:f-div-DRO}. The exponentially high costs brought about by KL-divergence DRO decisions is immediately apparent from Panels (a) and (d): When $\Prob$ is light-tailed, the KL-divegence formulation with radius $\eta = 0.1$ results in decisions which are $33\%$ and $120\%$ more expensive when the reliability targets are set at 95\%  and 120\% respectively (that is, when $1-\alpha = 0.95$ and $0.992$). The same formulation produces decisions which are, respectively, $235\%$ and $9067\%$ more expensive when $\Prob$ is heavy-tailed.
Thus, KL-divergence DRO  can be seen to yield prohivitively expensive decisions even at not so-stringent reliability targets. 

DRO employing $\chi^2$-divergence, on the other hand, produces decisions which, in Panels (b) and (e), grow at  the same scaling rate when $\Prob$ is light-tailed and at a polynomially faster rate when $\Prob$ is heavy-tailed. Besides empirically illustrating the results in Table \ref{tab:f-div-scalingrates}, Fig \ref{fig:DRO} brings out that the extent of conservativeness has a complex, often prohibitive, dependence on the target reliability level $1-\alpha,$ the radius $\eta,$ and the tails of the baseline distribution $\Prob$ under KL and $\chi^2$-divergences. 

In the right-most panels Panels (c) and (f), we take  $f(x) = (x-1)^2\exp(x)$ and examine the effect of exponentially growing $f$ on the cost of optimal decisions. This choice of $f$ is in contrast to the logarthmic growth $f(x) = x\log x$ for KL-divergence and $f(x) = 0.5\vert x -1 \vert^2$ for $\chi^2$-divergence. For such an exponential divergence measure, we find in Panels (c) and (f) that the resulting decisions are conservative within the same orders of magnitude of the optimal costs of $\CCPalpha$ uniformly  for all target reliability levels $1-\alpha.$ In particular, tuning the choice of radius $\eta$ allows one to fulfil its role as a knob  controlling  conservativeness, without letting the target reliability level $1-\alpha$ and  the tails of the baseline distribution $\Prob$ render the resulting decisions unduly expensive. \hfill $\Box$
\end{num_example}

The observations in Table \ref{tab:f-div-scalingrates}, and Numerical Illustration \ref{num-eg:DRO} reinforce the extreme conservativeness of KL-divergence DRO observed in the specific example of quantile estimation in \cite{blanchet2020distributionally,birghila2021distributionally}. Building on 
the empirical observation in  Numerical Illustration \ref{num-eg:DRO}, Proposition \ref{prop:exp-div} below rigorously establishes the desirable scale-preserving nature of $f$-divergences defined via $f(x)$ growing exponentially in $x.$ 
\begin{proposition}
  Suppose that the Assumptions in Theorem \ref{thm:f-div-DRO} hold and the function $f$  defining $f$-divergence in \eqref{eqn:f-divergence-ambiguity} is exponentially growing in $x$ as $x \rightarrow \infty.$ Then as $1-\alpha \rightarrow 1,$ we have $\log c(\xx^{\tt f}_\alpha) \sim \log v_\alpha^\ast$ satisfied under both Assumptions $(\mathscr L)$ and $(\mathscr H)$ for $\Prob.$
  \label{prop:exp-div}
 \end{proposition}

Summarizing, Table \ref{tab:f-div-scalepreserving} below   precisely delineates the $f$-divergence DRO formulations which preserve the scaling from those whose cost grow  much faster than the nominal counterpart $v_\alpha^\ast.$ 
 Please refer \ref{sec:DRO-proofs} in the  appendix for a proof for the properties in Table \ref{tab:f-div-scalepreserving}.

 \begin{table}[h!]
    \centering
    \caption{Delineation of $f$-divergence DRO formulations: 
      Scale-preserving (SP) implies $c(\xx^{\tt f}_\alpha ) = O(v_\alpha^\ast)$ and weakly-scale preserving (w-SP) implies $\log c(\xx^{\tt f}_\alpha) \sim \log v_\alpha^\ast,$ as $1-\alpha \rightarrow 1.$ ``No'' implies $\varliminf_{\alpha \rightarrow 0} \frac{\log c(\xx^{\tt f}_\alpha)}{\log v_\alpha^\ast} > 1$} 
    \begin{tabular}{|c||c|c|c|}
      \hline
    & \multicolumn{3}{c|}{Growth rate of $f(x)/x$ as $x \rightarrow \infty$}\\
      $\Prob$ satisfies  & logarithmically growing & polynomially growing & exponentially growing \\[0.5ex]
      \hline \hline 
       Assumption $(\mathscr L)$ & No & SP & SP \\[0.5ex]
       Assumption $(\mathscr H)$ & No & No & w-SP\\
         \hline 
    \end{tabular}
    \label{tab:f-div-scalepreserving}
\end{table}

\subsection{Distortion of scaling under Wasserstein DRO} Suppose that the  ambiguity set in $\DROCCPalpha$ is given by 
\begin{equation}\label{eqn:Wasserstein_const}
    \mathcal P = \left\{\mathbb Q: d_W(\mathbb P,\mathbb Q)  \leq \eta\right\}, \text{ where } d_W(\mathbb P_1,\mathbb \Prob_2)=\inf_{\Pi(\Prob_1, \Prob_2)} \left[E_{\Pi}|\XX-\YY|^p \right]^{1/p}
\end{equation}
 is the Wasserstein distance of order $p \in [1,\infty)$ between any two probability measures $\mathbb P_1,\mathbb P_2$ supported on $\R^d$ and $\Pi(\Prob_1, \Prob_2)$ is the collection  of all joint distributions with marginals $\XX \sim \mathbb \Prob_1$ and $\YY \sim \mathbb \Prob_2$. Let $\xx_{\alpha}^{\tt W}$ be an optimal solution of the resulting distributionally robust formulation.
\begin{theorem}[\textbf{Conservativeness of Wasserstein DRO}]\label{lem:Wasserstein_CCP}
 Suppose that the constraint functions satisfy either Assumption \ref{assume:hom-safe-set} \underline{or} \ref{assume:non-hom-safeset} and Assumption \ref{assume:cost} holds. Let the  ambiguity set $\mathcal{P}$ in $\DROCCPalpha$  be the Wasserstein ball in \eqref{eqn:Wasserstein_const} for some $\eta \in (0,\infty)$ and  distribution $\Prob$  satisfying either (i) Assumption ($\mathscr{L}$) holds, or (ii) Assumption ($\mathscr{H}$) holds with $p < \gamma.$ Then for any $\varepsilon > 0,$ the cost of the optimal decision prescribed by $\DROCCPalpha,$ denoted by $c(\xx_{\alpha}^{\tt W}),$ satisfies
 \begin{align*}
    c(\xx_{\alpha}^{\tt W}) \geq v\alpha^{-(r/p+\varepsilon)},
 \end{align*}
 for all $\alpha$ sufficiently small and some positive constant $v$ that does not depend on $\alpha.$
 As a result, $\varliminf_{\alpha \rightarrow 0} \frac{\log c(\xx_{\alpha}^{\tt W})}{\log v_{\alpha}^\ast} = \infty$ when $\Prob$ satisfies Assumption $(\mathscr L)$ and
 $\varliminf_{\alpha \rightarrow 0} \frac{\log c(\xx_{\alpha}^{\tt W})}{\log v_{\alpha}^\ast} \geq \gamma/p$ when $\Prob$ satisfies Assumption $(\mathscr H).$ 
\end{theorem}


An immediate implication of Theorem \ref{lem:Wasserstein_CCP} is that the decisions prescribed by Wasserstein DRO formulations are exponentially more expensive than its non-robust counterpart when $\Prob$ is light-tailed and polynomially more expensive when $\Prob$ is heavy-tailed. 
Intuitively speaking, this behavior arises because there always exist probability distributions within the distributional ambiguity set $\mathcal{P}$ whose marginal tail distributions are as heavy as that of a Pareto (power law) distribution with tail parameter $p+\varepsilon$ for any $\varepsilon > 0.$ Thus, regardless of whether the baseline distribution $\Prob$ is light-tailed or heavy-tailed, the cost of decisions prescribed by Wasserstein DRO  grows at a drastically different rate from its non-robust counterpart  as the target service level $1-\alpha$ is raised to 1.

\subsection{Scaling of optimal costs under inclusion of marginal distributions in $\mathcal{P}$}
When equipped with the knowledge of marginal distributions $(F_1,\ldots,F_d)$ of the random vector $\xxi = (\xi_1,\ldots,\xi_d),$ it is natural to define the distributional ambiguity set  via $\mathcal{P}(F_1,\ldots,F_d);$ recall that  $\mathcal{P}(F_1,\ldots,F_d)$ is the collection of all joint  distributions of $\xxi$ under which the CDF  of $\xi_i$ equals $F_i,$ for $i=1,\ldots,d.$  The use of ambiguity sets based on marginal distributions is well-known in the broader DRO literature: see  \cite{rachev_mass_vol1, rachev2006mass, natarajan2021optimization, joint_mixability, pass2024robust} and references therein for  detailed accounts of its tractability and desirable properties in broader operations and risk management applications. 

To describe the scaling of optimal decisions prescribed by marginal distributions based DRO formulations, define
    $r(\xx,\mv{z}) = \inf\left\{t > 0: g^\ast(\xx,t\zz) \geq 0 \right\},$
for any $\zz$ on the unit sphere $\mathcal{S} = \{\zz \in \R^d: \Vert \zz \Vert_\infty = 1 \}.$ Letting $\mathcal{S}_i = \{\zz \in \mathcal{S}: \vert z_i \vert = 1\},$ further define $a_i = \int_{\mathcal{S}_i} r^{-\gamma}(\xx,\zz)d\zz_{-i}.$

\begin{theorem}
[\textbf{Preservation of scaling under marginal distribution based DRO}]
\label{prop:DRO-marg_sol_LT}
        Suppose that the constraint functions satisfy either Assumption \ref{assume:hom-safe-set} \underline{or} \ref{assume:non-hom-safeset} and Assumption \ref{assume:cost} holds. Let the  ambiguity set $\mathcal{P}$ in $\DROCCPalpha$  be  $\mathcal{P}(F_1,\ldots,F_d),$ where  the marginal complementary CDFs $\bar{F}_i(x_i) = 1 - F_i(x_i),$ for $i  =1,\ldots,d,$ satisfy either Assumption ($\mathscr{L}$) or ($\mathscr{H}$). Then, as $\alpha \rightarrow 0,$
\begin{itemize}
    \item[a)] the optimal value of 
$\DROCCPalpha$, call it  $v_\alpha^{\tt M},$ satisfies   
     $v_\alpha^{\tt  M} \sim v^* s_\alpha^{r},$
     where $v^\ast$ is the optimal value of the deterministic optimization problem  $\Dapprox$ in which the constraint function $I(\xx)$ is obtained by setting the specific choice of $\varphi^\ast(\cdot)$ below in the definition of $I$ in \eqref{eqn:I_star_limit}: 
     \begin{align*}
         \varphi^\ast(\zz) = \begin{cases}
             \max_{i \in [d]} c_i \vert z_i \vert^\gamma \quad &\text{if } \xxi \text{ satisfies Assumption } (\mathscr{L})\\ 
             \gamma  r^{-\gamma}(\xx,\zz/\Vert \zz \vert)\sum_{i=1}^d  c_i/a_i \mathbf{1}(\vert z_i \vert = \Vert \zz \Vert) &\text{if } \xxi \text{ satisfies Assumption } (\mathscr{H});
         \end{cases}
     \end{align*}    
    \item[b)] any optimal solution of $\DROCCPalpha,$ call it $\xx_\alpha^{\tt M},$ satisfies  $ d(s_\alpha^{-r}\xx_{\alpha}^{\tt M}, \mathcal{X}^\ast) \rightarrow 0,$ where 
    $\mathcal{X}^\ast$ denotes the set of optimal solutions for the respective $\Dapprox.$  In particular, if $\Dapprox$ has a unique solution $\xx^\ast,$  then $\xx_\alpha^{\tt M} \sim  \xx^* s_\alpha^r,$ as $\alpha \to 0;$ 
\end{itemize}
\end{theorem}

\begin{corollary}
    Consider any distributional ambiguity set $\mathcal{P}$  contained in $\mathcal{P}(F_1,\ldots,F_d),$ where the marginal complementary CDFs $\bar{F}_i(x_i) = 1 - F_i(x_i),$ for $i  =1,\ldots,d,$ satisfy either Assumption ($\mathscr{L}$) or ($\mathscr{H}$). Then the resulting DRO formulation $\DROCCPalpha$ is scale-preserving in the sense  that 
\begin{align*}
     \varlimsup_{\alpha \rightarrow 0} \   \frac{v_\alpha^{\tt DRO}}{v_\alpha^\ast(\mathbb{Q})} < \infty \quad \text{ for every } \mathbb{Q} \in \mathcal{P}.
\end{align*}
\end{corollary}
Thus a modeler seeking to utilize marginal information, when available, in $\DROCCPalpha$ can arrive at a decision that is more expensive only by a constant factor irrespective of the target level $1\ - \alpha.$
While the use of marginals based ambiguity sets  is well-known in the broader DRO literature  (see, eg., \citealt{natarajan2021optimization}), its effectiveness in chance-constrained formulations compared to the more commonly used DRO counterparts 
is a relatively novel observation, albeit intuitive in light of the scaling rate $s_\alpha^r$ depending fundamentally only on marginal distributions. Coincidentally, estimating marginal distributions of $\xxi$ from data does not suffer from curse of dimensions unlike the estimation of  joint distribution, which reinforces the effectiveness of marginals based DRO.  

\subsection{Distortion of scaling under   moments-based DRO} 
Consider a  mean-dispersion based distributional ambiguity set of the form
\begin{equation}\label{eqn:moment_ambiguity}
    \mathcal P = \left\{ \mathbb Q: 
E_{\mathbb Q}[\xxi] = \mv{\mu}, \  
E_{\mathbb Q}\left[\mv{d}(\xxi )\right] \preceq_{\mathcal{D}}  \mv\sigma
 \right\}
\end{equation}
treated in \cite{grani_ccp}, where $\mathcal{P}$ is the set of all Borel probability distributions on $\R^d$ with a given mean vector $\mv{\mu} \in \R^d,$ and $\mv{\sigma}$ is an upper bound on the dispersion measure corresponding to a suitable dispersion function $\mv{d}: \R^d \rightarrow \mathcal{D}.$ Here $\mathcal{D}$ is a proper cone which depends on the dispersion measure chosen by the modeler, and the notation $\mv v \preceq_{\mathcal{D}} \mv w$ means $\mv w - \mv v \in \mathcal{D}.$ The choice $d_{cov}(\zz) = (\zz- \mv{\mu}) (\zz- \mv{\mu})^\intercal$ 
is related to the well-known Chebyshev ambiguity set constraining the mean and covariance of $\xxi:$ Indeed in this case, $\E[d_{cov}(\xxi)] = \E[(\xxi-\mv\mu)(\xxi-\mv\mu)^\intercal)] = \text{Cov}[\xxi]$ and $\mathcal{D}$ is the cone of positive semidefinite matrices. For any $p \in [1,\infty),$ other common choices of the dispersion function include (i) $d_{ad,p}(\zz) = \vert \zz - \mv{\mu}\vert^p = (\vert z_1  - \mu_1 \vert^p, \ldots, \vert z_d  - \mu_d \vert^p ),$ for $p \in (1,\infty),$ which bounds the component-wise mean absolute deviations or its higher moment counterparts, (ii) the choice $d_{sd,p}(\zz) = ((z_i  - \mu_i)_+^p, (\mu_i  - z_i)_+^p): i = 1,\ldots,d)$ which bounds the component-wise upper and lower mean semi-deviations, (iii) $d_{norm,p}(\zz) = \Vert \zz - \mv{\mu} \Vert^p,$ among others. Please refer \cite{el_ghaoui_robust_var, grani_ccp, KUCUKYAVUZ2022100030} and references therein for additional information on choices for dispersion function $d(\cdot).$ 

Theorem \ref{lem:moment_ambiguity} below provides a lower bound of the cost of an optimal decision prescribed by mean-dispersion based DRO formulations. For uniformity in notation, take $p = 2$ if $d = d_{cov}.$
 
\begin{theorem}[\textbf{Conservativeness of moment based DRO}]
\label{lem:moment_ambiguity}
 Suppose that the constraint functions satisfy either Assumption \ref{assume:hom-safe-set} \underline{or} \ref{assume:non-hom-safeset} and Assumption \ref{assume:cost} holds. Let the  collection $\mathcal{P}$ in  $\DROCCPalpha$   be the  mean-dispersion distributional ambiguity set in \eqref{eqn:moment_ambiguity}, and the dispersion function $d(\cdot)$ is one among $d_{cov}, d_{ad,p}, d_{sd,p},$ or $d_{norm,p}$ defined above. Then for any $\varepsilon > 0,$ the cost of the optimal decision $\xx_{\alpha}^{\tt md}$ prescribed by $\DROCCPalpha$  satisfies
 \begin{align*}
    c(\xx_{\alpha}^{\tt md}) \geq v\alpha^{-(r/p+\varepsilon)},
 \end{align*}
 for all $\alpha$ sufficiently small and some positive constant $v$ depending on the problem instance. The conclusion remains the same even if $\mathcal P = \left\{ \mathbb Q: 
E_{\mathbb Q}[\xxi] = \mv{\mu}, \  
E_{\mathbb Q}\left[\mv{d}(\xxi )\right] =  \mv\sigma
 \right\}.$  
\end{theorem}

Thus, among the prominent DRO formulations we have considered in this section, we find only the marginal distribution based DRO to be scale-preserving unconditionally. Excluding this choice, we find $f$-divergence DRO  in which $f(x)/x$ is growing at an exponential rate to be weakly scale-preserving unconditionally. Due to the express availability of tractable reformulations, we find this particular family of $f$-divergence DRO formulations to be  blending tractability with a desirable scale-preserving property which ensures that the optimal decisions produced by the DRO formulation are not unduly conservative. 

\section{Application II: Quantifying \& reducing the conservativeness of safe  approximations}
\label{sec:conv_rel}
Inner approximations to chance constraints, based on Conditional Value at Risk (CVaR) and Bonferroni inequality, have served as   prominent vehicles for tackling the computational challenges in solving  $\CCPalpha$. 
In this section, we demonstrate how the scaling machinery developed in this paper leads to a thoroughly novel approach for quantifying  and reducing the conservativeness of these popular inner  approximations of $\CCPalpha.$ 

\subsection{CVaR-based inner approximations}
\label{sec:cvar-approx}
We begin by recalling the definition of CVaR: For any random variable $Z,$ its CVaR at level $1-\alpha$ is defined as  
   $ \text{CVaR}_{1-\alpha}[Z] =  E[Z\mid Z\geq \text{VaR}_{1-\alpha}(Z)],  \text{ where } \text{VaR}_{1-\alpha}[Z] = \inf\{u: P(Z\geq u)\leq \alpha\}$
is the $(1-\alpha)$-th quantile of $Z.$ 
If a constraint function $g(\xx,\xxi) \leq 0$ is convex in $\xx$, then $\text{CVaR}_{1-\alpha}[g(\xx,\xxi)]$ retains the convexity (\citealt{rockafellar2000optimization}).  In this case, it is well-known that the constraint $\{\xx \in \mathcal{X}: \text{CVaR}_{1-\alpha}[g(\xx,\xxi)] \leq 0\}$ serves as the tightest convex inner approximation to the chance constraint $\Prob\{g(\xx,\xxi) \leq 0\} \geq 1-\alpha;$  see, example, \cite{nemirovski2007convex}. 

For the generic $\CCPalpha$ formulation in \eqref{eqn:CCP}, one may similarly consider the following inner approximation (see \citealt{chen2010cvar}): 
\begin{align}
         \inf_{\xx\in\mathcal X}{c(\xx)} \quad \text{ s.t. }\quad \text{CVaR}_{1-\alpha}\left[ \max_{k \in [K]} \,\eta_k g_k(\xx,\xxi)\right] \leq 0,
        \label{eqn:cvar-approx}
\end{align}
where $\eta_1,\ldots,\eta_K$ are positive constants. Indeed in this case, $\text{VaR}_{1-\alpha}\left[ \max_{k \in [K]} \,\eta_k g_k(\xx,\xxi)\right] \leq 0$ when a decision $\xx$ satisfies the constraints in \eqref{eqn:cvar-approx}, and hence $\Pr\{\max_{k \in [K]} \,\eta_k g_k(\xx,\xxi) \leq 0 \} \geq 1-\alpha.$ This, in turn, is equivalent to the probability constraint in  $\CCPalpha,$ as $\eta_1,\ldots,\eta_K$ are positive constants. Consequently, any solution to \eqref{eqn:cvar-approx} satisfies the probability  constraint in $\CCPalpha$ and hence can be considered to be prescribing ``safe" decisions. Though the constants $\eta_1,\ldots,\eta_K$ may seem superfluous, the flexibility of tuning these parameters has been found to be helpful in improving the quality of subsequent approximations: see \cite{chen2010cvar,Zymler_2011}.  The  variational representation of CVaR due to \cite{rockafellar2002conditional} renders the following to be a computationally attractive way of rewriting the formulation in \eqref{eqn:cvar-approx}:
\begin{align*}
   \inf_{\xx\in\mathcal X,\  u \in \R}{c(\xx)} \qquad \text{ s.t. }\quad u + \frac{1}{\alpha}\E \left[ \left(\max_{k \in [K]} \ \eta_k g_k(\xx,\xxi) - u \right)^+\right] \leq 0.
\end{align*}

As CVaR-based approximations involve computation of expectations, they may become vacuous if the constraint functions are extended real-valued and can take $+\infty$ with positive probability. To avoid this, it is natural to restrict the constraint functions to be real-valued and work with an assumption slightly less general than the epigraphical convergence in Assumption \ref{assume:non-hom-safeset}. 

\begin{assumption}
    \label{assume:constraint-functions-cvar}
\em 
  There exist constants $r \neq 0,\,\rho \geq 0,$ and  a  function $g^\ast:\R^m \times \R^d \rightarrow \bar{\R}$ such that the following are satisfied: 
\begin{itemize}
\item[i)] for any $\xx \in \mathcal{X}$ and $t > 1,$ we have $t^r\xx \in \mathcal{X};$
 \item[ii)] for any sequence $(\xx_n,\zz_n) \rightarrow (\xx,\zz)$ and $k \in [K],$ 
 \begin{align*}
 \lim_{n \rightarrow \infty} \frac{g_k(n^r\xx_n,n\zz_n)}{n^\rho} = g^\ast_k(\xx,\zz); 
 \end{align*}
 \item[iii)] the set-valued  map $\mathcal{S}^\ast: \R^m \rightrightarrows \R^d$ defined by  $\mathcal{S}^\ast(\xx) = \{\zz \in \mathcal{E}: g^\ast(\xx,\zz) \leq 0\}$  is non-vacuous  for the support $\mathcal{E}$ and the decision set $\mathcal{X}.$
\end{itemize}  
\end{assumption}

\begin{theorem}[\textbf{Scaling of the optimal value and solution of CVaR approximation}]
\label{thm:convex_rel}
    Suppose Assumption~\ref{assume:cost}, Assumption~\ref{assume:constraint-functions-cvar},  and either Assumption ($\mathscr{L}$) or ($\mathscr{H}$) hold. 
    Let $v_\alpha^{\tt apx}$ and $\xx_\alpha^{\tt apx}$ denote, respectively, the optimal value and an optimal solution to the CVaR approximation in \eqref{eqn:cvar-approx}. Further assume that $\E\left[ \max_{k \in [K]} g_k(\xx,\xxi)\right] < \infty.$ Then, as $\alpha \rightarrow 0,$
\begin{itemize}
    \item[a)] the optimal value $v_\alpha^{\tt apx}$ satisfies
     $v_\alpha^{\tt apx} \sim cv^* s_\alpha^{r},$
     where $c \in [1,\infty)$ and $v^\ast$ is the optimal value of the deterministic optimization problem $\Dapprox$
     in \eqref{eq:d-approx};    
    \item[b)] any optimal solution $\xx_\alpha^{\tt apx}$ satisfies  $ d(s_\alpha^{-r}\xx_{\alpha}^{\tt apx}, c\,\mathcal{X}^\ast) \rightarrow 0,$ where 
    $\mathcal{X}^\ast$ denotes the set of optimal solutions for $\Dapprox.$  In particular, if $\Dapprox$ has a unique solution $\xx^\ast,$  then $\xx_\alpha^{\tt apx} \sim  c\,\xx^* s_\alpha^{r}.$
    
    \item[c)] the conservativeness of the CVaR approximation relative to $\CCPalpha,$ for all values of $\alpha$ sufficiently small, is thus essentially determined by the constant $c$ which is given by, 
    \begin{align*}
        c = 
        \begin{cases}
            1 \quad&\text{ under Assumption } (\mathscr{L}),\\
            (1-1/\gamma)^{-r}  &\text{ under Assumption } (\mathscr{H}). 
        \end{cases}
    \end{align*}
\end{itemize}   
\end{theorem}

Theorem \ref{thm:convex_rel} reveals that the extent of conservativeness introduced by CVaR-based  approximations remains bounded, irrespective of the 
target reliability level $1-\alpha.$ Notably, the conservativeness  is determined fundamentally by the heaviness of the tails of the distribution of $\xxi$ and is larger for heavier-tailed random variables. In particular, the gap in approximation is inversely proportional to $(1-1/\gamma)^r,$ when the distribution of $\xxi$  has polynomially decaying heavy-tails, such as,  $P(\xi_i > z) = O(z^{-\gamma})$ as $\vert z \vert \rightarrow \infty.$  Conversely, for light-tailed distributions meeting Assumption $(\mathscr L)$, two key observations emerge:
(i) the conservativeness vanishes, and (ii) the solutions $\xx_{\alpha}^{\tt apx}$ of the CVaR approximation and the target solution $\xx_\alpha^\ast$ of the $\CCPalpha$  formulation coincide asymptotically. 

Under both Assumptions $(\mathscr L)$ and $(\mathscr H),$ we  have $\xx^{\tt apx}_{\alpha}/\xx^\ast_\alpha \rightarrow c$ componentwise, where $\xx^\ast_\alpha$ is  optimal for the original $\CCPalpha$ formulation \eqref{eqn:CCP}.  Refer the Panels (c) and (f) in Figure \ref{fig:cvar-line-search} for an empirical demonstration of this finding. 
As we shall see in Section \ref{sec:line-search}, this property can be used to develop an elementary algorithm  capable of improving $\xx^{\tt apx}_\alpha$ into a
nearly optimal solution for $\CCPalpha.$ 


\subsection{Inner approximations based on Bonferroni inequality}
\label{sec:bfi-approx}
As joint constraints are computationally more challenging than individual chance constraints, Bonferroni inequality has served as an intuitive and popular approach towards breaking up multiple constraints, such as in \eqref{eqn:CCP}, into individual probability constraints as below: For any positive constants $\eta_1,\ldots,\eta_d$ satisfying $\sum_{i=1}^d \eta_i = 1,$  consider
\begin{align}
   \qquad \inf_{\xx \in \mathcal{X}}\,c(\xx)  \quad \text{ s. t. } \quad \Prob\left\{ g_k(\xx,\xxi) \leq 0 \right\} \geq 1-\eta_k \alpha, \quad \forall k \in [K]. 
     \label{eqn:bfi-approx}
\end{align}
Indeed, for any $\xx$ satisfying the $K$ individual probability constraints in \eqref{eqn:bfi-approx}, we have $\Prob(\cup_{k \in [K]}\{ g_k(\xx,\xxi) > 0\}) \leq \sum_{k=1}^K \Prob(g_k(\xx,\xxi) > 0) \leq \sum_{k=1}^K \eta_k \alpha = \alpha,$ due to the Bonferroni's inequality. Consequently, any solution to \eqref{eqn:bfi-approx} in turn satisfies the probability  constraint in $\CCPalpha.$ 

\begin{theorem}[\textbf{Scaling for Bonferroni inequality based approximation}]
\label{thm:bfi-approx}
    Suppose Assumption~\ref{assume:cost}, Assumption~\ref{assume:constraint-functions-cvar},  and either Assumption ($\mathscr{L}$) or ($\mathscr{H}$) hold. 
    Let $v_\alpha^{\tt apx}$ and $\xx_\alpha^{\tt apx}$ denote, respectively, the optimal value and an optimal solution to the Bonferroni inequality based approximation in \eqref{eqn:bfi-approx}. Then the following hold as $\alpha \rightarrow 0$:
\begin{itemize}
    \item[a)] Under Assumption $(\mathscr{L}),$ the optimal value $v_\alpha^{\tt apx}$ and any optimal solution $\xx_\alpha^{\tt apx}$ satisfy
     $v_\alpha^{\tt apx} \sim {v}^\ast s_\alpha^{r}$ and $d(s_\alpha^{-r}\xx_{\alpha}^{\tt apx}, {\mathcal{X}}^\ast) \rightarrow 0,$ where ${v}^\ast$ and $\mathcal{X}^\ast$ are the optimal value and optimal solutions of the deterministic optimization problem $\Dapprox$ in \eqref{eq:d-approx};

    \item[b)] Under Assumption $(\mathscr{H}),$ the optimal value $v_\alpha^{\tt apx}$ and any optimal solution $\xx_\alpha^{\tt apx}$ satisfy
     $v_\alpha^{\tt apx} \sim \tilde{v} s_\alpha^{r}$ and $d(s_\alpha^{-r}\xx_{\alpha}^{\tt apx}, \tilde{\mathcal{X}}) \rightarrow 0,$ where $\tilde{v}$ and $\tilde{X}$ are the optimal value and optimal solutions of 
    \begin{align}
        \min_{\xx \in \mathcal{X}} \ c(\xx) \quad \text{ s.t. } \quad I_k(\xx) \leq \eta_k \quad \forall k \in [K], 
        \label{eqn:bfi-dapprox}
    \end{align}   
    with $I_k(\xx) = \int_{g_k^\ast(\xx,\zz) \geq 0} \varphi^\ast(\zz) d\zz.$
    \end{itemize}   
\end{theorem}

Theorem \ref{thm:bfi-approx} reveals that the inner approximation based on Bonferroni inequality features the same scaling as $\CCPalpha$. Hence the extent of conservativeness introduced by the formulation \eqref{eqn:bfi-approx}
remains bounded even as the problem becomes more challenging with the  target  level $1-\alpha$ approaching 1. Similar to CVaR approximation, the conservativeness of \eqref{eqn:bfi-approx} vanishes if the distribution of $\xxi$ is light-tailed and the resulting solution $\xx_{\alpha}^{\tt apx}$ coincides with a target solution $\xx_\alpha^\ast$ of the $\CCPalpha$  formulation  asymptotically. 
However, in the case of heavy-tailed distributions,  the limiting deterministic optimization problem in \eqref{eqn:bfi-dapprox} does not  match with \eqref{eq:d-approx}. In particular, note that the  constraint sets $\cap_{k \in [K]} \{\zz: \int_{g_k^\ast \geq 0} \varphi^\ast(\zz) d\zz \leq \eta_k \}$ and $\{\zz:\int_{g^\ast  \geq 0} \varphi^\ast(\zz) d\zz \leq 1\}$ are of very different nature even when we set $\eta_1 = \cdots = \eta_K.$ Hence one cannot expect the solution $\xx_{\alpha}^{\tt apx}$ to be proximal to the scaled solution set $cs_\alpha^r \mathcal{X}^\ast,$ which was a nice feature of CVaR-based inner approximations that cannot be guaranteed with Bonferroni inequality based inner approximations.

\subsection{A line-search procedure to obtain  asymptotically optimal solutions for $\CCPalpha$}
\label{sec:line-search}
Consider a target reliability level $1-\alpha \in (0,1)$ and a solution $\xx^{\tt apx}_\alpha$ obtained by solving either the CVaR-based inner approximation in \eqref{eqn:cvar-approx} or the Bonferroni inequality based inner approximation in \eqref{eqn:bfi-approx}. Assuming the ability to evaluate the probability of constraint violation $p(\xx) = P(\max_{k \in [K]} \, g_k(\xx,\xxi) > 0)$ at any $\xx \in \mathcal{X},$ Algorithm  \ref{algo:Vanish_Regret} below outlines an  elementary  line search procedure which returns a strictly improved feasible solution to $\CCPalpha.$ This strict improvement is  guaranteed by the  observations in Lemma \ref{lem:line-search-properties} below.
\begin{lemma}
    \label{lem:line-search-properties}
    Suppose $\mathcal{X}$ is a cone and $p(\xx)$ is  continuous.  
    Then under the assumptions of Thm. \ref{thm:convex_rel}, 
    \begin{itemize}
        \item[a)] the cost of the decisions parameterized by $\xx(t) = t^r\xx_\alpha^{\tt apx}$ is strictly increasing in $t > 0;$ 
        \item[b)] the decision  $\xx(t)$ is feasible for $\CCPalpha$ for any $t \geq t_\alpha,$ where $t_\alpha = \inf\{t \in [0,1]: p(t^r\xx_\alpha^{\tt apx}) \leq \alpha\};$
        \item[c)] we have $t_\alpha < 1$ and consequently $c(\xx^\prime_\alpha) <  c(\xx^{\tt apx})$ for the assignment $\xx^\prime_\alpha = t_\alpha^r \xx^{\tt apx}_\alpha.$
    \end{itemize}
\end{lemma}

\begin{algorithm}[h!] 
  \caption{Line search for improving an inner approximation's solution to be asymptotically optimal for $\CCPalpha$}
  \ \vspace{-2pt}\\
  \KwIn{Target reliability level $1-\alpha$} 
     \noindent \textbf{1. Solve an inner approximation to $\CCPalpha$}:  If $\xxi$ satisfies Assumption $(\mathscr L),$ obtain $\xx_\alpha^{\tt apx}$ by either solving \eqref{eqn:cvar-approx} or \eqref{eqn:bfi-approx}.
     If $\xxi$ satisfies Assumption $(\mathscr H),$ obtain $\xx_\alpha^{\tt apx}$ by solving \eqref{eqn:cvar-approx}.\\
     
  \noindent \textbf{2. Find the best feasible solution along the ray:}  Solve the one-dimensional line search  $t_\alpha = \inf\{t \in [0,1]: p(t^r\xx_\alpha^{\tt apx}) \leq \alpha\}.$

\noindent \textbf{3. Return} 
 $\xx_\alpha^\prime = t_\alpha^r \xx_\alpha^{\tt apx}$ as our candidate solution for solving $\CCPalpha.$ 
  \label{algo:Vanish_Regret}  
\end{algorithm}

Despite the simplicity of the line-search in Algorithm \ref{algo:Vanish_Regret},  Proposition \ref{prop:Approx_alg_guarantees} below demonstrates that the improved solution $\xx_\alpha^\prime$ returned by Algorithm \ref{algo:Vanish_Regret} is asymptotically optimal for $\CCPalpha$.
In Proposition \ref{prop:Approx_alg_guarantees} below,  recall that $v_\alpha^\ast$ denotes the optimal value of $\CCPalpha$ in \eqref{eqn:CCP}. 
\begin{proposition}
\label{prop:Approx_alg_guarantees}
Suppose the Assumptions in Theorem \ref{thm:convex_rel} and Lemma \ref{lem:line-search-properties} hold. Then the solution $\xx_\alpha^\prime = t_\alpha^r\xx_\alpha^{\tt apx}$ returned by Algorithm \ref{algo:Vanish_Regret} has vanishing relative optimality gap: Specifically, 
\begin{align*}
   \lim_{1-\alpha \rightarrow 1} \  \frac{c(\xx_\alpha^\prime) - v_\alpha^\ast}{v_\alpha^\ast} = 0, 
   \end{align*}
and Algorithm \ref{algo:Vanish_Regret} decreases the cost by the factor $\lim_{1-\alpha \rightarrow 1} c(\xx_\alpha^\prime)/c(\xx_\alpha^{\tt apx}) =  c.$
\end{proposition}

\begin{num_example}[Effectiveness of Algorithm \ref{algo:Vanish_Regret}]
\em 
Considering the same transportation example in Numerical Illustration \ref{num-eg:asymptotics}, we compare in Figures \ref{fig:cvar-line-search}(a) and \ref{fig:cvar-line-search}(d) the following optimal costs at various target levels $1-\alpha$: (i) the cost $c(\xx_\alpha^\ast)$ due to the optimal solution of $\CCPalpha,$ (ii) the cost $c(\xx_\alpha^{\tt apx})$ due to the CVaR inner approximation \eqref{eqn:cvar-approx} obtained by setting $\eta_k = 1$ for $k \in [K],$ and (iii) the cost $c(\xx_\alpha^\prime)$ due to the solution returned by the line search in Algorithm \ref{algo:Vanish_Regret}. 
\begin{figure}[h]
\caption{Effectiveness of Algorithm \ref{algo:Vanish_Regret} for target reliability levels $1-\alpha$ in $x$-axis of Panels (a)-(b), (d)-(e)}
\label{fig:cvar-line-search}
\begin{center}
\includegraphics[width=\textwidth]{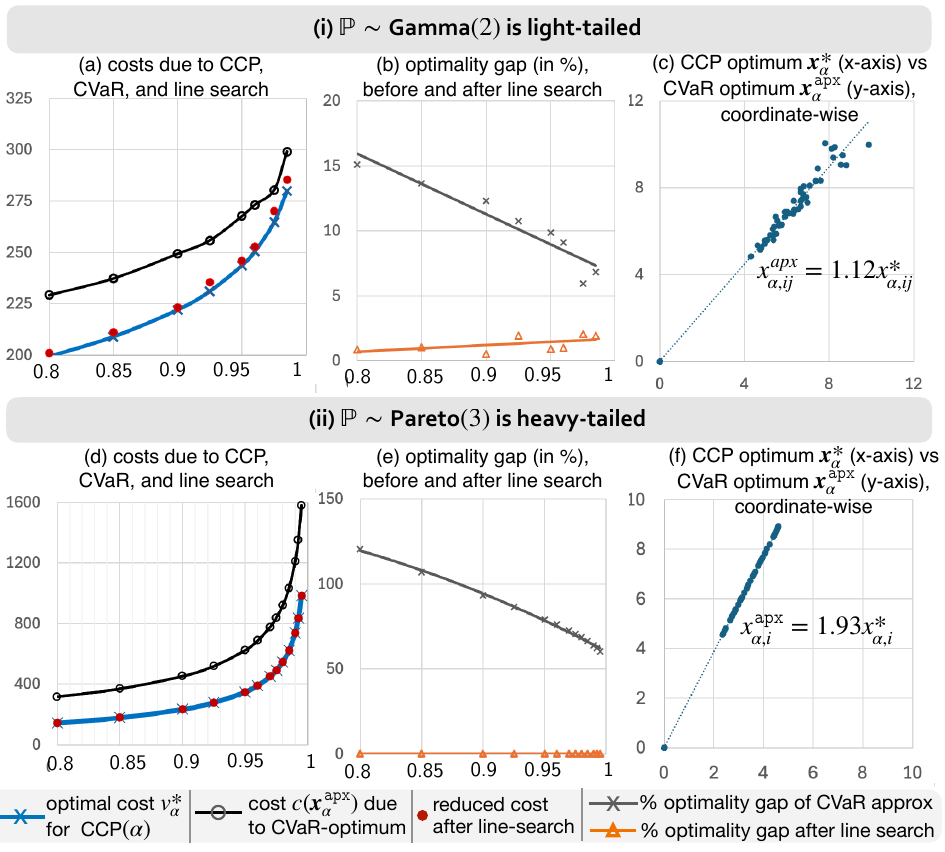}     
\end{center}
\end{figure}
Besides illustrating the near-optimality of $\xx_\alpha^\prime$ pictorially in Panels (a) and (d), we find from Panels (b) and (e) that the optimality gap gets reduced by as much as 50 to 120 percentage points  in the heavy-tailed setting and 5 - 15 percentage points in the light-tailed setting. Panels (c) and (f)  illustrate the reason behind the effectiveness of the line search: As predicted by Theorem \ref{thm:convex_rel}, we see that $\xx_\alpha^{\tt apx}/\xx_{\alpha}^\ast \approx c$ component-wise, which is brought about vividly by the scatter plots between the  components of the solutions output by $\CCPalpha$ and its respective CVaR inner approximation. \hfill$\Box$
\label{num-eg:cvar-line-search}
\end{num_example}


At a conceptual level, the near-optimality of the line search in Algorithm \ref{algo:Vanish_Regret} is primarily due to the  observation that $\xx^{\tt apx}_\alpha/\xx^\ast_\alpha \rightarrow c$ componentwise, where $\alpha \rightarrow 0$ and $\xx^\ast_\alpha$ is a solution to $\CCPalpha$ in \eqref{eqn:CCP}. As a result, conducting a line search for a less-expensive feasible solution along the ray passing through $\xx_\alpha^{\tt apx}$ turns out to be sufficient for obtaining a solution with vanishing relative optimality gap even though $\xx_\alpha^{\tt apx}$ may be sub-optimal by a constant factor. 

\section{Application III: Estimation of Pareto-efficient solutions from limited data}
\label{sec:Algo}
In this section, we use the  scaling in Theorem \ref{thm:solution_convergence}  as the basis to develop  a novel semiparametric method for estimating approximately Pareto-optimal  decisions from limited data. 

\subsection{The $\Omega(1/N)$  estimation barrier in handling  constraint violation probabilities}
\label{sec:dd-underestimate}
In practice, the underlying probability distribution $\Prob$ with which the problem of interest, $\CCPalpha$ in \eqref{eqn:CCP}, should be solved is seldom known and is  informed only by independent realizations $\xxi_1,\ldots,\xxi_N$ from $\Prob.$ In this case, one of the most basic approaches towards solving $\CCPalpha$ is to approximate the probability of constraint violation in \eqref{eqn:CCP} via its average computed from the $N$ samples as below: 
\begin{align}
    p_{_N}(\xx) = \frac{1}{N}\sum_{i=1}^N \mathbb{I}\big\{ g_k(\xx, \xxi_i) > 0 \text{ for some } k = 1,\ldots, K  \big\}.
    \label{eq:saa-prob}
\end{align}
See example, \cite{ahmed2008solving,pagnoncelli2009sample} for the considerations pertinent to approximating $\CCPalpha$ with this Sample Average Approximation (SAA) approach. Recalling the notation $p(\xx) = \Pr( \max_{k \in [K]} g_k(\xx, \xxi_i) > 0),$ we have the following  due to Chebyshev's inequality:
\begin{align}
    \Pr \big( \vert p_N(\xx) - p(\xx) \vert  > \varepsilon p(\xx) \big) \leq \frac{1-p(\xx)}{N \varepsilon^2 p(\xx)}. 
    \label{eq:chebyshev}
\end{align}
Hence, we need at least $N > \alpha^{-1} (1-\alpha) \varepsilon^{-2} \delta^{-1}$ samples in order to well-approximate the constraint violation probability at a boundary point $\xx \in \mathcal{X}$ satisfying $p(\xx) = \alpha,$ with at least $(1-\delta) \times 100\%$ confidence. Similar sample requirement arises with the well-known scenario approximation approach as well, see \cite{calafiore2006scenario}. 
Consequently, both SAA and scenario approximation approaches are well-suited for solving $\CCPalpha$ only when the number of available samples $N$  and the considered  constraint violation probability level $\alpha^{(N)}$  are such that 
$N \times \alpha^{(N)}$ is not small. 
Conversely, if $\alpha^{(N)} \notin \Omega(1/N)$ as $N \rightarrow \infty,$ constraint violations are rarely observed in-sample, leading to hazardous underestimation of the  probability of constraint violation. 

\begin{num_example}
\em 
Here we consider the joint capacity sizing and transportation example in Numerical Illustration \ref{num-eg:asymptotics} with $30$ distribution centers. Taking $N = 1000$ independent observations of the demands realized at the 30 distribution centers as our dataset, we first illustrate the properties of the SAA optimal solutions $\hat{\xx}_\alpha$ obtained by solving  $\min_{\xx \in \mathcal{X}}\{c(\xx): p_N(\xx) \leq \alpha\}$ at different target reliability levels $1-\alpha \in [0.8, 1-10^{-5}].$ The out-of-sample reliability levels $1-p(\hat{\xx}_\alpha)$ met by the  SAA solutions for the joint chance-constrained formulation and  their respective costs $c(\hat{\xx}_\alpha)$ are presented for the case of independent light-tailed Gamma(2) distributed demands in Figure \ref{fig:dd_underestimate_LT} and Table \ref{tab:dd-underestimate} 
below. Corresponding results for the individual chance-constrained heavy-tailed  counterparts are presented in Figure  \ref{fig:dd_underestimate} in the Introduction and Table \ref{tab:dd-underestimate}. The  costs and constraint violation probabilities are reported together with their  respective 95\% confidence intervals computed from  32 independent experiments. 

From Panels (a) of both  Figures \ref{fig:dd_underestimate} and \ref{fig:dd_underestimate_LT}, we  observe that the SAA optimal solutions fall significantly short when the target reliability level $1-\alpha$ is such that $N \times \alpha < 10.$ The shortfall is severely more pronounced when $N \times \alpha \ll 1,$ as typically zero samples fall in the constraint violation regions for the identified SAA optimal decisions.  Insufficient in-sample constraint violations often leads the SAA to incorrectly conclude a decision $\xx$ whose $p(\xx) \in (\alpha,1/N)$ to be feasible, even though it  patently violates the requirement $p(\xx) \leq \alpha$. As a result, we see from  Panels (a)-(b) of  Figures \ref{fig:dd_underestimate} and \ref{fig:dd_underestimate_LT} that the cost $c(\hat{\xx}_\alpha)$ and  reliability $1-p(\hat{\xx}_\alpha)$ offered by SAA solutions taper and do not improve for  $\alpha \ll 1/N.$ This example serves to illustrate why SAA cannot be relied upon to yield a decision $\xx$ satisfying $p(\xx) \ll 1/N.$  
\hfill$\Box$
\begin{figure}[h]  \includegraphics[width=0.95\textwidth]{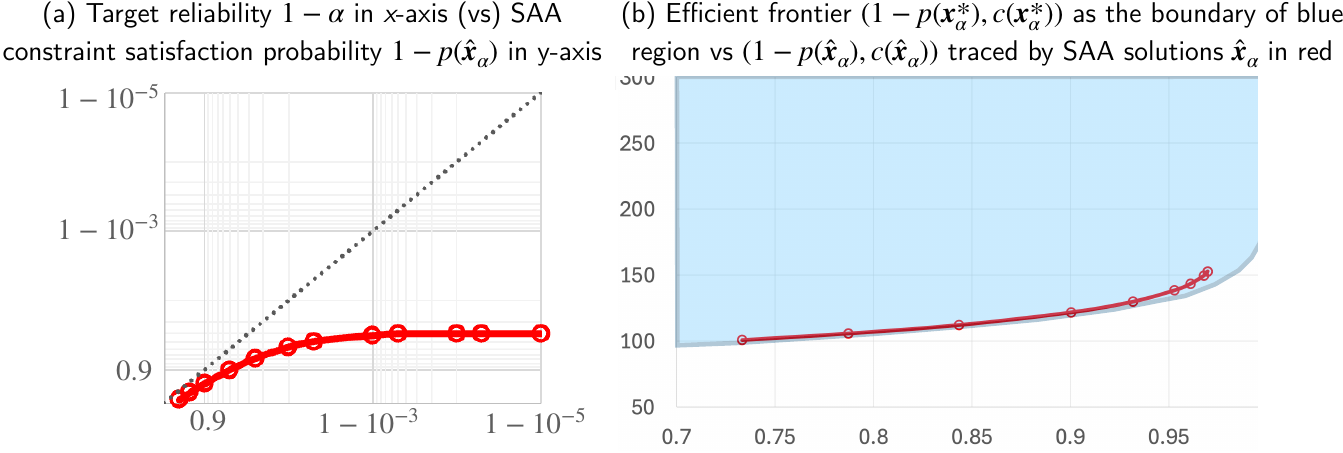}
  \caption{\textnormal{Performance of SAA optimal solutions $\hat{\xx}_\alpha$ obtained at target reliability levels $1-\alpha \in [0.8, 1-10^{-5}]$  using $N = 1000$ samples
    from independent Gamma(2) distributed demands}}
\label{fig:dd_underestimate_LT}
\end{figure}
\begin{table}[h!]
    \centering
    \caption{Constraint satisfaction probabilities and the costs of SAA optimal solutions $\hat{\xx}_\alpha$}
    \begin{tabular}{|c|c|c|c|c|}
    \hline 
        & \multicolumn{2}{c|}{$\Prob$ is light-tailed} &    \multicolumn{2}{c|}{$\Prob$ is heavy-tailed}\\ 
    \hline 
      Target $1-\alpha$   & $1-p(\hat{\xx}_\alpha)$ & $c(\hat{\xx}_\alpha)$ & 
          $1-p(\hat{\xx}_\alpha)$ & $c(\hat{\xx}_\alpha)$\\
    \hline \hline 
    0.8 & 0.733 $\pm$ 0.005 & 100.720 $\pm$ 0.422 & 0.772 $\pm$ 0.003 & 44.184 $\pm$ 0.170 \\
    0.85 & 0.787 $\pm$ 0.004 & 105.687 $\pm$ 0.484 & 0.826 $\pm$	0.002 & 54.923 $\pm$ 0.172 \\
    0.9 & 0.843 $\pm$ 0.004 & 112.072 $\pm$ 0.591 & 0.878 $\pm$ 0.002 & 71.831 $\pm$ 0.241 \\
    0.95 & 0.900 $\pm$ 0.004 & 121.650 $\pm$ 0.815 & 0.933 $\pm$ 0.002 & 106.432 $\pm$ 0.451\\ 
    0.975 & 0.932 $\pm$ 0.003 & 129.775 $\pm$ 1.023 & 0.963 $\pm$ 0.001 & 149.790 $\pm$ 0.970 \\
    0.99 & 0.953 $\pm$ 0.003 & 138.345 $\pm$ 1.307 & 0.981 $\pm$ 0.001 & 223.545 $\pm$ 1.529\\
     $1-5 \times 10^{-3}$ & 0.961 $\pm$ 0.002  & 143.199 $\pm$ 1.456  & 0.988 $\pm$ 0.0007 & 294.528 $\pm$ 3.710 \\    
   $1 - 10^{-3}$ & 0.968 $\pm$ 0.002  & 149.483 $\pm$  1.601
    & 0.994 $\pm$ 0.0005 & 503.952 $\pm$ 15.540 \\
    $1-5 \times 10^{-4}$ & 0.970 $\pm$ 0.002  & 152.612 $\pm$ 1.685  & 0.996 $\pm$ 0.0005 & 648.958 $\pm$ 34.310 \\    
    $1-10^{-4}$ & 0.970 $\pm$ 0.002  & 152.612 $\pm$ 1.685  & 0.995 $\pm$ 0.0005 & 764.963 $\pm$ 53.252 \\ 
        $1-5 \times 10^{-5}$ & 0.970 $\pm$ 0.002  & 152.612 $\pm$ 1.685  & 0.996 $\pm$ 0.0005 & 779.464 $\pm$ 55.678 \\    
    $1 - 10^{-5}$ & 0.970 $\pm$ 0.002  & 152.612 $\pm$ 1.685  & 0.996 $\pm$ 0.0005 & 791.064 $\pm$ 57.624\\
    \hline 
    \end{tabular}
    \label{tab:dd-underestimate}
\end{table}

\label{num-eg:DD-Pareto-Front}
\end{num_example}

The statistical bottleneck demonstrated in Numerical Illustration \ref{num-eg:DD-Pareto-Front} is acute in applications requiring very high levels of reliability, such as power distribution and telecommunication networks (where $\alpha \ll 10^{-4})$, regulatory risk management (where $\alpha \leq 1/100$), etc. Due to the low tolerance for constraint violation (specified via small $\alpha$) in these settings, obtaining sufficient number of representative samples to ensure $p(\hat{\xx}_\alpha) < \alpha$ becomes often an impossible proposition unless additional parametric distributional assumptions are made.   
The need for effective data-driven algorithms that perform well even when the  target $\alpha \ll 1/N$ has been widely acknowledged, particularly in the example of  quantile estimation.  Note that estimating quantile of a random variable $\xi$ at a level $1-\alpha$ can, in fact, be seen as one of the simplest instances of $\CCPalpha$ via its definition  
\begin{align}
q_{\alpha} = \min\{x: \mathbb{P}(\xi \leq x) \geq 1-\alpha\}.    
\label{eq:quantiles}
\end{align}
The field of Extreme Value Theory (EVT) in Statistics (see, eg., \citealt{deHaan,Resnickbook}) has bridged this critical need in the particular context of  estimation of quantiles $q_\alpha,$ when $\alpha N $ is small, corresponding to  levels far beyond those captured via the empirical distribution  from the  data samples. 

\subsection{Basis for extrapolation of solution trajectories inspired from EVT}
\label{sec:evt-analogy}
The central ingredient in extreme value theory (EVT) which enables estimation of quantile $q_\alpha,$ at levels $\alpha \ll 1/N,$ is a limiting characterization of the form, 
\begin{align}
    \lim_{\alpha_0 \rightarrow 0} \frac{q_{\alpha_0/t} - q_{\alpha_0}}{s_{\alpha_0}} = \frac{t^{\gamma} - 1}{\gamma}, 
    \label{eq:quantile-extrapolation}
\end{align}
akin to \eqref{eqn:solution_convergence}, which allows the extrapolation of extreme tail quantile $q_{\alpha}$  from an intermediate quantile $q_{\alpha_0}$  which is observed sufficiently in the available data; see \cite[Thm. 1.1.6]{deHaan}. Here $\gamma$ is a positive constant,  $s_\alpha$ is a  scaling function depending on the underlying distribution, and one takes $\alpha = \alpha_0/t$ to facilitate an extrapolation via \eqref{eq:quantile-extrapolation}. The enduring appeal of this extrapolation, revolving around the limiting  relationship \eqref{eq:quantile-extrapolation}, stems from the fact it does not require the modeler to make parametric assumptions about distribution tail regions which are not sufficiently witnessed in data. Its practical effectiveness and minimally assuming  semiparametric nature have made this approach a cornerstone for estimating quantiles  in diverse engineering and scientific disciplines.  

Momentarily viewing our chance-constrained formulation $v_\alpha^\ast = \inf\{c(\xx): P\big( \max_{k \in [K]} g_k(\xx,\xxi) \leq 0\big) \geq 1-\alpha\}$ as a sophisticated generalization of the quantile estimation task  in \eqref{eq:quantiles}, we propose to utilize the following consequence of Theorem \ref{thm:solution_convergence} in a similar way the limiting relationship \eqref{eq:quantile-extrapolation} has powered quantile extrapolation in EVT. In particular, Proposition \ref{prop:solution-set-convergence} below serves as the counterpart of \eqref{eq:quantile-extrapolation} that allows direct data-driven estimation of solutions at constraint violation probability levels  far below $1/N$ without requiring the modeler to make parametric assumptions on the distribution of $\xxi.$ To state Proposition \ref{prop:solution-set-convergence}, let $\mathcal{X}^\ast(\alpha) =\arg\min \{c(\xx): P\big( \max_{k \in [K]} g_k(\xx,\xxi) \leq 0\big) \geq 1-\alpha\} $ denote the optimal solution set of $\CCPalpha$ at any $\alpha \in (0,1).$

\begin{proposition}
    Suppose the assumptions in Theorem \ref{thm:solution_convergence} are satisfied.  Suppose also that $\mathcal X^*$ is  singleton.
     Then for any $t \geq 1,$ we obtain 
    \begin{align*}
    &s_{\alpha_0}^{-r} d\left( t^{\frac{r}{\gamma}}\xx_{\alpha_0}^*, \  \mathcal{X}^\ast \left(\alpha_0^t \right) \,\right) \ \rightarrow \ 0 \quad \text{ under Assumption } (\mathscr L); and \\
   &s_{\alpha_0}^{-r} d\left( t^{\frac{r}{\gamma}}\xx_{\alpha_0}^*, \  \mathcal{X}^\ast \left(\frac{\alpha_0}{t} \right) \,\right) \ \rightarrow \ 0 \quad \text{ under Assumption } (\mathscr H), 
    \end{align*}
    as $\alpha_0\to 0$.
    \label{prop:solution-set-convergence}
\end{proposition}

Suppose, for a moment, we have access to an optimal solution  $\xx_{\alpha_0}^\ast$  obtained by solving ${\tt CCP}(\alpha_0)$ at an appropriate base constraint violation probability level $\alpha_0.$ Then 
Proposition \ref{prop:solution-set-convergence} brings out a desirable property of the extrapolated trajectory of decisions,    $\xx(t) =         t^{\frac{r}{\gamma}}\xx_{\alpha_0}^\ast$ for $t \geq 1,$ constructed from the optimal solution $\xx_{\alpha_0}^\ast$ at the base level. In particular,  Proposition \ref{prop:solution-set-convergence} asserts that such an extrapolated trajectory $({\xx}(t): t \geq 1)$ gets vanishingly close, in a relative sense, point-wise to the trajectory of optimal solution sets $(\mathcal{X}^\ast(\alpha): \alpha = \alpha_0^t, \, t \geq 1)$ in the case of light-tailed distributions and $(\mathcal{X}^\ast(\alpha): \alpha = \alpha_0/t, t \geq 1)$ in the case of heavy-tailed distributions. Thus, even when $\alpha = \alpha_0^t$  is such that $\alpha N$ is small and constraint violations at that level are  observed  insufficiently in the dataset $\{ \xxi_1,\ldots,\xxi_N\}$ as a consequence, the extrapolated trajectory $(\xx(t): t \geq 1)$ and  Proposition \ref{prop:solution-set-convergence} can be seen as empowering us with a pathway for approximating the solutions of $\CCPalpha$ to which have no access otherwise in the statistically challenging setting where $\alpha N \ll 1.$.


\subsection{Estimation of weakly Pareto-efficient solutions which breach the $\Omega(1/N)$  barrier}
Given a dataset $\{\xxi_1,\ldots,\xxi_N\}$ comprising $N$ i.i.d. samples of $\xxi,$ consider the trajectory,
\begin{align}
    \bar{\xx}^{(N)} = \left( t^{r} \hat{\xx}_{\alpha_0}^{(N)}: \ t \geq 1 \right),
    \label{eq:extrapolated-trajectory}
\end{align}
where $\hat{\xx}_{\alpha_0}^{(N)}$ is an optimal solution estimate  obtained either from  the sample average approximation (or) the scenario approximation of ${\tt CCP}(\alpha_0)$ constructed using the dataset $\{\xxi_1,\ldots,\xxi_N\}.$ 
For the base SAA solution $\hat{\xx}_{\alpha_0}^{(N)}$ to be accurate, it is necessary from the application of Chebyshev's inequality in \eqref{eq:chebyshev} that  $\alpha_0 =   \Omega(N^{-1}\varepsilon^{-2} \delta^{-1}),$ as $N \rightarrow \infty.$  With $\bar{\xx}^{(N)}$ defined in \eqref{eq:extrapolated-trajectory} serving as the key construct enabling estimation of approximately Pareto optimal decisions, we first empirically bring out its  properties in Numerical Illustration \ref{num-eg:DD-Pareto-Front} below before proceeding to its theoretical properties in Theorem 
 \ref{thm:WPE-data-driven} and the following Section \ref{sec:P-Model}.


\noindent 
\textbf{Numerical Illustration \ref{num-eg:DD-Pareto-Front} continued.} Taking the base constraint violation probability level $\alpha_0 = 0.2,$ Figure \ref{fig:dd_pareto_fronts} and Table \ref{tab:dd-Pareto-Fronts} below present the trajectory of (out-of-sample constraint satisfaction probability $1-p(t^r\hat{\xx}_{\alpha_0}^{(N)}),$ cost $c(t^r\hat{\xx}_{\alpha_0}^{(N)})$) pairings traced by the  decisions $\bar{\xx}^{(N)} = (t^r\hat{\xx}_{\alpha_0}^{(N)}: t \geq 1).$
Unlike the tapering of constraint satisfaction probabilities and costs observed with SAA in Table \ref{tab:dd-underestimate} and Figure \ref{fig:dd_underestimate}, we observe from Figure \ref{fig:dd_pareto_fronts} that the extrapolated trajectory $\bar{\xx}^{(N)}$ is able to produce decisions whose constraint satisfaction probabilities can be significantly larger than the empirically observed barriers of $0.97$ in the light-tailed joint chance-constrained setting and $0.995$ in the heavy-tailed individual chance-constrained setting. In particular, even when equipped with only $N = 1000$ samples, the trajectory  $\bar{\xx}^{(N)}$ is able to produce decisions $t^r\hat{\xx}_{\alpha_0}^{(N)}$ whose constraint  violation probabilities $p(t^r\hat{\xx}_{\alpha_0}^{(N)})$ can be as small as $10^{-5}$ or even smaller when larger values of $t$ are used. We also observe that the associated confidence intervals in Table \ref{tab:dd-Pareto-Fronts} to be significantly narrower than those reported for SAA in Table \ref{tab:dd-underestimate}. More interestingly in Figure \ref{fig:dd_pareto_fronts}, the extrapolated decisions in $\bar{\xx}^{(N)}$ offer reliability, cost pairings $(1-p(t^r\hat{\xx}_{\alpha_0}^{(N)}),c(t^r\hat{\xx}_{\alpha_0}^{(N)}))$ that lie close to the efficient frontier capturing the best possible (reliability, cost) pairings attainable for the problem. \hfill$\Box$

\begin{figure}[h]
\caption{\textnormal{(reliability $1-p(t^r\hat{\xx}_{\alpha_0}^{(N)}),$ cost $c(t^r\hat{\xx}_{\alpha_0}^{(N)})$) pairings offered by the  trajectory of decisions $\bar{\xx}^{(N)}$ depicted in black, to be compared with   $(1-p(\hat{\xx}_\alpha), c(\hat{\xx}_\alpha))$  offered by SAA optimal solutions (in red, data from Table \ref{tab:dd-underestimate}). All attainable $((1-p(\xx), c(\xx)): \xx \in \mathcal{X})$ pairings are shaded in blue, with its boundary representing the efficient frontier}}
\label{fig:dd_pareto_fronts}
\includegraphics[width=\textwidth]{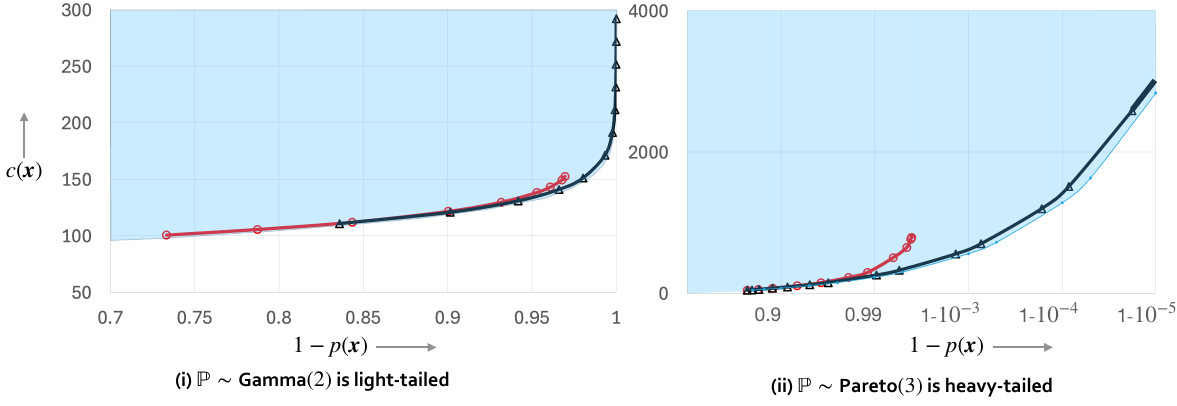}
\end{figure}

\begin{table}[h!]
    \centering
    \caption{Constraint satisfaction probabilities and the costs offered by decisions in the trajectory $\bar{\xx}^{(N)}$}
    \begin{tabular}{|c|c|c|c|c|c|}
    \hline 
         \multicolumn{3}{|c|}{$\Prob$ is light-tailed} &    \multicolumn{3}{|c|}{$\Prob$ is heavy-tailed}\\ 
    \hline 
       $t$   & $1-p(t^r\hat{\xx}_{\alpha_0}^{(N)})$ & $c(t^r\hat{\xx}_{\alpha_0}^{(N)})$ & $t$  & 
          $1-p(t^r\hat{\xx}_{\alpha_0}^{(N)})$ & $c(t^r\hat{\xx}_{\alpha_0}^{(N)})$\\
    \hline \hline 
    1.0 & 0.733 $\pm$ 0.005 & 100.720 $\pm$ 0.422 & 1.0 & 0.772 $\pm$ 0.003 & 44.184 $\pm$ 0.170 \\
    1.1 & 0.836 $\pm$ 0.003 & 110.792 $\pm$ 0.465 & 1.1 & 0.797 $\pm$	0.003 & 48.631 $\pm$ 0.189 \\
    1.2 & 0.902 $\pm$ 0.002 & 120.864 $\pm$ 0.506 & 1.3 & 0.829 $\pm$ 0.002 & 55.669 $\pm$ 0.217 \\
    1.3 & 0.942 $\pm$ 0.002 & 130.936 $\pm$ 0.549 & 1.6 & 0.877 $\pm$ 0.002 & 70.139 $\pm$ 0.273\\ 
    1.4 & 0.966 $\pm$ 0.001 & 141.008 $\pm$ 0.591 & 2.0 & 0.915 $\pm$ 0.001 & 88.369 $\pm$ 0.344 \\
    1.5 & 0.980 $\pm$ 0.0006 & 151.080 $\pm$ 0.634 & 2.7 & 0.951 $\pm$ 0.001 & 119.936 $\pm$ 0.467\\
   1.7 & 0.993 $\pm$ 0.0003  & 171.224 $\pm$ 0.718 & 3.4
    & 0.969 $\pm$ 0.0006 & 151.109 $\pm$ 0.589 \\
    1.9 & 0.998 $\pm$ 0.0001  & 191.368 $\pm$ 0.803 & 5.9 & 0.991 $\pm$ 0.0002 & 258.393 $\pm$ 1.006 \\ 
    2.1 & $1- 8 \times 10^{-4} \pm$ 0.00004  & 211.512 $\pm$ 0.887 &  7.4 & 0.995 $\pm$ 0.0001 & 325.555 $\pm$ 1.268\\
   2.3 & $1- 3 \times 10^{-4} \pm$ 0.00002  & 231.656 $\pm$ 0.972 &  12.6 & 0.999 $\pm$ 0.00003 & 556.691 $\pm$ 2.168\\
   2.5 &  $1- 9 \times 10^{-5} \pm$ 0.000007  & 251.800 $\pm$ 1.056 &  15.9 &  $1-7\times 10^{-4} \pm$ 0.00002 & 701.387 $\pm$ 2.732\\
   2.7 & $1- 3 \times 10^{-5} \pm$ 0.000003  & 271.944 $\pm$ 1.141 &   27.1 &  $1-2\times 10^{-4}\pm$ 0.000005 & 1199.355 $\pm$ 4.671\\
    2.9 & $1- 10^{-5} \pm$ 0.000001  & 292.088 $\pm$ 1.225 & 34.2  &  $1-9 \times 10^{-5} \pm$ 0.000002 & 1511.092 $\pm$ 5.886\\

    \hline 
    \end{tabular}
    \label{tab:dd-Pareto-Fronts}
\end{table} 

Numerical Illustration \ref{num-eg:DD-Pareto-Front} empirically demonstrates the ability of the  trajectory $\bar{\xx}^{(N)}$ to produce approximately Pareto optimal decisions whose reliability levels $1-p(t^r\hat{\xx}_{\alpha}^{(N)})$  can be significantly smaller than the $1/N$ barrier. This empirically observation can be quite surprising when viewed  from the following standpoint: 
If $\alpha \ll 1/N,$ note that it is  statistically impossible to even just verify the hypothesis that $p(\xx) \leq \alpha$,  when one is, say,  presented with a   decision  $\xx \in \mathcal{X}$ satisfying $p(\xx) \in (\alpha/2,2\alpha).$ 
Theorem \ref{thm:WPE-data-driven} and Corollary \ref{cor:P-Model-data-driven} below examine  the  above  intriguing  empirical observation on approximate Pareto optimality  and argue that it is not a coincidence: In particular, Theorem \ref{thm:WPE-data-driven} below and Corollary \ref{cor:P-Model-data-driven} and the following section   establish  the weak Pareto efficiency of the solutions traced by the trajectory $\bar{\xx}^{(N)}$ and its optimality for the well-known $P$-model (see \citealt{Charnes_Cooper_PModel,He_Sim_Zhang}). 

\begin{definition}[\textbf{Weak Pareto efficiency}]
\em 
    A sequence of decisions $(\xx_n: n \geq 1) \subset \mathcal{X}$ is weakly Pareto efficient in balancing the cost and constraint violation probability if, for any given $\varepsilon > 0,$ there exists $\delta \in (0,1)$ and $n_0 \geq 1$ such that the following holds for all $n \geq n_0$: 
    \begin{align}
        \big\{ \xx \in \mathcal{X} :  \ c(\xx) < c(\hat{\xx}_n) - \varepsilon \vert c({\xx}_n)\vert \big\} \subseteq \big\{\xx \in \mathcal{X}:\  \log p(\xx) \geq (1-\delta) \log p({\xx}_{n}) \big\}.
        \label{eq:WPE}
    \end{align}
    \label{defn:weak-pareto-efficiency}
\end{definition}
In words, \eqref{eq:WPE} implies that for every decision $\xx \in \mathcal{X}$ whose cost is smaller than $c(\xx_n)$  as in $c(\xx) < c({\xx}_n) - \varepsilon \vert c({\xx}_n)\vert,$ it is necessary that the respective constraint violation probability $p(\xx)$ is larger by 
        $\log p(\xx) \geq (1-\delta) \log p({\xx}_{n}).$
In Theorem \ref{thm:WPE-data-driven} below, we  require $\max_{k \in [K]}g_k$ to be Caratheodory, besides the mild structural assumptions we have been requiring on the constraint functions $(g_k: k \in [K])$ and the distribution of $\xxi.$   A function $F:\Omega\times \mathcal E \to \R$ is said to be Caratheodory, if for almost every $\xx\in\mathcal E$, the mapping $\xxi\mapsto F(\xx,\zz)$ is measurable, and if for almost every $\zz \in \Omega$, the mapping $\xx \mapsto F(\xx,\zz)$ is continuous. Conventionally, Caratheodory functions have played a key role in establishing statistical consistency of sample average approximations of stochastic optimization problems (see \cite{shapiro2021lectures}, Chapter 5). 
\begin{theorem}[\textbf{Weak Pareto efficiency of the extrapolated trajectory sequence}]
Suppose that the assumptions in Theorem \ref{thm:solution_convergence} 
are satisfied and the constraint function $\max_{k\in[K]}g_k$ is Caratheodory. Further suppose that the trajectory $\bar{\xx}^{(N)}$ in \eqref{eq:extrapolated-trajectory} is constructed from base level $\alpha_0 = \alpha_0^{(N)}$ satisfying $\alpha_0^{(N)} \rightarrow 0$ and $N\alpha_0^{(N)} \rightarrow \infty,$ as the number of observations  $N \rightarrow \infty.$ Then, any sequence of decisions $(\hat{\xx}^{(N)} : N \geq 1) \subset \mathcal{X}$ considered from the trajectories $(\bar{\xx}^{(N)}: N \geq 1)$ such that 
\begin{itemize}
    \item[i)] $\hat{\xx}^{(N)} \in \bar{\xx}^{(N)}$ for every $N \geq 1,$ and 
    \item[ii)] $p(\hat{\xx}^{(N)}) =  \Omega(N^{-b})$ for some $b \geq 1,$
\end{itemize}
is weakly Pareto efficient in the sense of Definition \ref{defn:weak-pareto-efficiency}. 
\label{thm:WPE-data-driven}
\end{theorem}
Since the exponent $b$ can be greater than 1, note  from Theorem \ref{thm:WPE-data-driven} that the weak Pareto efficiency holds even if the decisions $\hat{\xx}^{(N)}$ considered from the trajectory sequence $\{\bar{\xx}^{(N)}: N \geq 1\}$ are such that  $p(\hat{\xx}^{(N)}) \times N \rightarrow 0.$

\subsection{Optimality for the P-Model}
\label{sec:P-Model}
This section is devoted to discussing how the weak Pareto efficiency established in Theorem \ref{thm:WPE-data-driven} translates into  asymptotic optimality for the $P$-model counterpart of $\CCPalpha.$ 
Originally proposed by \cite{Charnes_Cooper_PModel}, the $P$-Model in \eqref{eq:P-Model}  minimizes the probability of constraint violation $p(\xx) = \Prob\left\{ \max_{k \in [K]} g_k(\xx,\xxi) > 0 \right\}$ while meeting a given cost target $B \in \R$:
\begin{align}
    {\tt P-Model}(B): \qquad \nu^\ast(B) := \min_{\xx \in \mathcal{X}} \ \log p(\xx) \quad \text{ s. to } \quad c(\xx) \leq B.
    \label{eq:P-Model}
\end{align}

Observe that any  $\xx \in \mathcal{X}^\ast(\alpha)$ is also optimal for the $P$-model counterpart  \eqref{eq:P-Model} when the budget parameter $B = v_{\alpha}^\ast.$ Now, given access to the distribution of $\xxi$ only via a limited number of independent observations $\{\xxi_1,\ldots,\xxi_N\},$ how should we solve for  an optimal solution to the $P$-Model($B$) at a given budget level $B$? This question becomes statistically challenging  when the budget parameter $B$ is  large enough to allow feasible decisions whose constraint violation probabilities $p(\xx) \leq 1/N.$ The sample based approximation $p_N(\xx)$ severely underestimates the true constraint violation probability $p(\xx)$ in this setting, primarily due to the same underestimation issue reported in Section \ref{sec:dd-underestimate}. Building on the desirable Pareto efficiency properties of the extrapolated trajectory $\bar{\xx}^{(N)},$ Algorithm \ref{algo:P-Model} below overcomes this underestimation issue by prescribing a decision $\hat{\xx}^{(N)}$ from the trajectory whose cost $c(\hat{\xx}^{(N)})$ precisely meets the given cost target $B.$

\begin{algorithm}[h!] 
  \caption{Extrapolation-based solution the $P$-Model($B$) when the budget $B$ is sufficient large to allow feasible decisions $\xx$ with $p(\xx) \leq 1/N$ }
  \ \vspace{-2pt}\\
  \KwIn{Dataset comprising $N$ i.i.d. samples $\{ \xxi_1,\ldots,\xxi_N\},$ cost target $B,$  $\alpha_0 \in (0,1)$ such that $N > C\alpha_0^{-1}(1-\alpha_0)$ for a large constant $C$} 
     \noindent \textbf{1. Solve ${\tt CCP}(\alpha_0)$}  via the SAA $\min\{c(\xx): p_{N}(\xx) \leq \alpha_0, \xx \in \mathcal{X}\},$ where $p_N(\xx)$ is  defined in \eqref{eq:saa-prob}. Let $\hat{\xx}_{\alpha_0}^{(N)}$ denote the resulting optimal solution.\\ 
     
  \noindent \textbf{2. Find the solution in  the trajectory $\bar{\xx}^{(N)}$ tightly meeting the cost target:}  Solve the one-dimensional line search  $t_\ast = \inf\{t \geq 1:c(t^r \hat{\xx}_{\alpha_0}^{(N)}) \geq B\}.$ For a cost  $c(\cdot)$ satisfying Assumption \ref{assume:cost}, $t_\ast$ is explicitly given by $t_\ast = [B/c(\hat{\xx}_{\alpha_0}^{(N)})]^{1/r}.$

\noindent \textbf{3. Return} 
 $\hat{\xx}^{(N)} = t_\ast^r \xx_{\alpha_0}^{(N)}$ as our candidate solution for  $P$-Model($B$). 
  \label{algo:P-Model}  
\end{algorithm}

\begin{corollary}[\textbf{Optimality of Algorithm \ref{algo:P-Model} for the P-model}] 
\label{cor:P-Model-data-driven}
Suppose that the conditions in Theorem \ref{thm:WPE-data-driven} are satisfied.  Letting  $\hat{\xx}^{(N)} =  t^{r} \hat{\xx}_{\alpha_0}^{(N)}$ for any fixed $t \geq 1,$ we obtain  
\begin{align*}
    \log p(\hat{\xx}^{(N)})  \sim \min\{\log p(\xx): c(\xx) \leq c(\hat{\xx}^{(N)}), \xx \in \mathcal{X}\}, \quad \text{ as } N \rightarrow \infty. 
\end{align*}
\end{corollary}
In words, Corollary \ref{cor:P-Model-data-driven} implies that the $P$-model's objective evaluated at $\hat{\xx}^{(N)}$ is asymptotically as good as that evaluated at any  $\xx \in \mathcal{X}$ whose cost does not exceed   $c(\hat{\xx}^{(N)}).$ 

\begin{num_example}
\em 
We consider solving the $P$-Model in \eqref{eq:P-Model} for the  data-driven transportation planning setting in Numerical Illustration \ref{num-eg:DD-Pareto-Front} with $30$ distribution centers. Given $N = 1000$ independent observations of the demand vector $\xxi,$ we compare in Figure \ref{fig:p-model} below the quality of solutions returned by Algorithm \ref{algo:P-Model} and the SAA $\min\{ \log p_{N}(\xx): c(\xx) \leq B, \xx \in \mathcal{X}\}$ for the $P$-Model. Considering different values for $B,$ the box plots in Figure \ref{fig:p-model} summarize the distribution of the out-of-sample constraint violation probability $p(\hat{\xx}^{(N)})$ observed over 32 independent experiments. 
\begin{figure}[h]
\includegraphics[width=\textwidth]{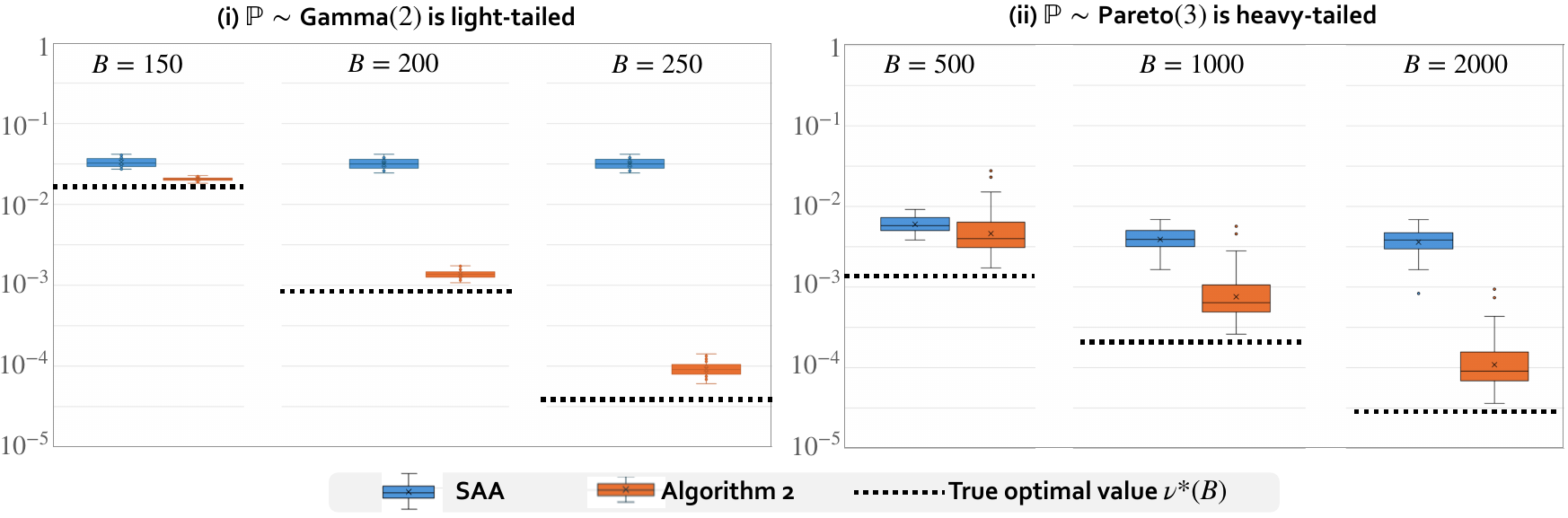}
  \caption{Out-of-sample constraint violation probabilites due to solutions returned by SAA and Algorithm \ref{algo:P-Model}}
\label{fig:p-model}
\end{figure}

From Figure \ref{fig:p-model}, we first observe that the constraint violation probabilities due to SAA solutions do not reduce below a certain level when the sample size $N$ is fixed, irrespective of how large the budget parameter $B$ is. As discussed before with SAA, this phenomenon is primarily due to the paucity of observations violating the constraints when $B$ is suitably large. The extrapolation based approach in Algorithm \ref{algo:P-Model}, however, can be seen to yield solutions whose constraint violation probabilities  improve significantly when the budget $B$ is increased. More interestingly, comparison with the smallest  constraint violation probability attainable for a given budget $B$ (indicated via the dashed line in Figure \ref{fig:p-model}) reveals the near-optimality of the solutions $\hat{\xx}^{(N)}$ even in the statistically challenging setting where  $N \times p(\hat{\xx}^{(N)}) \leq 1000 \times 10^{-4.1} \leq 0.08.$ \hfill $\Box$

\label{num-eg:p-model}
\end{num_example}



Despite the desirable properties of the extrapolated trajectories $\bar{\xx}^{(N)}$ noted in Theorem \ref{thm:WPE-data-driven} and Corollary \ref{cor:P-Model-data-driven}, it is important to note that  guaranteeing that a constraint of the form, $p(\hat{\xx}^{(N)}) \leq \alpha,$ is satisfied with high confidence remains a challenging problem when $\alpha \ll 1/N.$ Given $N$ samples and an acceptable constraint violation level $\alpha \leq 1/N$, it is statistically impossible to ensure a high-confidence guarantee on the probability of constraint violation  unless significantly stronger distribution assumptions are made; see \cite[Thm. 4.4.1]{deHaan} in the simpler univariate setting of quantile estimation.  Therefore, in the regime where we are interested in decisions whose constraint violation probabilities are smaller than $1/N,$ $P$-Model may be a more practical alternative than $\CCPalpha$ due to its relatively modest objective of  minimizing the constraint satisfaction probability for a given cost target $B.$  
In this case, the extrapolation procedure outlined  in Algorithm \ref{algo:P-Model} can serve as the basis for obtaining solutions with asymptotically the smallest  constraint violation probability attainable for the given budget $B.$  

\section*{Acknowledgments} The authors are grateful to Melvyn Sim for  helpful discussions. The research of K. Murthy is supported by the MOE AcRF Tier 2 award T2EP20223-0008.

\bibliographystyle{authordate1}
\bibliography{references}
\appendix

\section{Proofs of Main Results} 
\textbf{Notations.} We introduce the following notation  before discussing the proofs: For any $\mv{a} = (a_1,\ldots,a_d) \in \Real^d,$
$\mv{b} = (b_1,\ldots,b_d) \in \Real^d$ and $c \in \Real,$ we have
$\vert \mv{a} \vert = (\vert a_1 \vert, \ldots, \vert a_d \vert),$
$\mv{a}\mv{b} = (a_1b_1,\ldots,a_d b_d),$
$\mv{a}/\mv{b} = (a_1/b_1,\ldots,a_d/b_d),$
$\mv{a} \vee \mv{b} = (\max\{a_1,b_1\},\ldots,\max\{a_d,b_d\}),$
$\mv{a}^{\mv{b}} = (a_1^{b_1},\ldots,a_d^{b_d}),$
\sloppy{$\mv{a}^{-1} = (1/a_1,\ldots,1/a_d),$}
$\log\mv{a} = (\log a_1,\ldots,\log a_d),$
$c^{\mv{a}} = (c^{a_1},\ldots,c^{a_d}),$
denoting the respective component-wise operations. Further, we use $B(\xx,r)$ to denote the  ball $\{\yy: \Vert \xx - \yy \Vert_\infty \leq r\}.$ If the centre $\xx=0$, we abbreviate the norm-ball $B(\mv{0},r)$ as $B_r$. 

\subsection{Technical Background and Preliminary Lemmas}
We outline some additional notation used for the proofs. For a set $A$ let $[A]^{1+\delta}$ denote all the points that are at a distance of at most $\delta$ from $A$, that is $[A]^{1+\delta} = \{\xx: \exists \yy\in A: \|\xx-\yy\|\leq \delta\}$. Let $\partial A$ denote the boundary of $A$. Let ${\tt int}(A)$ and ${\tt cl}(A)$ respectively denote interior and closure of a set $A$. 
For a point $\yy$, and a set $A$ define $d(\yy,A) = \inf\{d(\xx,\yy):\xx\in A\}$. 
Recall that
\begin{align}\label{eqn:lim-of-set_equiv}
 \limsup_{n\to\infty} S_n &= \left\{\zz: \forall \varepsilon>0 \ \exists \text{ a sub-sequence } \{n_{i_k} :k\geq 1\} \ \text{such that } \forall k, \ \zz \in S_n+\varepsilon B_1\right\}\nonumber \\
    \liminf_{n\to\infty} S_n &= \left\{\zz: \forall \varepsilon>0 \ \exists n_0: \forall n\geq n_0, \ \zz \in S_n+\varepsilon B_1 \right\}
\end{align} 
\begin{definition}[Set Convergence]\label{def:set-lim}\em
We say that a sequence of sets $\{S_n:n\geq 1\}$ converges to a limit $S$ if as $n\to\infty$ $\limsup_{n\to\infty} S_n=\liminf_{n\to\infty} S_n=S$.
\end{definition}
We refer the reader to \cite[Chapter 4]{rockafellar2009variational} for additional examples and sufficient conditions required for set convergence. The following definition (see \cite{dembo2009large}, Pg. 6) is critical to this the subsequent proofs: 
\begin{definition}[Large Deviations Principle]\label{def:LDP}\em
    We say that a sequence of random vectors $\{\xxi_n:n\geq 1\}$ satisfies a large deviations principle (LDP) with rate function $J$ if as $n\to\infty$ 
    \begin{equation}\label{eqn:LDP}
       \limsup_{n\to\infty} n^{-1}
\log \mathbb P(\xxi_n\in F) \leq -\inf_{\xx\in F} J(\xx) \quad \text{and}\quad \liminf_{n\to\infty}n^{-1}\log \mathbb P(\xxi_n\in G) \geq  -\inf_{\xx\in G} J(\xx) 
    \end{equation}
for every closed set $F$ and open set $G$. We denote this by $\xxi_n\in {\tt LDP}(J)$.    
\end{definition}
\begin{proposition}\label{prop:ldp}
Suppose the density of $\xxi$ satisfies Assumption ($\mv{\mathscr{L}}$). Then $\{\xxi/\lambda^{\leftarrow}(t):t> 0\}\in {\tt LDP}(\varphi^*)$.     
\end{proposition}

When $\xxi$ satisfies Assumption ($\mv{\mathscr{H}}$), Proposition~\ref{thm:HTT_RV} which follows directly from \cite{Resnickbook}, Proposition 5.20 gives a characterization of the tail probability: 
\begin{proposition}\label{thm:HTT_RV}
     Suppose the density of $\xxi$ satisfies Assumption ($\mv{\mathscr{H}}$). Then for any set $A$ not containing the origin, we have
    \[
    \mathbb P\left(\xxi/t \in A\right)    \sim \lambda(t) \int_{\zz\in A }\varphi^*(\zz)d\zz :=\lambda(t) \nu(A)
    \]
\end{proposition}
\begin{lemma}\label{lem:RVs}\em
For $t>0$ and $\xx$, $\varphi^*(t\xx) = t^{c}\varphi^*(\xx)$ where $c = \gamma$ if Assumption~($\mathscr{L}$)  and $c= -(\gamma+d)$ if Assumption~($\mathscr{H}$) holds.    
\end{lemma}

\begin{lemma}\label{lem:I_homogenity}
    The function $I(\cdot)$ defined in \eqref{eqn:I_star_limit} satisfies $I(c\xx) = c^{-\gamma/r}I(\xx)$ for all $\xx$ and $c>0$.
\end{lemma}



\subsection{Proofs from Section~\ref{sec:Modelling Assumptions}} 

\noindent \textbf{Proof of Lemmas~\ref{lem:marginals-light}
and \ref{lem:marginals-heavy}: }First suppose that Assumption~($\mathscr{L}$) holds. Then, from Proposition~\ref{prop:ldp}, the sequence of random vectors $\{\xxi/\lambda^{\leftarrow}(n) :n\geq 0\}\in {\tt LDP}(\varphi^*)$. With $n=\lambda(t)$, note that then for each $i$, 
\begin{align*}
    \log \mathbb P(\xi_i \geq tz) &\sim -\lambda(t) \inf_{z_i \geq z}\varphi^*(\zz)\\
    &=-\lambda(t) \inf_{p_i>1} \varphi^*(\pp/z) = -\lambda(t) z^{-\gamma} \inf_{p_i>1}\varphi^*(\pp)
\end{align*}
The first statement above follows directly from the large deviations principle, and the second from Lemma~\ref{lem:RVs}. Setting the infimum above as $c_i$ completes the proof. Now, we have that 
\[
\bar F_i(tz)  =\exp\left\{-c_i z^{-\gamma} \lambda(t) [1+o(1)]\right\}
\]
and taking the maximum over all $i$,
\[
\bar F(tz) =\exp\left\{-\min_{i\in[d]} c_i z^{-\gamma} \lambda(t) [1+o(1)]\right\}
\]
Setting $\lambda(t) = -\log \bar F(t)$ above,  it follows that the multiplying constant $\min_i c_i =1$. Next suppose that Assumption~($\mathscr{H}$) holds. Then from Proposition~\ref{thm:HTT_RV}, note that 
\begin{align*}
    \mathbb P(\xi_i \geq tz) &\sim \bar F(t) \int_{z_1\geq z}\varphi^*(\zz)d\zz \\
    & = \bar F(t) z^{-\gamma}\int_{z_1\geq 1}\varphi^*(\zz)d\zz \quad 
\end{align*}
The second equality above follows upon changing variables from $z_1$ to $\zz/z$ and then using Lemma~\ref{lem:RVs}.
Setting $\nu(A_1)$ above as $c_i$ gives the first part of the proof. Taking the maximum over all $i$,
\[
\bar F(tz) = \max_{i\in[d]} c_i z^{-\gamma} \lambda(t) [1+o(1)]
\]
Setting $\lambda(t) = \bar F(t)$ above,  it follows that the multiplying constant $\max_i c_i =1$.
\qed
\subsection{Proofs from Section~\ref{sec:Asymptotics}} 
We introduce some simplifying notation which we will carry forward through the rest of the proofs.
Denote the function 
$g_t(\xx,\zz) = \max_{i} t^{-k} g_i(t^r\xx,t\zz)$.
Let 
\begin{equation}\label{eqn:level_sets_of_functions}
    \mathcal S_{t,\alpha,M}(\xx) = \{\zz: g_t(\xx,\zz) \leq \alpha\} \cap B_M \quad \text{and} \quad \mathcal S^*_{\alpha,M}(\xx) = \{\zz: g^*(\xx,\zz) \leq \alpha \}\cap B_M. 
\end{equation}
denote the $\alpha-$level set of $g_t(\xx,\cdot)$, and $g^*(\xx,\cdot)$ restricted to the $M-$norm ball. If $M=\infty$, we simply write $\mathcal S_{t,\alpha}$ and $\mathcal S^*_{t,\alpha}$ respectively. In what follows, let $\Upsilon(\xx) = \mathcal [S(\xx)]^{c}$ denote the set of realisations when the chance constraint is violated. Observe now, that \eqref{eqn:CCP} can be equivalently re-stated as
\begin{equation}\label{eqn:CCP_unsafe}
\min c(\xx) \quad \text{ s.t }\quad \mathbb P[\xxi\in \Upsilon(\xx) ] \leq \alpha.
\end{equation}

Defining $\Upsilon_t(\xx) := t^{-1}\Upsilon(t^r \xx)$, further note that
\begin{align*}
    \Upsilon_t(\xx) &= \left\{\zz/t: \max_i g_i(t^r\xx,\zz) > 0\right\} \\ 
    & = \left\{\zz: g_t(\xx,\zz) > 0\right\} \text{ replacing $\zz/t $ by $ \zz$ and multiplying by $t^{-\rho}$}\\
    & = \mathcal S_{t,0}^c(\xx) .
\end{align*}

The following technical lemmas are required to proceed with the subsequent proofs:
\begin{lemma}\label{lem:level_sets_containment}
    Suppose that Assumption~\ref{assume:hom-safe-set} or \ref{assume:non-hom-safeset} holds and that $\xx_t\to\xx \neq \mv 0$. Then, for every $(\varepsilon, M)>0$, there exists an $t_0$ such that for all $t>t_0$, 
\begin{equation}\label{eqn:level_sets_containment}
        \mathcal S_{t, 0,M}(\xx_t) \subseteq [\mathcal S^*_{0, M}(\xx)]^{1+\varepsilon} \quad \text{ and }\quad  \mathcal S_{-\varepsilon, M}^*(\xx) \subseteq \mathcal S_{t,0,M}(\xx_t)  
    \end{equation}
\end{lemma}
\begin{lemma}\label{lem:compact_containment}
    There exists $\alpha_0>0$  and a compact set $K$ such that for all $\alpha<\alpha_0$, $\mathcal Y^*_\alpha \subset K$.
\end{lemma}
Let 
\[
I(\xx) = \begin{cases}
    \inf{\varphi^*(\xx) : g^*(\xx,\zz) \geq 0}\}]^{-1}\ [\min_i c_i] \quad\quad \quad \text{ if Assumption~($\mathscr{L}$) holds}\\
    \int_{g^*(\xx,\zz)\geq 0}\varphi^*(\zz)d\zz\ [\max_i c_i]^{-1} \quad\quad \quad\quad\quad\quad \text{ if Assumption~($\mathscr{H}$) holds}
\end{cases}
\]
\begin{lemma}\label{lem:extended_functional}
    Define the extended functionals:
\begin{equation}\label{eqn:Ext_functional}
    \bar{c}_t(\xx) = c(\xx) + \chi_{\Lev_1(I_t)}(\xx) \quad \text{ and }\quad c^*(\xx) = c(\xx) + \chi_{\Lev_1(I^*)}(\xx)
\end{equation}
Then, $\bar c_t\to c^* $ epigraphically.
\end{lemma}

\noindent \textbf{Proof of Proposition~\ref{prop:prob-convergence_LT}:} 
Denote the function 
\[
G_t(\xx) = \frac{\log\mathbb P \left(g_t(\xx, \xxi_t) > 0\right) }{\lambda(t)}
\]
Note that proving the uniform convergence in Proposition~\ref{prop:prob-convergence_LT} is then equivalent to proving the continuous convergence of $G_t \to -G_0$ (see \cite{rockafellar2009variational}, Theorem 7.14). To this end, rewrite using the notation defined previously 
\begin{equation}\label{eqn:I_t}
    G_t(\xx) = \frac{\log \mathbb P\left( \xxi_t \in \Upsilon_t(\xx) \right)}{\lambda(t)}.
\end{equation}

 \noindent \textbf{Step 1) Upper Bound: } Recall that $\Upsilon_t(\xx_t) = \mathcal S_{t,0}^c(\xx_t)$ and note the following sequence of set operations:
 \begin{align*}
      \mathcal S_{t,0}^c(\xx_t) &=   (\mathcal S_{t,0}^c(\xx_t) \cap B_{2M})  \cup (\mathcal S_{t,0}^c(\xx_t) \cap B_{2M}^c) \\
      & \subseteq ([(\mathcal S_{t,0}(\xx_t) \cap B_M) \cup (\mathcal S_{t,0}(\xx_t) \cap B_M^c)]^c \cap B_{2M})\cup B_{2M}^c \\
      &\subseteq (\mathcal S_{t,0,M}^c(\xx_t) \cap B_{2M}) \cup B_{2M}^c \quad \,\text{ dropping the union above}.
 \end{align*}
 Then, observe that
 \begin{equation}\label{eqn:contain_upper}
\mathbb P\left(\xxi_{t} \in \Upsilon_{t}(\xx_t) \right) \leq \mathbb P(\xxi_{t} \in [\mathcal S_{0,t,M}(\xx)]^c \cap B_{2M}) + \mathbb P(\xxi_{t} \in  B_{2M}^c).
 \end{equation}
The second containment in \eqref{eqn:level_sets_containment} implies that for a fixed $\delta>0$,
\[
    \mathbb P(\xxi_{t} \in [\mathcal S_{0,t,M}(\xx)]^c \cap B_{2M}) \leq \mathbb P(\xxi_{t} \in  [\mathcal S^*_{-\delta,M}(\xx)]^{c} \cap B_{2M})  \text{ for all large enough $t$}. 
\]
Consequently, for all large enough $t$, one has that 
\[
G_t(\xx_t ) \leq \frac{\log \left\{\mathbb P(\xxi_{t} \in \mathcal [\mathcal S^*_{-\delta,M}(\xx)]^{c} \cap B_{2M}) +  \mathbb P(\xxi_{t} \in  B_{2M}^c)\right\}}{\lambda(t)}.
\]
Recall that from Proposition~\ref{prop:ldp}, the sequence $\xxi_{t}\in {\tt LDP}(\varphi^*)$. Upon applying the LDP upper bound from \eqref{eqn:LDP}, 
\begin{equation*}
\limsup_{t\to\infty} G_t(\xx_t)\leq - \min\left\{  \inf_{\zz\in [\mathcal{S}^*_{-\delta,M}(\xx)]^c \cap B_{2M}} \varphi^*(\zz), \inf_{\zz \in B_{2M}^c} \varphi^*(\zz)\right\}. 
\end{equation*}
Note that as $\delta\to 0$, the sets $[\mathcal S_{-\delta,M}^*(\xx)]^{c}\cap B_{2M}$ decrease to $[\mathcal S_{0,M}^*(\xx)]^{c}\cap B_{2M}$ and are relatively compact in the sense that the closure of
\[
\bigcup_{\delta>0} [\mathcal S_{-\delta,M}^*(\xx)]^{c}\cap B_{2M} \text{ is compact} 
\]

Following \cite{Langen} Theorem 2.2, for a fixed $\varepsilon>0$, whenever  $\delta<\delta_0$, $\inf_{\zz\in [\mathcal S_{-\delta,M}^*(\xx)]^{c}\cap B_{2M}} \varphi^*(\zz) \geq \inf_{\zz\in [\mathcal S_{0,M}^*(\xx)]^{c}\cap B_{2M}} \varphi^*(\zz) -\varepsilon$. Next, observe that $\mathcal [S_{0,M}^*(\xx)]^c\cap B_{2M} = [\Upsilon^*(\xx) \cup B_{M}^c] \cap B_{2M}  $. Consequently,
\[
 \inf_{\zz\in [\mathcal S_{0,M}^*(\xx)]^{c}\cap B_{2M}} \varphi^*(\zz)  = \min\left\{\inf_{\zz\in  \Upsilon^*(\xx) \cap B_{2M} } \phi^*(\zz), \inf_{\|\zz\|\in [M,2M]} \varphi^*(\zz) \right\}
\]
For the above choice of $\delta$, 
\begin{equation}\label{eqn:LDP_upper}
\limsup_{t\to\infty} G_t(\xx_t) \leq -\min\left\{  \inf_{\zz\in \Upsilon^*(\xx)} \varphi^*(\zz) -\varepsilon, \inf_{\zz \in B_{M}^c} \varphi^*(\zz)\right\}. 
\end{equation}
where we can replace $\Upsilon^*(\xx) \cap B_{2M}$ by $\Upsilon^*(\xx)$ in the infimum since $\Upsilon^*(\xx) \cap B_{2M} \subset \Upsilon^*(\xx)$.  

Note that from Lemma~\ref{lem:RVs}, $\varphi^*$ is a $\gamma$-homogeneous function. Therefore, by selecting $M$ sufficiently large, the second infimum above can be made as large as desired. 
Whenever the non-vacuousness condition holds,
 note that $\Upsilon^*(\xx)\neq \emptyset$. Then, since $\varphi^*$ is $\gamma-$homogeneous, for all $\xx$, $G_0(\xx) = \inf\{\varphi^*(\zz):\zz\in\Upsilon^*(\xx)\}< \infty $.
Further, 
Hence, there exists an $t_0$ such that for all $t>t_0$, 
\begin{equation}\label{eqn:lower_I}
    G_t(\xx_t)  \leq  - G_0(\xx) +\varepsilon.
\end{equation}

\noindent \textbf{ii) Lower Bound: }We now prove a lower bound to match \eqref{eqn:lower_I}. For any $(\delta,M)>0$,
\begin{align*}
    \mathbb P(\xxi_{t} \in \Upsilon_{t}(\xx_t)) &\geq  \mathbb P(\xxi_{t} \in \Upsilon_{t}(\xx_t)\cap B_{M})\\
    &= \mathbb P(\xxi_{t} \in \mathcal [S_{t,0,M}(\xx_t)]^c\cap B_{M}) 
\end{align*}
To see the last part, note that 
\[
\mathcal  S_{t,0,M}^c(\xx_t) \cap B_M= [\mathcal S_{t,0}(\xx_t)\cap B_M]^{c} \cap B_M = [\mathcal S_{t,0}^c(\xx_t)\cup B_M^c]\cap B_M = \mathcal S_{t,0}^c(\xx_t)\cap B_M
\]
Utilizing the first containment in \eqref{eqn:level_sets_containment}, 
\begin{align}\label{eqn:lb}
    \mathbb P(\xxi_{t} \in \Upsilon_{t}(\xx_t)) &\geq \mathbb P(\xxi_{t} \in ([\mathcal S_{0,M}^*(\xx)]^{1+\delta})^c \cap B_{M}), 
\end{align}
where the statement above is for all $t>t_0$, chosen in accordance with Lemma~\ref{lem:level_sets_containment}.
Invoke the large deviations lower bound to get that 
\[
\liminf_{n\to\infty} \frac{\log \mathbb P(\xxi_{t} \in ([\mathcal S_{0,M}^*(\xx)]^{1+\delta})^c \cap B_{M})}{\lambda(t)}  \geq -\inf_{\zz\in [\mathcal S_{0,M}^*(\xx)]^{1+\delta})^c \cap B_{M}}\varphi^*(\zz).
\]
To conclude notice the convergence $([\mathcal S_{0,M}^*(\xx)]^{1+\delta})^c \cap B_{M} \nearrow \Upsilon^*(\xx) \cap B_M$ as $\delta\downarrow 0$, and that the right hand side in this convergence is relatively compact.
Therefore, upon an application of \cite{Langen}, Theorem 2.2(ii), for a small enough $\delta$,
\[
\inf_{\zz\in [([\mathcal S_{0,M}^*(\xx)]^{1+\delta})^c \cap B_{M}]}\varphi^*(\zz) \leq \inf_{\zz\in \Upsilon^*(\xx) \cap B_{M}}\varphi^*(\zz) +\varepsilon.
\]
Plugging the above display into the large deviations lower bound: 
\[
\liminf_{t\to\infty} G_t(\xx_t) \geq -\inf_{\zz\in \Upsilon^*(\xx) \cap B_{M}}\varphi^*(\zz) -\varepsilon = -G_0(\xx)-\varepsilon.
\]
Since $M$ and $\varepsilon$ are arbitrary, combining the above convergence with \eqref{eqn:lower_I} yields that  $G_t(\xx_t) \to -G_0(\xx)$, or $G_t$ converges continuously and therefore uniformly to $-G_0$.  \qed

\noindent \textbf{Proof of Proposition~\ref{prop:prob-convergence_HT}: }
Define 
\[
G_t(\xx) = \frac{\mathbb P(\xxi \in \Upsilon_t(\xx))}{\lambda(t)}
\]
Similar to the proof of Proposition~\ref{prop:prob-convergence_LT}, it is sufficient to 
show that $G_t\to G_0$. 

\noindent \textbf{i) Upper Bound:} Recall that from \eqref{eqn:contain_upper}, 
\[
\mathbb P\left(\xxi_{t} \in \Upsilon_{t}(\xx_t) \right) \leq \mathbb P(\xxi_{t} \in  [\mathcal S^*_{-\delta,M}(\xx)]^{c} \cap B_{2M})  + \mathbb P(\xxi_t\in B_{2M}^c)\text{ for all large enough $t$} 
\]
Observe that the non-vacuousness of the set-valued map $\mathcal S^*$ implies that $\mv 0 \not \in \mathcal S^*_{0,M}(\xx)$. Note now that for all small enough $\delta$, $\mv 0 \not \in [\mathcal S_{-\delta, M}(\xx)]^{c}$. Using Proposition~\ref{thm:HTT_RV},
\[
[\lambda(t)]^{-1} \mathbb P(\xxi_{t} \in  [\mathcal S^*_{-\delta,M}(\xx)]^{c} \cap B_{2M}) \to \nu([\mathcal S^*_{-\delta,M}(\xx)]^{c} \cap B_{2M}) \text{ as  } t \to \infty
\]
Further note that $[\mathcal S^*_{-\delta,M}(\xx)]^{c}\cap B_{2M} \uparrow [\mathcal S^*_{0,M}(\xx)]^c\cap B_{2M}$. Then, one has that for all small enough $\delta$,
$\nu(\mathcal S^*_{-\delta,M}(\xx)) \leq \nu (\mathcal [S^*_{0,M}(\xx)]^{c} \cap B_{2M}) + \varepsilon$. 
This implies that
\[
\limsup_{t\to\infty} [\lambda(t)]^{-1}\mathbb P\left(\xxi_{t} \in \Upsilon_{t}(\xx_t) \right)\leq \nu([\mathcal S^*_{0,M}(\xx)]^c\cap B_M) + \nu(B_M^c) +\varepsilon
\]
Finally, since $(\varepsilon,M)$ are arbitrary $\limsup_{t\to\infty} [\lambda(t)]^{-1}\mathbb P\left(g(t^r\xxi,\xx_t)>0\right) \leq \nu(\Upsilon^*(\xx))$. 

\noindent \textbf{ii) Lower Bound: }For the matching lower bound, we start with \eqref{eqn:lb}: 
\[
P(\xxi_t\in \mathcal U_t(\xx_t)) \geq P(\xxi_t \in ([\mathcal S_{0,M}^*(\xx)]^{1+\delta})^{c} \cap B_M)
\]
Now, apply Proposition~\ref{thm:HTT_RV} to obtain that 
\[
\liminf_{t\to\infty } [\lambda(t)]^{-1} P(\xxi_t\in \mathcal U_t(\xx_t)) \geq \nu([\mathcal S_{0,M}^*(\xx)]^{1+\delta})^c \cap B_M)
\]
Since $(\delta,M)$ are arbitrary, the above display implies that $\liminf_{t\to\infty } [\lambda(t)]^{-1} P(\xxi_t\in \mathcal U_t(\xx_t)) \geq \nu(\Upsilon^*(\xx))$. Noting that $\nu(\Upsilon^*(\xx)) = \int_{g^*(\xx,\zz)\geq 0}\varphi^*(\zz)d\zz$ completes the proof.\qed

\noindent \textbf{Proof of Theorem~\ref{thm:solution_convergence}: } 
With $t = s_\alpha$ (we will suppress the dependence of $t$ on $\alpha$ to simplify notation), recall that from Lemmas~\ref{lem:marginals-light}-\ref{lem:marginals-heavy}, 
\[
\lambda(t) = \begin{cases}
- \log \alpha  [\min_i c_i](1+o(1)) \quad\quad  \text{ if Assumption~($\mathcal L$) holds and}\\
 \alpha(1+o(1)) / [\max_i c_i] \quad  \quad\quad\quad\text{ if Assumption~($\mathcal H$) holds }
\end{cases}
\]
 Now, note that the \eqref{eqn:CCP}
can be re-written as (using the homogeneity of $c(\cdot)$)
\[
t^r \min\{ c(\xx/t^r) : I_t(\xx/t^r) \leq 1\} \text{ where } I_t(\xx) =\begin{cases}
    [\min_i c_i+o(1)]/ \times G_t(\xx) \quad \text{ if Assumption~($\mathcal L$) holds and}\\
    G_t(\xx)/[\max_i c_i] \quad \ \ \ \ \text{ if Assumption~($\mathcal H$) holds}
\end{cases}
\]
Let $ y_\alpha^*$ and $\mathcal Y_\alpha^*$ be the optimal value and solutions of the optimization problem 
\[
O_\alpha: \min\{c(\pp) : I_t(\pp) \leq 1\}
\]

Lemmas~\ref{lem:compact_containment} and \ref{lem:extended_functional} now imply that $\{O_\alpha:\alpha>0\}$ satisfy the conditions of \cite{bonnans2013perturbation}, Proposition 4.4 (see the discussion below Proposition 4.6 therein), with the limiting function
\[
I(\xx) = \begin{cases}
    \inf{\varphi^*(\xx) : g^*(\xx,\zz) \geq 0}\}]^{-1}\ [\min_i c_i] \quad\quad \quad \text{ if Assumption~($\mathscr{L}$) holds}\\
    \int_{g^*(\xx,\zz)\geq 0}\varphi^*(\zz)d\zz\ [\max_i c_i]^{-1} \quad\quad \quad\quad\quad\quad \text{ if Assumption~($\mathscr{H}$) holds}
\end{cases}
\]
Therefore, $y_\alpha^* \to v^*$ and $d(\yy_\alpha^*,\mathcal Y^*) \to 0$ for any $\yy_\alpha\in \mathcal Y_\alpha^*$.
To conclude the proof of the main theorem, note that $\mathcal X_{\alpha}^* = t^r \mathcal Y_\alpha^*$ and that $v_\alpha^* = t^r y_\alpha $. Observing that as a result, $\xx_\alpha\in \mathcal X_\alpha^*$ can be written as $t^r \yy_\alpha^*$ for some $\yy_\alpha^*\in\mathcal Y_\alpha^*$, concludes the proof.\qed

\noindent \textbf{Proof of Lemma~\ref{lem:salpha-growth-rate}: }Under Assumption~($\mathscr{L}$) note that for $s>0$,
\[
-\log F_{1}(st) \sim c_1 s^{\gamma}\lambda(t) \quad \text{ and } -\log F_{1}(st) \sim c_1 \lambda(st)
\]
This implies that $\lambda(st)/\lambda(t) \to s^\gamma$ for all $s>0$ or that $\lambda\in \RV(\gamma)$. Further note that $s_\alpha = \lambda^{\leftarrow}(\log(1/\alpha))$. From \cite{deHaan}, Proposition B.1.9, $\lambda^{\leftarrow}\in \RV(1/\gamma)$, and therefore, we have that 
\[
s_\alpha = [\log(1/\alpha)]^{1/\gamma} L(\log(1/\alpha)), \text{ where } L\in \RV(0)
\]
Taking logarithms on both sides, one gets that 
\[
\log s_\alpha \sim \frac{1}{\gamma} \log \log(1/\alpha)
\]
Next, suppose that Assumption~($\mathscr{H}$) holds. Here, following the above steps, $\bar F(t) \in \RV(-\gamma)$. Consequently, $s_\alpha \in \RV(1/\gamma)$ and the conclusion of the Lemma follows. \qed

\subsection{Proofs from Section~\ref{sec:Robust_Asymototics}} 
\label{sec:DRO-proofs}
The following technical lemma is used in the proofs in this section
\begin{lemma}\label{lem:solution_of_s}
    Let $s(t)$ solve 
    \[
    {\eta} = t\frac{f(s)}{s} + \left(1-\frac{t}{s}\right)f\left(\frac{1-t}{1-t/s}\right).
    \]
Then, $s(t) = g(\eta/t+\kappa(t))$ where $\kappa(t)\to 0$ as $t\to 0$.    
\end{lemma}

\noindent \textbf{Proof of Theorem~\ref{thm:f-div-DRO}:} The proof of Theorem~\ref{thm:f-div-DRO} relies crucially on \cite{blanchet2020distributionally}, Corollary 1, which in our context states that the worst case probability evaluation $p_{\tt wc}(\xx) = \sup_{P\in\mathcal P}P(F(\xx,\xxi) >0)$ is related to $p(\xx)$ by the equations:
\begin{align}
    \eta &= p_{\tt wc}(\xx) \frac{f(s)}{s} + \left(1-\frac{p_{\tt wc}(\xx)}{s}\right)f\left(\frac{1-p_{\tt wc}(\xx)}{1-p_{\tt wc}(\xx)/s}\right) \text{ and }\label{eqn:1st_condition_fdiv}\\
    p_{\tt wc}(\xx) &= s p(\xx)\label{eqn:2nd_condition_fdiv}
\end{align}
where $s=s(\xx)$ depends on the decision $\xx$, and is unique.
Upon substituting $t = p_{\tt wc}(\xx)$ in Lemma~\ref{lem:solution_of_s}, write $s(p_{\tt wc}(\xx)) = g(\eta/p_{\tt wc}(\xx)+\kappa(p_{\tt wc}(\xx))$. From \eqref{eqn:2nd_condition_fdiv}, 
\[
p(\xx) = \frac{p_{\tt wc}(\xx)}{g\left({\eta}/{p_{\tt wc}(\xx)}+\kappa(p_{\tt wc}(\xx)\right)}
\]
Recall that over solution which are feasible to $\DROCCPalpha$, one has that $\{p_{\tt wc}(\xx) \leq \alpha\}$. Since the map $x \mapsto x/g(\eta /x+\kappa(x))$  is (for all large enough $x$) monotone increasing, this implies that the feasible region of $\DROCCPalpha$ can be written as
\[
\left\{\xx: p(\xx) \leq \frac{\alpha}{g(\eta/\alpha+\kappa(\alpha))}\right\}
\]
To complete the rest of the proof we repeat the crucial steps from the proof of Propositions~\ref{prop:prob-convergence_LT}-\ref{prop:prob-convergence_HT} and  Theorem~\ref{thm:solution_convergence}. First, note  $p_{\tt wc}(\xx) \leq \alpha$ over the feasible region of $\DROCCPalpha$. 
Note that $\kappa(\alpha) = o(1)$, and therefore, $g(\eta/\alpha +\kappa(\alpha)) = g(\eta/\alpha)(1+o(1))$. Now,
define $\alpha_0 = \alpha/g(\eta/\alpha)$, and note that with $t=s_{\alpha_0}$ in the notation of Theorem~\ref{thm:solution_convergence}, $\DROCCPalpha$ can be re-written as
\[
t^r \min\{c(\xx/t^r) : I_t(\xx/t^r)(1+o(1)) \leq 1\} \text{ where $I_t(\xx)$ is as defined in Theorem~\ref{thm:solution_convergence}}.
\]
Observe that using the same arguments as from Propositions~\ref{prop:prob-convergence_LT}-\ref{prop:prob-convergence_HT}, the RHS above converges uniformly over compact sets to $I(\cdot)$.
Now, following the proof of  Theorem~\ref{thm:solution_convergence}, the optimal value of $\DROCCPalpha$ satisfies $v_{\alpha}^{f} \sim s_{\alpha_0}^{r} v^*$ and that the distance $d(s_{\alpha_0}^{-r}\xx_{\alpha}^{f}, \mathcal X^*)\to 0$. Noting that $s_{\alpha_0} = \bar F^{\leftarrow}(\alpha_0)$ completes the proof.\qed


\noindent \textbf{Proof of Proposition~\ref{prop:exp-div} and Table~\ref{tab:f-div-scalepreserving}: }
Observe that from Theorem~\ref{thm:f-div-DRO}, the scaling rate $t_\alpha$ governs the rate of increase of costs. To prove the strong (resp.) weak rate preserving property, it is therefore sufficient to prove that for all $\mathbb Q$ in $\mathcal P,$ $\limsup_{\alpha\to 0} t_\alpha/s_\alpha <\infty$ (resp. $\limsup_{\alpha\to 0} \log t_\alpha/\log s_\alpha \leq 1$), where $t_\alpha$ and $s_\alpha$ are as defined previously.
First, note that 
\begin{equation}\label{eqn:t_ratio}
    \frac{t_{\alpha}}{s_{\alpha}} = \frac{\bar F^{\leftarrow}(\alpha/g(\eta/\alpha)) }{\bar F^{\leftarrow}(\alpha)} \text{ and } \frac{\log t_{\alpha}}{\log s_{\alpha}}  = \frac{\log \bar F^{\leftarrow}(\alpha/g(\eta/\alpha))}{\log \bar F^{\leftarrow}(\alpha)}
\end{equation}
 
\noindent  \textbf{a)} Here, suppose that $f(x)/x =\Omega(x^p)$ for some $p>0$. Then, note that $g(x) = o(x^{1/p})$. When Assumption~($\mathscr{L}$) holds, recall that from Lemma~\ref{lem:salpha-growth-rate}, 
$\log \bar F^{\leftarrow}(t) \sim \gamma^{-1}\log \log(1/t)$ as $t\to 0$. Then,
\[
\bar F^{\leftarrow}(\alpha/g(\eta/\alpha)) = \exp(\gamma^{-1}\log \log (o((\eta/\alpha)^{1/p})/\alpha) ) = \exp(\gamma^{-1}\log \log (o((1/\alpha)^{1/p+1}))
\]
Substituting this into \eqref{eqn:t_ratio}, 
\begin{equation}\label{eqn:lt_lim}
    \varlimsup_{\alpha\to 0}\frac{t_\alpha}{s_\alpha} \leq 
\varlimsup_{\alpha\to 0} \frac{\exp(\gamma^{-1}\log \log ((1/\alpha)^{1/p+1})}{\exp(\gamma^{-1}\log \log (1/\alpha))} <\infty 
\end{equation}
Suppose now that Assumption~($\mathscr{L}$) holds.
Note that \eqref{eqn:lt_lim} also demonstrates when $f(x)/x \sim \exp(x)$, $t_\alpha = O(s_\alpha)$, since then $\bar F^{\leftarrow}(\alpha/g(\eta/\alpha)) = \exp(\gamma^{-1}\log \log (1/\alpha)^{1+o(1)})$.

\noindent \textbf{b)} Here, suppose $f(x)/x \sim \exp(x^p)$ for some $p$. This implies that $g(x) = o((\log x)^{1/p})$. If Assumption~($\mathscr{H}$) holds in addition, then note further that 
by Lemma~\ref{lem:salpha-growth-rate}, $\log \bar F^{\leftarrow}(t) \sim \gamma^{-1}\log (1/\alpha)$. 
Then,
\[
\log \bar F^{\leftarrow}(\alpha/g(\eta/\alpha))  = \gamma^{-1} \log(g(\eta/\alpha)/\alpha)(1+o(1)) = \gamma^{-1}\log(1/\alpha) (1+o(1)),
\]
where the last statement is because $g(x) = o((\log x)^{1/p})$. 
Substituting this into \eqref{eqn:t_ratio}, 
\[
\varlimsup_{\alpha\to 0}\frac{\log t_\alpha}{\log s_\alpha} = 
\varlimsup_{\alpha\to 0} \frac{\gamma^{-1}(1+o(1))\log(1/\alpha)}{\gamma^{-1}(1+o(1))\log(1/\alpha)} \leq 1 
\]
Suppose instead Assumption~$(\mathscr{L})$ holds. From Lemma~\ref{lem:salpha-growth-rate}, $\log \bar F^{\leftarrow}(t) \sim \gamma^{-1}\log \log (1/\alpha)$. 
\noindent \textbf{c) }Suppose that $u(x) :=f(x)/x \sim x^{p}$ for some $p>0$. Then, note that $g = u^{\leftarrow} \in \RV(1/p)$. If Assumption~($\mathscr{H}$) holds in addition, then note further that 
by Lemma~\ref{lem:salpha-growth-rate}, $\log \bar F^{\leftarrow}(t) \sim \gamma^{-1}\log (1/\alpha)$. Substituting this into \eqref{eqn:t_ratio}, 
\[
 \bar F^{\leftarrow}(\alpha/g(\eta/\alpha))) =   \exp(\gamma^{-1}\log(g(\eta/\alpha)/\alpha))(1+o(1)) ) = \exp((1/p+1)\gamma^{-1})\log(\eta/\alpha))(1+o(1)))
\]
Consequently, since $p>0$,
\[
\varliminf \frac{\log t_{\alpha}}{\log s_{\alpha}} = \frac{((1/p+1)\gamma^{-1})\log(1/\alpha)(1+o(1)))}{(\gamma^{-1} \log(1/\alpha)(1+o(1)))} >1  
\]

\noindent \textbf{d) }Next suppose that $u(x) \sim (\log x)^{p}$ for some $p>0$ and Assumption~($\mathscr{L}$) holds. Since $g = u^{\leftarrow
}$, one has that $g(x) = \exp(x^{1/p}(1+o(1)))$. From Lemma~\ref{lem:salpha-growth-rate} note that here 
$\log \bar F^{\leftarrow}(t) \sim \gamma^{-1}\log \log(1/t)$ as $t\to 0$. 
Substituting this into \eqref{eqn:t_ratio}, 
\begin{align*}
    \bar F^{\leftarrow}(\alpha/g(\eta/\alpha))& = \exp(\gamma^{-1}\log \log (\exp(\eta/\alpha)^{1/p}/\alpha) (1+o(1))) \\
    &= \exp((p\gamma)^{-1}\log(\eta/\alpha)(1+o(1))) \text{ for every }\mathbb Q\in \mathcal P
\end{align*}

Consequently, 
\[
\frac{t_{\alpha}}{s_{\alpha}}  = \frac{\exp((p\gamma)^{-1}(1+o(1))\log(\eta/\alpha))}{\exp(\gamma^{-1}(1+o(1))\log \log (1/\alpha))} \to \infty \text{ as }\alpha\to 0. 
\]
\qed

\noindent \textbf{Proof of Theorem~\ref{lem:Wasserstein_CCP}: }
Note that for any $\mathbb Q_0 \in \mathcal P$, the cost of $\DROCCPalpha$ is at least as much as $v_\alpha^*(\mathbb Q)$.
In order to complete the proof, we hence demonstrate that so long as Assumption~($\mathscr{L}$) or ($\mathscr{H}$) is satisfied, there exists a measure $\mathbb Q_0$ which satisfies Assumption~($\mathscr{H}$) with $\gamma= p+\varepsilon$ for any $\varepsilon>0$. Consequently, one has that $c(\xx_{\alpha}^{\tt W}) \geq v \alpha^{r/p +\varepsilon}$. 

We use the following approach to show the existence of such as measure $\mathbb Q_0$: suppose there exists a coupling $\mv{\Pi}$ of $\xxi$ and $\xxi_0$ (where $\xxi_0$ has the distribution $\mathbb Q_0$) such that
$E_{\Pi}[\|\xxi-\xxi_0\|_p]\leq \eta$ and $v_\alpha^*(\mathbb Q_0)\geq v \alpha^{r/p +\varepsilon}$. . 
Then, 
\[
d_{W}(\mathbb P,\mathbb Q_0) \leq E_{\Pi}[\|\xxi-\xxi_0\|_p], \text{ from the definition of Wasserstein distance}.
\]
Then, $v_{\alpha}^{\tt DRO} \geq v \alpha^{r/p +\varepsilon}$. 

The rest of the proof is focused on demonstrating the existence of such a coupling, and is completed in two steps. Suppose the random variable $\xxi_0$ is defined as the following mixture:
\[
\xxi_0 = (1-B_\delta ) \xxi + B_\delta \mv{T},
\]
where $B_\delta$ is a Bernoulli random variable with $P(B_\delta =1) =\delta$ (to be chosen later), and $\mv{T}$ is a standard multivariate-t distribution with $(p+\varepsilon)$ degrees of freedom which are independent of $\xxi$ and each other. Note that the joint density $\xxi_0$ can be expressed as
\begin{equation}\label{eqn:xi0_density}
    \tilde f(\zz_1,\zz_2) = \delta f_{\xxi}(\zz_1) f_{\mv{T}}(\zz_2)  + (1-\delta)f_{\xxi}(\zz_1) \mv{1}(\zz_1=\zz_2)
\end{equation}
\textbf{i) Show that $\xxi_0$ satisfies Assumption~$\mathscr{H}$ with $\gamma =p+\varepsilon$: }Note first that by  construction, $\mv T$ satisfies Assumption~($\mathscr{H}$) with $\gamma =(p+\varepsilon)$. Now,
\begin{align*}
    \bar F_{0,i}(tz) &=  (1-\delta) \bar F_{i}(tz)+ \delta \bar F_{\mv T,i}(tz)\\
   & \sim \delta \bar F_{\mv T,i}(tz) \\
    &\sim \delta c_i z^{-\gamma}\bar F_{\mv T}(t) \text{ since $\mv T$ satisfies Assumption~($\mathscr{H}$)}
\end{align*}
Similarly, 
\begin{align*}
f_{0}(t\zz) &= (1-\delta)f_{\xxi}(t\zz) + \delta f_{\mv T}(t\zz)\\
& = (1-\delta) \exp(- \lambda(t) \varphi^*(\zz)(1+o(1))) +\delta  \bar F_{\mv T}(t)(1+o(1)) \varphi^*_0(\zz)\\
&\sim \delta  \bar F_{\mv T}(t)\varphi^*_0(\zz) 
\end{align*}
Combining the above two displays shows that $\xxi_0$ satisfies Assumption~($\mathscr{H}$) with $\gamma = (p+\varepsilon)$.

\noindent \textbf{ii) Bound Wasserstein distance: } 
Note that 
\begin{align*}
    E_{\Pi}[\|\xxi-\xxi_0\|_p] &= \int_{(\zz_1,\zz_2)} \|\zz_1-\zz_2\|^p \, d\Pi(\zz_1,\zz_2) = 
    \delta \int_{(\zz_1,\zz_2)}\|\zz_1-\zz_2\|^p f_{\xxi}(\zz_1)f_{\mv T}(\zz_2) d\zz_1d\zz_2\\
    &\leq \delta \left(\left(\int_{\zz_1} \|\zz_1\|^p f_{\xxi
    }(\zz_1)d\zz_1 \right)^{1/p}+  \left(\int_{\zz_2} \|\zz_2\|^p f_{\mv{T}}(\zz_2) d\zz_2\right)^{1/p}\right)^p
\end{align*}
To see the last statement, recall that from Minkowski's inequality, for two functions $\|f+g\|_p\leq \|f\|_p+\|g\|_p$. Setting $f(\zz_1,\zz_2) =\zz_1$ and $g(\zz_1,\zz_2) = -\zz_2$ to get the last inequality. To demonstrate that both the integrals above are finite, we perform the polar co-ordinate transformation: $\zz =r\pphi$, and note that $d\zz = r^{d-1}drd\pphi$. Substituting this above, 
\begin{equation}\label{eqn:integrals}
    \int_{\zz_1}\|\zz_1\|^p f_{\xxi}(\zz_1) d\zz_1 = \int_{(r,\pphi)} r^{d-1+p}f_{\xxi}(r\pphi)drd\pphi \text{ and } \int_{\zz_2}\|\zz_2\|^p f_{\mv{T}}(\zz_2) d\zz_2 = \int_{(r,\pphi)} r^{d-1+p}f_{\xxi}(r\pphi)drd\pphi
\end{equation}
To demonstrate the integrals are finite, we show that their tails integrate to $0$. To this end,
note that as $r\to\infty,$ $f_{\xxi}(r\pphi)  = \exp(-r^\gamma \varphi^*(\pphi)(1+o(1)))$, uniformly over $\pphi\in \mathcal E\cap \mathcal S^{d-1}$. Now,  
\[
\int_{r>m,\pphi} r^{d-1}f_{\xxi}(r\pphi) drd\pphi \leq (2\pi)^d \int_{r>m} r^{d-1}\exp(-0.5\varphi^*r^{\gamma})dr \to 0 \text{ as } r\to\infty 
\]
where $\varphi^*= \inf_{\pphi\in \mathcal E\cap \mathcal S^{d-1}} \varphi^*(\pphi)$. Therefore, the tail of the first expression in \eqref{eqn:integrals} goes to $0$.

For the second integral in \eqref{eqn:integrals}, note that from Assumption~($\mathscr{H}$), $f_{\mv T}(r\pphi) \sim r^{-(d+p+\varepsilon)}\varphi^*_0(\pphi)$. Then, similar to the previous case:
\begin{align*}
    \int_{r>m,\pphi} r^{d-1}f_{\mv{T}}(r\pphi)drd\pphi &\leq  \varphi^*_0 (2\pi)^d \int_{r>m} r^{d-1}r^{-(d+p+\varepsilon)} dr \to  0 \text{ as }r\to\infty.
\end{align*} 
where $\varphi^* = \sup_{\pphi\in \mathcal S^{d-1}}\varphi^*_0(\pphi)$.
Putting the above two displays together, the integrals in \eqref{eqn:integrals} are bounded, and upon choosing $\delta$ sufficiently small, it is possible to make $E_{\Pi}[\|\xxi-\xxi_0\|^p]\leq \eta$. \qed

\noindent\textbf{Proof of Theorem~\ref{lem:moment_ambiguity}: }We demonstrate that for each choice of dispersion measure in the theorem statement, there exists a distribution $\mathbb Q_0$ which satisfies the conditions of Assumption~($\mathscr{H}$) with $\gamma = p+\varepsilon$, where $\varepsilon>0$ can be taken to be arbitrary. Then $v_\alpha^{\tt DRO}$ is at least as much as the cost of $\CCPalpha$, where $\xxi\sim \mathbb Q_0$, where the latter quantity grows as $v^*\alpha^{-(r/(p+\varepsilon))}$ due to Theorem~\ref{thm:solution_convergence}.

\noindent \textbf{Case i) $d =d_{\tt cov}$: }Here, let $\mathbb Q_0$ be a multivariate $t$ distribution with $2+\varepsilon$ degrees of freedom, mean $\mv{\mu}$ and covariance matrix $\mv{\sigma}$. Clearly $\mathbb Q_0\in \mathcal P$, and from Example~\ref{eg:Elliptical_densities}, $\mathbb Q_0$ satisfies Assumption ($\mathscr{H}$) with $\gamma = p+\varepsilon$. Consequently, the theorem statement holds as a result of the discussion above.

\noindent \textbf{Case ii) $d =d_{\tt ad,p}$ (or) $d =d_{\tt sd,p}$: } Once again, take $\mathbb Q_0$ to be a multivariate $t$-distribution with $p+\varepsilon$ degrees of freedom, whose mean is $\mv{\mu} = (\mu_1,\ldots,\mu_d)$, and whose marginals have a scale of $\mv{\tau} = (\tau_1\ldots,\tau_d)$. Note that here, $E[|\xi_i-\mu_i|^p]<\infty$ and that this quantity decreases with $\tau_i$. Hence, when $\tau_i$ is such that for all $i$, $E[|\xi_i-\mu_i|^p]<\sigma_i$, $\mathbb Q_0\in \mathcal P$, and the conditions of the theorem hold.
For $d= d_{\tt sd,p}$, note that $E[|\xi_i-\mu_i|^p] > \max\{E[(\xi_i-\mu_i)_+^p], E[(\mu_i-\xi_i)_+^p]\}$, and repeat the same steps as before.

\noindent \textbf{Case iii) $d =d_{\tt norm,p}$: } Observe that
\[
E\left[\|\xxi-\mv{\mu}\|^p \leq d^p \left(E[|\xi_1 -\mu_1|^p] +\ldots+ E[|\xi_d -\mu_d|^p] \right)\right]
\]
We now choose $\mv{\tau}$ in part ii) above, so that each of $E\left[|\xi_i-{\mu}_i|^p\right] \leq d^{-(p+1)}\sigma$. Then, repeating the calculation of part ii), $E\left[\|\xxi-\mv{\mu}\|^p\right] \leq \eta$, and $\mathbb Q_0$ satisfies Assumption ($\mathscr{H}$) with $\gamma=p+\varepsilon$.\qed

\subsection{Proofs from Section~\ref{sec:conv_rel}}
The following lemmas are essential to prove results from this section. Let
\[
w(u,\xx) = \begin{cases}\left[\inf\{\varphi^*(\zz) : \max_{k}\eta_k g_{k}^*(\xx,\zz) \geq u\}\right]^{-1} \quad \quad \ \ \text{ if Assumption ($\mathscr{L}$) holds}\\
    \int_{\max_k \eta_k g_k^*(\xx,\zz) \geq u} \varphi^*(\zz)d\zz \quad \quad \quad \quad \quad \ \  \quad \quad \quad \ \ \ \text{ if Assumption ($\mathscr{H}$) holds}
\end{cases}
\]
\begin{lemma}\label{lem:var_asymptotic}
    Suppose the conditions of Theorem~\ref{thm:solution_convergence} hold. Then, uniformly over compact sets, 
    $
    v_\alpha(\xx) \sim u^*(\xx)$ 
    where $u^*(\xx)  = \inf\{u: w(u,\xx) \leq 1\}$.
\end{lemma}

\noindent \textbf{Proof of Theorem~\ref{thm:convex_rel}: } 
Our proof strategy is identical to the one used in Theorem~\ref{thm:solution_convergence}. Consider the scaled optimization problem:
\begin{equation}\label{eqn:Oalpha}
    O_{\alpha}^{\tt rel}:  \min\{c(\xx): J_\alpha(\xx) \leq 0\},\text{ where }J_\alpha(\xx) = {\tt CVaR}_{1-\alpha}\left[\max_{k}\eta_k g_{k,t}(\xx,\xxi_t)\right] 
\end{equation}
Note that as in Theorem~\ref{thm:solution_convergence}, the value and solutions of the  CVaR constrained problem \eqref{eqn:cvar-approx} can be written as $s_\alpha^r$ times those of $O_\alpha^{\tt rel}$
To complete the proof, we demonstrate that the optimal value and solutions of \eqref{eqn:Oalpha} converge to $cv^*$ and $c\mathcal X^*$ respectively, in the sense defined in Theorem~\ref{thm:solution_convergence}. To this end, define the function $v_\alpha(\xx) = {\tt  VaR}_{1-\alpha}\left[\max_k \eta_k g_{k,t}(\xx,\xxi_t)\right]$ and note that
\begin{align}\label{eqn:CVaR_representation}
    J_\alpha(\xx) = v_\alpha(\xx) + \alpha^{-1}E\left[ \max_{k}\eta_k g_{k,t}(\xx,\xxi_t) - v_\alpha(\xx)\right]^+
\end{align}

\noindent \textbf{(i) Assumption~$(\mathscr{L})$ holds: } 
We demonstrate that $J_\alpha(\xx) \sim v_\alpha(\xx)$ uniformly over compact sets. 
first To this end, note that by definition for all $\varepsilon>0$ and $\xx\in \mathcal X$, $J_\alpha(\xx) \geq v_\alpha(\xx)$. We therefore prove a matching upper bound. To this end, it suffices to show that
\[
E\left[\max_k \eta_k g_{k,t}(\xx,\xxi_t) - v_\alpha(\xx)\right]^+=o(\alpha)
\]
Fix $\varepsilon>0$, and decompose the above expression as
\begin{align}\label{eqn:decompose_expc}
    E\left[\max_k \eta_k g_{k,t}(\xx,\xxi_t) - v_\alpha(\xx)\right]^+&\leq E\left[(\max_k \eta_k g_{k,t}(\xx,\xxi_t) - v_\alpha(\xx))^+ ; \psi_\alpha (\xx,\xxi_t)\leq \varepsilon\right]\nonumber \\
    &+ E\left[\max_k \eta_k g_{k,t}(\xx,\xxi_t) ; \psi_\alpha (\xx,\xxi_t)> \varepsilon\right] 
\end{align}
where 
\[
\psi_\alpha(\xx,\xxi_t) = \max_{k}\eta_k g_{k,t}(\xx,\xxi_t) - {v_\alpha(\xx)}.
\]
Note that the first term above can be bounded by $\varepsilon \alpha $. For the second term note that as a result of the continuous convergence of $g_{k,t}$
\[
\max_{k}\eta_k g_{k,t}(\xx,\zz) \to \max_{k}\eta_k g^*_k(\xx,\zz), \text{ uniformly over compact subsets of $\mathcal X$ }
\]
This suggests that whenever $\xx_t\to \xx$, 
\begin{align*}
    \max_{k}\eta_k g_{k,t}(\xx_t,\zz) &= \max_{k}\eta_k g_{k,t}(\xx_t,\hat\zz\|\zz\|)\\
    &= \max_{k}g^*_{k}(\xx,\hat \zz)\|\zz\|^{k}(1+o(1))
\end{align*}
as $\|\zz\|\to\infty$. Consequently, we have that whenever $\|\zz\|$ is large enough,
\begin{align*}
    \log  \max_{k}\eta_k g_{k,t}(\xx,\zz)  &\leq 2k\log \|\zz\| + \log \max_k g_k^*(\xx,\hat \zz)\\
    & \leq\varepsilon \|\zz\|^{\gamma+d} \text{ uniformly over $\xx$ in compact sets}
\end{align*}
Now, write the second term of \eqref{eqn:decompose_expc} as
\[
E\left[\max_k \eta_k g_{k,t}(\xx,\xxi_t) ; \psi_\alpha (\xx,\xxi_t)> \varepsilon\right]   \leq \exp\left(\log(g_\infty) \right) \mathbb P\left[\psi_{\alpha}(\xx,\xxi_t) > \varepsilon  \right]
\]
where $g_\infty = \sup_{k,\pp\in \mathcal S^{d-1}, \xx\in K} g_k^*(\xx,\pp)<\infty$. Use Proposition~\ref{prop:prob-convergence_LT} to obtain that uniformly over $\xx$ in compact sets, with $u = u^*(\xx) + \varepsilon$
\begin{align*}
    \log \mathbb P\left[\max_{k} \eta_k g_{k,t}(\xx,\xxi_t) > u\right] &\leq  \log \mathbb P\left[\max_{k}  g_{k,t}(\xx,\xxi_t) > u\right]  \text{ since }\eta_k\in(0,1)\\
    & \sim (w(\xx, u^*(\xx)+\varepsilon))^{r/\gamma} \log\alpha \\
    & \sim (1+\varepsilon)^{r/\gamma}\log \alpha 
\end{align*}
where the last step above follows upon noting that at $u^*(\xx)$, by definition, $w(\xx,u(\xx)) =1$, and that as $u$ increases, so does $w^*(\xx,u)$. Conclude from the above that, 
\[
E\left[\max_k \eta_k g_{k,t}(\xx,\xxi_t) ; \psi_\alpha (\xx,\xxi_t)> \varepsilon\right]  = O(\alpha^{1+\varepsilon}).
\]
Plug this observation in \eqref{eqn:decompose_expc} to note that 
\[
 E\left[\max_k \eta_k g_{k,t}(\xx,\xxi_t) - v_\alpha(\xx)\right]^+ \leq \alpha \varepsilon + O(\alpha^{1+\varepsilon})
\]
Since $\varepsilon$ above was arbitrary we have that 
\[
\alpha^{-1} E\left[\max_k \eta_k g_{k,t}(\xx,\xxi_t) - v_\alpha(\xx)\right] = o(1),
\]
Now apply Lemma~\ref{lem:var_asymptotic} to further infer that $C_\alpha(\xx) \sim u^*(\xx) $/ Following the proof of Theorem~\ref{thm:solution_convergence}, this  implies that any optimal solutions $\{\xx_\alpha^*: \alpha>0\}$ of the $\Dapprox$ satisfies $d(s_\alpha^{-r}\xx_\alpha^*,\mathcal X^*_0) \to 0$, where  $\mathcal X^*_0 = \arg\min\{c(\xx) : \max_k \eta_k g_k^*(\xx) \geq 0\}$. Noting that 
\[
\arg\min\{c(\xx) : \max_k \eta_k g_k^*(\xx) \geq 0\} = \arg\min \{c(\xx) : \max_k g_k^*(\xx) \geq 0\} = \mathcal Y^*
\]
concludes the proof.


\noindent \textbf{Proof of Theorem~\ref{thm:bfi-approx}: }
Define the functions 
\begin{align*}
    \tilde G_{k,t}(\xx) &=
 \frac{\log\mathbb P \left(g_{k,t}(\xx, \xxi_t) > 0\right) }{ \log \eta_k \lambda(t)} \quad\quad\quad \text{ if Assumption~$(\mathscr{L})$ holds} \\
  &=\frac{\mathbb P\left( g_{k,t}(\xx, \xxi_t) > 0 \right)}{\eta_k  \lambda(t)} \quad\quad \ \ \ \ \ \quad\text{ if Assumption~$(\mathscr{H})$ holds}.
\end{align*}

Similar to the proof of Theorem~\ref{thm:solution_convergence}, note that given $\{\eta_k:k\leq K \}$, \eqref{eqn:bfi-approx} can be rewritten as 
\[
t^r \min\{c(\xx/t^r): I_{k,t}(\xx/t^r) \leq 1\, \forall k\in[K]\}, \text{ where } I_{k,t}(\xx) =\begin{cases}
    1/ G_{k,t}(\xx) \quad \text{ if Assumption~($\mathcal L$) holds and}\\
    G_{k,t}(\xx)\quad \ \ \ \ \text{ if Assumption~($\mathcal H$) holds and}
    \end{cases}
\]

Note that following the proof of Theorem~\ref{thm:solution_convergence}, the following continuous convergence holds:
\[
I_{k,t} \to I_{0,k} \text{ where } I_{0,k} = \begin{cases}
\left[\inf\{\varphi^*(\zz) : g^\ast_k(\xx,\zz) \geq 0\}\right]^{-1}  \quad \quad \quad  \ \ & \text{ if } \xxi \text{ satisfies Assumption~(${\mathscr L}$)}\\
\eta_k^{-1} \int_{g^\ast_k(\xx,\zz) \geq 0}  \varphi^*(\zz) d\zz  \quad & \text{ if } \xxi \text{ satisfies Assumption~(${\mathscr H}$)},  
\end{cases}
\]
Note that given two sequences $\{a_{n,1},\ldots a_{n,K}:n\geq 1\}$ with $a_{n,k}\to a_k$ for all $k$, $\max\{a_{n,k}: k\leq K\}\to \max\{a_k:k\leq K\}$. Set $I_{k,t}(\xx_t) :=a_{k,t}$ and let $\xx_t\to\xx$. Therefore, with $I_{k,t}(\xx_t)\to I_{0,k}(\xx)$for all $k$, one has that $\max_k I_{k,t} \to \max_k I_{0,k}$ continuously. Observe that the latter evaluates to 
\begin{equation}\label{eqn:I_0k}
    \max_k I_{0,k}(\xx) = \begin{cases}
\left[\inf\{\varphi^*(\zz) : \max_k g^\ast_k(\xx,\zz) \geq 0\}\right]^{-1}  \quad \quad \quad  \ \ & \text{ if } \xxi \text{ satisfies Assumption~(${\mathscr L}$)}\\
\max_k \eta_k^{-1} \int_{g^\ast_k(\xx,\zz) \geq 0}  \varphi^*(\zz) d\zz  \quad & \text{ if } \xxi \text{ satisfies Assumption~(${\mathscr H}$)},  
\end{cases}
\end{equation} 
Define the optimization problems:
\[
O_\alpha^{\tt bf}: \min\{c(\pp) : \max_kI_{k,t}(\pp) \leq 1\}
\]
Then, the optimal value and set of solutions of $O_\alpha^{\tt bf}$ converge in the sense given by Theorem~\ref{thm:bfi-approx} to those of 
\[
\min\{c(\pp): \max_kI_{0,k}(\xx)\leq 1\} \text{ where } \max_{k} I_{k,0} \text{ is defined in \eqref{eqn:I_0k}}
\]
This suggests the results of parts (a) and (b). \qed

\noindent \textbf{Proof of Lemma~\ref{lem:line-search-properties}:} 
\textbf{a) }Note that the cost $c(\xx^{\tt app}_\alpha)$ is positive and homogeneous. Then,
$c(t^r \xx_{\alpha}^{\tt app}) = t^rc(\xx_{\alpha}^{\tt app}) $. 
Let us now consider two cases: first when $r>0$, the cost is positive,
and consequently, $c(\xx(t))$ is an increasing function of $t$. Next if $r<0$ then $c<0$, so as $t$ increases, once again $c(t^r\xx_\alpha^{\tt app})$ is increasing in $t$.

\noindent \textbf{b)} Note that $p(\xx_{\alpha}^{\tt app})\leq \alpha$, since $\xx_{\alpha}^{\tt app}$ is feasible to the safe approximation. Consequently, the infimisation problem defining $t_\alpha$ has a feasible solution. Further, by definition, for all $t\geq t_\alpha$, one has that $p(t^r\xx_{\alpha}^{\tt app}) \leq \alpha$,  or that $\xx(t)$ is feasible to $\CCPalpha$.

\noindent \textbf{c) }Note that $p(\cdot)$ is a continuous function, and by definition, $p(\xx_{\alpha}^{\tt app}) < \alpha$. Owing to continuity, there then exists a $t_\alpha^{\prime}<1$ such that $p(\xx(t))< \alpha$ for all $t\in (t_\alpha^\prime,1]$. Therefore, $t_\alpha = \inf\{t\in [0,1]: p(\xx(t)) \leq \alpha\} <t_\alpha^\prime <1$.
From part \textbf{a)} it now follows that $c(\xx_\alpha^\prime) < c(\xx_\alpha^{\tt app})$.\qed

\noindent \textbf{Proof of Proposition~\ref{prop:Approx_alg_guarantees}: }Note that for any $s\in(0,1]$, 
\begin{align*}
    p(s\xx_{\alpha}^{\tt app}) &=\mathbb P\left(\max_k g_k(s^r\xx_{\alpha}^{\tt app},\xxi) >0 \right)\\
    & = \mathbb P\left(g_t(s^rt^{-r}\xx_\alpha^{\tt app}) >0 \right) \\
    & = \mathbb P\left(g_t(s^r\yy_t^*(1+o(1))) >0 \right)  \text{ for a suitable $\yy_t^*\in c\mathcal X$}
\end{align*}
The last statement above follows from Theorems~\ref{thm:convex_rel}-\ref{thm:bfi-approx}, where it is demonstrated that any sequence of solutions  $\xx_\alpha^{\tt app}$  of the relaxed problem satisfy
$d(s_\alpha^{-r} \xx_{\alpha}^{\tt app},c\mathcal X^*) \to0 $ and therefore, $ t^r\xx_{\alpha}^{\tt app} = \yy_t^*+o(1)$ for some $\yy_t^*\in c\mathcal X$. Using the above observation, note that the minimization problem in Algorithm~\ref{algo:Vanish_Regret} can be re-cast as
\begin{equation}\label{eqn:t-alpha-value}
    t_\alpha = \inf\left\{s\in[0,1]: I_t(s\yy_t^*(1+(o(1))) \leq 1\right\}
\end{equation}
where $I_t$ is as defined in the proof of Theorem~\ref{thm:solution_convergence}. 
Define the function $w_t(s) =I_t(s\yy_t^*(1+(o(1))) $ and $w(s) = (cs)^{
-\gamma/r}$. 
Observe that with $\yy_t^*\in c\mathcal X^*$, one has that $I(\yy_t^*) = 1$.
Then, using the homogeneity of $I$ from Lemma~\ref{lem:I_homogenity}, note that whenever $\yy_t^*\in c\mathcal X$, $I(s\yy_t^*) =(cs)^{-\gamma/r} = w(s)$.
\[
|w_t(s_t) - w(s) |\leq |I_t(s_t\yy_t^*(1+o(1)) - I(s_t \yy_t^*(1+o(1)))| + |I(s_t \yy_t^*(1+o(1))) - w(s)| \leq \varepsilon
\]
for all $t$ sufficiently large. To see this, note that the convergence of $I_t\to\ I$ is uniform over compacts sets using Propositions~\ref{prop:prob-convergence_LT}-\ref{prop:prob-convergence_HT}. The first term above can therefore be made as small as desired. Further note that since $I$ is continuous (and therefore uniformly continuous on compact sets), the second term can also be made arbitrarily small. 
This suggests that $w_t\to w$ continuously, and following the proof of Theorem~\ref{thm:solution_convergence}, we have that the optimal value of \eqref{eqn:t-alpha-value} converges to 
\[
\inf\{s: (cs)^{-r/\gamma} \leq 1\} = 1/c
\]
Now, write $\xx_\alpha^\prime = t_\alpha^r\xx_{\alpha}^{\tt app} = c(1+o(1))\xx_{\alpha}^{\tt app}$. Using the homogeneity of the cost function $c(\cdot)$, we have that 
\[
\frac{c(\xx_\alpha^\prime) - v_\alpha^*}{v_{\alpha}^*} = \frac{c^{-1}\times c v^*s_\alpha^{r}(1+o(1)) - v^* s_\alpha^r (1+o(1)) }{v^* s_\alpha^r (1+o(1))} \to 0
\]
as $\alpha\to 0.$ Finally, the decrease in suboptimality is given  by
\[
\lim_{\alpha\to 0} \frac{cv_\alpha^*s_\alpha(1+o(1))}{v^*s_\alpha^r} = c
\]
\qed
\subsection{Proofs from Section~\ref{sec:Algo}}
The following consequences are immediate from Lemma~\ref{lem:RVs} - for any $M>0$, the statements below hold uniformly over  $t\in[1,M]$:
\begin{align}\label{eqn:get_rvs_from_s}
     s_{\alpha_0^{t}} &= s_{\alpha_0} t^{1/\gamma}(1+o(1)) \text{ if Assumption~($\mathscr{L}$) holds}\nonumber\\
     s_{\alpha_0/t} &= s_{\alpha_0} t^{1/\gamma}(1+o(1)) \text{ if Assumption~($\mathscr{H}$) holds}
\end{align}

The following lemmas are required to proceed
\begin{lemma}\label{lem:saa_sol_convergence}
    Suppose the conditions of Theorem~\ref{thm:WPE-data-driven} hold. Then, as $N\to \infty$, $d(s_{\alpha_0}^{-r} {\hat \xx}^{(N)}_{\alpha_0} ,\mathcal Y^*) \to 0$ almost surely.
\end{lemma}
\begin{lemma}\label{lem:I_value}
    Let $v^* = \min\{c(\xx) : I(\xx)\leq 1\}$, where $I(\cdot)$ is as defined in Theorem~\ref{thm:solution_convergence}. Then $\min\{I(\xx) : c(\xx) \leq v^*\} =1$.
\end{lemma}


\noindent \textbf{Proof of Proposition~\ref{prop:solution-set-convergence}: } Now, recall that if $\mathcal Y^*=\{\xx^*\}$ is a singleton set,  then in addition $\xx_{\alpha}^* =\xx^*s_\alpha^{r}(1+o(1)). $ Note that this implies $\|\xx_\alpha^*\| = s_\alpha^r(1+o(1)) \|\xx|$
Now, suppose that Assumption ($\mathscr{L}$) holds. Then, 
\[
\xx^*_{\alpha^t} = \xx^* s_{\alpha^t}^r(1+o(1)) = \xx^* t^{r/\gamma}(1+o(1)),  \text{ from \eqref{eqn:get_rvs_from_s}}.
\]
From the definition of set distances, recall that
\[
d(t^{r/\gamma} \xx^*_{\alpha_0} ,\mathcal X(\alpha^t_0)) = \inf\{d(t^{r/\gamma} \xx^*_{\alpha_0}, \yy): \yy\in  \mathcal X(\alpha^t_0)\}
\]
Let $\yy_\alpha^*$ be the arg-min of the above problem. Consider the following:
\begin{align*}
   \|t^{r/\gamma}\xx_{\alpha_0} - \xx_{\alpha_t}^*\| &=   \|t^{r/\gamma}\xx_{\alpha_0} - \xx^* s_{\alpha_0^t}^r(1+o(1))\| \\
  &=\|\xx^* s_{\alpha_0}^{r} t^{r/\gamma} - \xx^* s_{\alpha_0}^r t^{r/\gamma}(1+o(1))\| \quad \text{ from \eqref{eqn:get_rvs_from_s}}\\
  &= o(s_{\alpha_0}^{r}) \text{ since $t$ and $\xx^*$ are fixed}.
\end{align*}
Then, we have that $d(t^{r/\gamma} \xx^*_{\alpha_0} ,\mathcal X(\alpha^t_0)) = o(\|\xx_{\alpha_0}\|)$. Next, suppose that Assumption ($\mathscr{H}$) holds. Here, we have 
\begin{align*}
   \|t^{r/\gamma}\xx_{\alpha_0} - \xx_{\alpha_t}^*\| &=   \|t^{r/\gamma}\xx_{\alpha_0} - \xx^* s_{\alpha_0/t}^r(1+o(1))\| \\
  &=\|\xx^* s_{\alpha_0}^{r} t^{r/\gamma} - \xx^* s_{\alpha_0}^r t^{r/\gamma}(1+o(1))\| \quad \text{ from \eqref{eqn:get_rvs_from_s}}\\
  &= o(s_{\alpha_0}^{r}) \text{ since $t$ and $\xx^*$ are fixed}.
\end{align*}
and once again the conclusions of the proposition hold\qed

\noindent \textbf{Proof of Theorem~\ref{thm:WPE-data-driven}: }
To prove almost sure weak pareto efficiencies of the trajectories, it is sufficient to demonstrate that \eqref{eq:WPE} holds on some set $\Omega_0$ which has a probability of $1$. To this end, let's first suppose that $r>0$.  Let $\tilde \xx_N = s_{\alpha_0}^{-r}  \xx_{\alpha_0}^{(N)}$ and let $\Omega_0$ be the set from Lemma~\ref{lem:saa_sol_convergence} on which $d(\tilde \xx_N,\mathcal Y^*) \to 0$. Observe that any $\hat \xx^{(N)} \in \bar \xx_N$ has the representation $\hat \xx^{(N)} = s_{\bar\alpha} \tilde \xx_N$ where $\bar \alpha = \bar F(t_Ns_{\alpha_{0}})$. 

Note the following set of implications:
\begin{align}\label{eqn:cost_implications}
    c(\xx) <(1-\varepsilon)c(\bar{\xx}^{(N)}) &\overset{\textbf{(a)}}{\implies} (1-\varepsilon) s_{\bar\alpha}^r c(\tilde{\xx}^{(N)}) \overset{\textbf{(b)}}{\implies}  c\left(s_{\bar\alpha}^{-r}\frac{\xx}{1-\varepsilon}\right) \leq c(\tilde{\xx}_{N})
\end{align}
Here, implication \textbf{(a)} above follows using the homogeneity of $c(\cdot)$ and the definition of $\tilde{\xx}_{N}$. 
Now, note that for every $\kappa>0$ and  $\omega\in\Omega_0$, there exists a large enough $n_0(\omega)$ such that for all $n\geq n_0(\omega)$, $d(\tilde {\xx}^{(N)},\mathcal Y^*) \leq \kappa$. Since $c(\cdot)$ is continuous and $\mathcal Y^*$ is compact, given any $\varepsilon>0$ we can choose a $\kappa$ so small that $c(\tilde {\xx}^{(N)}) < (1+\varepsilon)c(\mathcal Y^*)$. Then, we have from \eqref{eqn:cost_implications} that
\[
c\left(s_{\bar \alpha}^{-r}\frac{\xx}{1-\varepsilon^2}\right) \leq c(\mathcal Y^*) \text{ for all } n\geq n_0(\omega)
\]
This now implies that 
\[
I\left(s_{\bar \alpha}^{-r}\frac{\xx}{1-\varepsilon^2}\right)  \geq 1 \text{ for all } n\geq n_0(\omega)
\]
since the minimum value of the function $c(\cdot)$ over $I(\xx) \leq 1$ is $c(\mathcal Y^*)$. Using the homogeneity of $I(\cdot)$ from Lemma~\ref{lem:I_homogenity}, we further have that 
$I(\xx s_{\bar\alpha}^{-r}) > (1-\varepsilon^2)^{-\gamma/r}$. 
Since $I(\cdot)$ is homogeneous with order $-\gamma/r$ and $r>0$, the above also implies that $\xx s_{\bar \alpha}^{-r}$ lies in a compact set. If $\bar t = s_{\bar\alpha}$,
then as $\bar \alpha \to 0$, $I_{\bar t}\to I$ uniformly over compact sets. Then, since $I_{\bar t}(\xx) \geq I(\xx) + \varepsilon $ uniformly over $\xx$ in compact sets. This implies that 
\[
I_{\bar t}(s_{\bar \alpha}^{-r}\xx) > (1-\varepsilon^2)^{-\gamma/r} -\varepsilon \quad \text{ for all $n$ sufficiently large}
\]

Now consider two cases:

\noindent \textbf{i) Assumption~($\mathscr{L}$) holds: } Recall that 
\[
I_{\bar t}(s_{\bar \alpha}^{-r}\xx) = -\frac{\log \bar F(s_{\bar\alpha}) }{\log p(\xx)} \implies \log p(\xx) \geq  \frac{\log \bar F(s_{\bar \alpha})}{(1-\varepsilon^2)^{-\gamma/r} -\varepsilon} \text{ for all } n>n_0
\]
Note now that $\log p(\hat \xx^{(N)})= \log  p(s_{\bar\alpha} \tilde \xx_N)$ where $d(\tilde \xx_N ,\mathcal Y^*) \to 0$ as $N\to\infty$. Consequently, we have that the tail sequence, $\{\tilde \xx_N: N>n_0\} $ lies in a compact set. Then, using the uniform convergence of log-probabilities from Proposition~\ref{prop:prob-convergence_LT},
\[
\left\vert\frac{-\log  p(s_{\bar\alpha}^{-r} \tilde \xx_N)}{\log \bar F(s_{\bar \alpha})} - I(\tilde \xx_N) \right\vert \leq \varepsilon \quad \text{ for all $n>n_1$ }. 
\]
Combining the above two displays, for all $n\geq \max\{n_0,n_1\}$
\[
\log p(\xx) \geq (1-\delta) \log p(\hat \xx^{(N)})
\]
where $\delta >0$. Note that by setting $t_N$ to be a sufficiently large constant, we can ensure the condition $p(\hat \xx^{(N)}) \sim aN^{-b}$ for $a>0$, $b\geq 1$. This concludes the proof in the light tailed case. 

\noindent \textbf{ii) Assumption~($\mathscr{H}$) holds:}
Here, we instead have that
\[
I_{\bar t}(s_{\bar \alpha}^{-r}\xx) = \frac{p(\xx)} {\bar F(s_{\bar\alpha}) }\implies p(\xx) \geq  \frac{ \bar F(s_{\bar \alpha})}{(1-\varepsilon^2)^{-\gamma/r} -\varepsilon} \text{ for all } n>n_0
\]
Note that from Proposition~\ref{prop:prob-convergence_HT}, that 
\[
\left\vert\frac{p(s_{\bar\alpha}^{-r} \tilde \xx_N)}{\bar F(s_{\bar \alpha})} - I(\tilde \xx_N) \right\vert \leq \varepsilon \quad \text{ for all $n>n_1$ }. 
\]
Combining the above two displays suggests that 
$p(\xx) \geq (1-\delta) p(\hat \xx^{(N)})$ for all $n\geq \max\{n_0,n_1\}$, which further suggests conclusion of the theorem.\qed

\noindent\textbf{Proof of Corollary~\ref{cor:P-Model-data-driven}: } First suppose that Assumption~($\mathscr{L}$) holds. 
Note the following for a arbitrary $\omega\in \Omega_0$.
\begin{align}\label{eqn:log_pn}
    \log p(\hat \xx^{(N)}) &= \log p(t^r s_{\alpha_0}^r \tilde \xx_N) = \log p(s_{\bar \alpha}^r \tilde \xx_N)\nonumber\\
    & = \lambda(s_{\bar \alpha}) I(\tilde \xx_N) (1+o(1)) = \lambda(s_{\bar \alpha}) I(\mathcal Y^*)(1+o(1)) \nonumber\\
    & = \lambda(s_{\bar \alpha}(1+o(1))) \text{ as $N\to\infty$}.
\end{align}
Next, note that from Theorem~\ref{thm:solution_convergence}, $t^rv_{\alpha_0}^* = (ts_{\alpha_0}^r v^*(1+o(1)))$. Then, write $\nu^*(t^rv_{\alpha_0}^*) $ as 
\[
\min_{\xx\in \mathcal X} \log p(\xx) \quad  \text{ s.t } \quad  c(\xx) \leq s_{\bar \alpha}^r v^*(1+o(1)).
\]
Use the homogeneity of $c(\cdot)$ and change variables to $\yy = s_{\bar\alpha}^{-r}(1+o(1)) \xx$ to reduce the above problem to
\[
\min_{\yy\in \mathcal X} \log p(s_{\bar\alpha}^{-r}\yy)  : c(\yy)\leq v^*
\]
Now, recall that from Proposition~\ref{prop:solution-set-convergence}, $\log p(s_{\alpha_0}^{-r}\yy) \sim -\lambda(s_{\alpha_0})[I(\yy)]^{-1}$. Therefore, an application of \cite{bonnans2013perturbation}, Proposition 4.4 yields that 
\[
\nu^*(t^rv_{\alpha_0}^*) \sim -\lambda(s_{\alpha_0}) \  \inf\{I(\xx) : c(\xx ) \leq 1\} \sim -\lambda(s_{\alpha_0})
\]
where the last asymptotic above follows from Lemma~\ref{lem:I_value}. Combining the above display with \eqref{eqn:log_pn} completes the proof for the light tailed case. 

If instead Assumption~($\mathscr{H}$) holds. Then, instead of \eqref{eqn:log_pn}, we have that $p(\hat \xx^{(N)}) =\lambda(s_{\bar \alpha})(1+o(1))$. Further, one has that $\nu^*(B)= \log (\mu^*(B))$ where  
\[
\mu^*(B) = \inf\{p(\xx) : c(\xx)\leq B\}
\]
Repeating the above steps,
$\mu^*(t^r v_{\alpha_0}^*) \sim \lambda(s_{\bar\alpha})$. Consequently, $\log p(\hat \xx^{(N)}) \sim \mu^*(t^r v_{\alpha_0}^*)$ as $N\to\infty$. Since $\omega\in \Omega_0$ was arbitrary, all the above asymptotics hold almost surely.\qed
\newpage

\newpage

\section{Proofs of Intermediate Results}
\noindent\textbf{Proof of Proposition~\ref{prop:ldp}: }  A sufficient
condition (see \cite{dembo2009large}, Theorem 4.1.11) to verify the existence
of LDP is to show that $\{\xxi_t:t\geq 1\}$ is exponentially tight, and that for all $\xx$, 
\begin{equation}\label{eqn:LDP-Equivalent}
-\varphi^*(\xx) = \inf_{\delta>0} \limsup_{t\to\infty} \frac{1}{t}\log\Prob\left( {\xxi}_{t} \in B_{\delta}(\xx)\right) = \inf_{\delta>0} \liminf_{t\to\infty} \frac{1}{t}\log\Prob\left( {\xxi_t} \in B_{\delta}(\xx)\right).
\end{equation}
Fix any $\varepsilon, M \in (0,\infty)$ and $\xx \in (0,M)^d.$ Since
$f_{\xxi}(\yy) = \exp(-\varphi(\yy)),$
\begin{align*}
  \Prob(\xxi_t \in B_{\delta}(\xx))
  = \int_{\xxi/g(t)\in B_{\delta}(\xx)}\exp(-\varphi(\yy))d\yy
    = [g(t)]^d \int_{\zz\in B_{\delta} (\xx)}  \e^{-\varphi(g(t)\zz)} d\zz.
\end{align*}

\noindent \textbf{Case (i): $\xx\not\in \mathcal E$: } Since $\mathcal E$ is closed, there exists a $\delta_1$ such that for all $\delta<\delta_1$, $B_{\delta}(\xx) \cap \mathcal E  = \emptyset$. Further,  since $\mathcal E$ is a cone, if $\xxi$ is supported on $\mathcal E$, so is $\xxi_t$ for all $t$. 
Consequently, for all $\delta<\delta_1$,  we have that $\mathbb P[\xxi_t\in   B_{\delta}(\xx)] = 0$, and \eqref{eqn:LDP-Equivalent} holds.

\noindent \textbf{Case (ii): $\xx\in\mathcal E$: } Recall that Assumption ($\mathscr{L}$) implies the following uniform
convergence over compact subsets of $\mathcal E$ not containing the origin: 
\begin{equation}
  \label{eqn:UC-expo}
n^{-1}\varphi(g(n)\xx) \xrightarrow{n\to\infty} \varphi^*(\xx).
\end{equation}
Due to this local uniform convergence and the continuity of $\varphi^*$, there
exist $\delta_0,t_0 \in (0,\infty)$ such that,
\begin{align*}
  \left\vert \frac{\varphi(g(t)\zz)}{t} - \varphi^*(\xx)  \right\vert \leq  \left\vert \frac{\varphi(g(t)\zz)}{t} - \varphi^*(\zz)\right\vert + \left\vert \varphi^*(\zz) - \varphi^*(\xx)\right\vert \leq \varepsilon/2, \textrm{ for all } \zz\in B_{\delta}(\xx)\cap \mathcal E
\end{align*}
whenever $t > t_0,$ $\delta < \delta_0$.
Thus, given $\varepsilon, M$ and $\xx \in (0,M)^d,$ there exist
$\delta_0,t_0 \in (0,\infty)$ such that for all $t>t_0$ and
$\delta \in (0, \delta_0),$
\begin{equation}
  \label{eqn:Y-Bounds}
  \exp \left(-t(\varphi^*(\xx) +\varepsilon)\right) \leq f_{\xxi}(g(t)\zz) \leq
  \exp\left({-t(\varphi^*(\xx) -\varepsilon)}\right), \text{ uniformly over
    $\zz \in B_{\delta}(\xx) \cap \mathcal E ;$}
\end{equation}
Then 
\[
[g(t)]^d \mathrm{Vol}(B_{\delta}(\xx)\cap \mathcal E) \exp({-t(\varphi^*(\xx) + \epsilon)}) \leq
\Prob(\xxi_t \in B_{\delta}(\xx)) \leq [g(t)]^d
\mathrm{Vol}(B_{\delta}(\xx)) \exp({-t(\varphi^*(\xx) -\epsilon)}).
\]
Since $\Prob(\xxi_t \in B_{\delta}(\xx)) $ is increasing in $\delta$
and these bounds hold for any $\delta < \delta_0,$
\[ - \varphi^*(\xx) - \epsilon \leq \inf_{\delta > 0}\liminf_{t\to\infty} \frac{1}{t} \log
  P(\xxi_t \in B_{\delta}(\xx)) \leq \inf_{\delta > 0} \limsup_{t\to\infty} \frac{1}{t} \log
  P(\xxi_t \in B_{\delta}(\xx)) \leq - \varphi^*(\xx) +\epsilon. 
\]
The last statement follows since 
$g\in \RV(1/\gamma)$, $\log g(t) = o(t)$ (see \cite{deHaan}, Proposition B.1.9).
Since the choices $\varepsilon, M \in (0,\infty)$ are arbitrary,
(\ref{eqn:LDP-Equivalent}) holds.

\noindent \textbf{Step 2: Prove exponential tightess: }To conclude, we demonstrate that  given $M>0$, there exists $\alpha=\alpha(M)$  such that $\limsup_n n^{-1}\log \mathbb P(\xxi_n\in B_\alpha^c) <-M$. Note that the hazard function of the radius $R=\|\xxi\|$ is
\begin{equation}\label{eqn:R_density}
    \Lambda_R(t) = -\log \left\{\int_{r\geq t,\pphi} r^{d-1}\exp(-\varphi(r\pphi)) drd\pphi\right\}
\end{equation}
Under Assumption ($\mathscr{L}$),
as $n\to\infty$, $\varphi(g(n)\xx) = n\varphi^*(\xx)(1+\varepsilon(n,\xx))$, where for any compact set $\mathcal K\subseteq \mathcal E$,
$\sup_{t\geq n,\xx\in \mathcal K}|\varepsilon(n,\xx)| \to 0$ as $n\to\infty$.
with $n=\lambda(r)$ the integral in  \eqref{eqn:R_density} becomes 
\begin{equation*}
\int_{r\geq t,\pphi} r^{d-1}\exp(-\lambda(r) \varphi^*(\pphi)(1+\varepsilon(r,\pphi))) drd\pphi = 
\int_{p\geq g^{\leftarrow}(t),\pphi} [\lambda(p)]^{d-1}\exp(-p \varphi^*(\pphi)(1+\varepsilon(\lambda(p),\pphi))) dpd\pphi    
\end{equation*}
Recall that from Lemma~\ref{lem:RVs}, $\lambda \in \RV(\gamma)$. Therefore, as $t\to\infty$, $\log \lambda(t) = o(t)$  (see \cite{deHaan}, Proposition B.1.9). In addition for every $\delta_0>0$, there exists a $n_0$ such that for $n>n_0$, $|\sup_{\pphi,n\geq n_0}  \varepsilon(n,\pphi)| \leq \delta_0$. The previous two statements imply that for all sufficiently large $t$,  
\begin{align*}
    \int_{p\geq\lambda(t),\pphi} [\lambda(p)]^{d-1}\exp(-p \varphi^*(\pphi)(1+\varepsilon(\lambda(p),\pphi))) dpd\pphi    &\leq \int_{p\geq \lambda(t)} \exp(-p \varphi^*(\pphi)(1-\delta_0))dpd\pphi\\
    &\leq {\tt Leb}(S^{d-1}) \exp(-\lambda(t) \inf_{\pphi\in \mathcal S^{d-1}\cap \mathcal E}\varphi^*(\pp))
\end{align*}
Consequently, 
\begin{equation}\label{eqn:LB_Lambda_R}
    \liminf_{t\to\infty} \frac{\Lambda_R(t)}{\lambda(t)} \geq \inf_{\pphi\in \mathcal S^{d-1}\cap \mathcal E} \varphi^*(\pphi)(1-\delta_0).
\end{equation}

To establish a lower bound, observe that the function $\varphi^*$ is continuous on $\mathcal E$. Let $\pphi_0$ be any point on $\mathcal S^{d-1}\cap \mathcal E$ where $\varphi^*$ attains a minimum. Then, fixing a $\delta_0$, there exists a neighbourhood of $\pphi_0$, call it $N_\delta(\pphi_0)$ 
such that $|\varphi^*(\pphi)-\varphi^*(\pphi_0)|\leq \delta_0$ for all $\pphi\in N_\delta(\pphi_0)$. Now, for all sufficiently large $t$, the integral 
\begin{align*}
    \int_{p\geq g^{\leftarrow}(t),\pphi} [\lambda(p)]^{d-1}\exp(-p \varphi^*(\pphi)(1+\varepsilon(\lambda(p),\pphi))) dpd\pphi    
    \end{align*}
    is lower bounded by
    \begin{align*}
    & \int_{p\geq \lambda(t), N_\delta(\pphi_0)} \exp(-p \varphi^*(\pphi_0)(1+\delta_0) dpd\pphi  \geq {\tt Leb}(N_\delta(\pphi_0)) \exp(-\lambda(t) (\delta_0+ (1+\delta_0) \inf_{\pphi\in \mathcal {S}^{d-1}\cap \mathcal E}\varphi^*(\pphi))
\end{align*}
Consequently, 
\begin{equation}\label{eqn:UB_Lambda_R}
    \liminf_{t\to\infty} \frac{\Lambda_R(t)}{\lambda(t)} \geq \inf_{\pphi\in \mathcal {S}^{d-1}\cap \mathcal E} \varphi^*(\pphi)(1+\delta_0) +\delta_0
\end{equation}
Since $\delta_0$ in \eqref{eqn:LB_Lambda_R}-\eqref{eqn:UB_Lambda_R} is arbitrary, conclude that $\Lambda_R(t) \sim \lambda(t) c_{\varphi^*}$,
where $\varphi^* = \inf_{\pphi\in \mathcal {S}^{d-1}\cap \mathcal E} \varphi^*(\pphi)$. To conclude the proof, note that 
\[
\log \mathbb P[\xxi/g(n) \in B_{\alpha}^c] = \log \mathbb P[R \geq g(n)\alpha]  \sim  - c_{\varphi^*} \lambda(g(n)\alpha)
\]
Since $\lambda\in \RV(\gamma)$, one has that $g\in \RV(1/\gamma)$ (see \cite{deHaan}, Proposition B.1.9). Now, further simplify the above expression as $n^{-1}\log \mathbb P[R \geq g(n)\alpha] \sim - c_{\varphi^*} \alpha^{1/\gamma} $
and thus,
\[
\limsup_{n\to\infty} \log n^{-1}\mathbb P[\xxi/g(n) \in B_{\alpha}^c] \leq - c_{\varphi^*} \alpha^{1/\gamma}.
\]
Set $\alpha = (M/c_{\varphi^*})^\gamma$ above to complete the proof.
\qed

\noindent \textbf{Proof of Proposition~\ref{thm:HTT_RV}: }
Note that since $0\not\in A$, $A\subset B_{\varepsilon}^c$ for some $\varepsilon>0$. Then \cite{Resnickbook}, Proposition 5.20, implies that
\[
\int_{\zz\in A}t^{d}f_{\xxi}(t\zz)d\zz \sim \bar \lambda(t) \nu(A)
\]
Next, observe that $\xxi/t$ has a density equal to $f_{\xxi/t}(\zz) = t^{d}f_{\xxi}(t\zz)$. Substituting this expression above implies that $\mathbb P(\xxi/t \in A )\sim \lambda(t)\nu(A)$.\qed

\noindent \textbf{Proof of Lemma~\ref{lem:RVs}: }First suppose Assumption~($\mathscr{L}$)  holds. Note that 
\begin{align*}
    -\log f((tc)\zz) \lambda(tc) \varphi^*(\zz) \quad \text{ and } -\log f((tc)\zz) \sim \lambda(t) \varphi^*(c\zz).
\end{align*}
Taking limits of both of the above implies that
\[
\frac{\varphi^*(c\zz)}{\varphi^*(\zz)} = \lim_{t\to\infty}\frac{\lambda(ct)}{\lambda(t)} = c^{\gamma}
\]
 Next suppose Assumption~($\mathscr{H}$) holds. Then, 
\begin{align*}
    f((tc)\zz) \sim (ct)^{-d}\bar F(tc)  \varphi^*(\zz) \quad \text{ and } f((tc)\zz) \sim \bar t^{-d} F(t) \varphi^*(c\zz).
\end{align*}
Once again, taking limits,
\[
\frac{\varphi^*(c\zz)}{\varphi^*(\zz)} = \lim_{t\to\infty}\frac{(ct)^{-d}\bar F(ct)}{t^{-d}\bar F(t)} = c^{-(d+\gamma)}
\]
\qed

\noindent \textbf{Proof of Lemma~\ref{lem:I_homogenity}: }
Note that under Assumption~\ref{assume:constraint-functions-cvar}, we have that for every $s>0$
\begin{align}\label{eqn:g_homogenity}
    g(s^r\xx,s\zz) &= \lim_{t\to\infty} t^{-k}\max_k g_{k,t}((ts)^r\xx,(ts)\zz) \nonumber\\
    &=(ts)^{-k} s^k\lim_{t\to\infty}\max_{k}g_{k,t}((ts)^r\xx,(ts)\zz) = s^{-\rho}g^*(\xx,\zz) 
\end{align}
We prove Lemma~\ref{lem:I_homogenity} in two parts. First suppose that Assumption~($\mathscr{L}$) holds, and 
\[
I(\xx) = \inf\left\{\varphi^*(\zz) : g^*(\xx,\zz) \geq 0\right\}^{-1}
\]
Then
\begin{align*}
 I(c^r\xx) &=    \inf\left\{\varphi^*(\zz) : g^*(c\xx,\zz) \geq 0\right\}^{-1} = c^{-\gamma} \inf\left\{\varphi^*(\zz/c) : g^*(c^r\xx,c(\zz/c)) \geq 0\right\}^{-1} \text{ using Lemma~\ref{lem:RVs}}\\
 & = c^{-\gamma} \inf\left\{\varphi^*(\pp) : g^*(c^r\xx,c\pp) \geq 0\right\}^{-1} \text{ setting $\zz/c=\pp$}\\
 &=c^{-\gamma} \inf\{\varphi^*(\pp) : g^*(\xx,\pp) \geq 0\} \text{ from \eqref{eqn:g_homogenity}}
\end{align*}
Note that the last infimum above equals $I(\xx)$. Then, $I(c\xx) = c^{-\gamma/r}I(\xx)$. Next, suppose that Assumption~($\mathscr{H}$) holds. Then, recall that 
\[
I(\xx)  = \int_{g^*(\xx,\zz) \geq 0}\varphi^*(\zz) d\zz
\]
Then,
\begin{align*}
I(c^r\xx) &= \int_{g(c^r\xx,\zz) \geq 0}\varphi^*(\zz) d\zz  = \int_{g(c^r\xx,c\pp) \geq 0} c^d\varphi^*(c\pp)d\pp \quad \text{ setting } \zz =c\pp\\
& = \int_{g(\xx,\pp)\geq 0} c^{-(\gamma+d)+d}\varphi^*(\pp) = c^{-\gamma}I(\xx)
\end{align*}
The result for heavy tails now follows.\qed

\noindent{\textbf{Proof of Lemma~\ref{lem:level_sets_containment}:} }
Note that under 
Assumption~\ref{assume:non-hom-safeset}, whenever $\xx_n\to \xx$, the functions $g_n(\xx_n,\cdot)\xrightarrow{e}g^*(\xx,\cdot)$. The first containment of \eqref{eqn:level_sets_containment} follows from Corollary 1, \cite{deo2023achieving} upon setting (i) $\alpha=0$ there, and (ii) replacing the ambient space to $\R^d$ from $\R^d_{++}$. For the 
second containment,  observe that in that paper, the authors deal the the upper level sets of the form $\{\zz: g_n(\zz) \geq 0\}$, and develop an inner approximation for the level set of the pre-limit using the level sets at $(\alpha+\varepsilon)$ for $g^*$. However, in the present problem, we  deal with lower level sets $\{\zz: g_n(\zz) \leq 0\}$, and therefore, get the same approximation for $(\alpha-\varepsilon)$ level sets of $g^*$ instead. To see this, note that in their proof, one can replace the sequence $\beta_n\nearrow (\alpha+\varepsilon)$ by $\beta_n\searrow (\alpha-\varepsilon)$, keeping rest of the steps the same. \qed

\noindent \textbf{Proof of Lemma~\ref{lem:compact_containment}: }
Note that since $\mathcal S^*$ is a non-vacuous limiting set, there exists $\xx_0$ for which $\inf_{\mathcal \zz\in \mathcal E}\{\|\zz\|: \ \zz\in \Upsilon(\xx_0)\} >0$. Using the homogeneity of $\varphi^*$ furnished by Lemma~\ref{lem:RVs}, $I^*(\xx_0) = \inf_{\zz\in \mathcal E}\{\varphi^*(\zz): \ \zz\in \Upsilon^*(\xx_0)\}^{-1}<\infty$, and thus, $\mathcal Y^*$ is a contained in a compact set. 

Recall that the feasible region for the CCP is $\Lev_1(I_t)$. To verify that the solution sets $\mathcal Y_\alpha^*$ are contained in the \textit{same} compact set for all $\alpha<\alpha_0,$ note that the continuous convergence of $I_n \to I^*$ implies that the level sets of $I_t$ satisfy  $[\Lev_{1}(I_t)\cap B_M] \subset [\Lev_1(I^*)\cap B_M]^{1+\delta}$ (see \cite{rockafellar2009variational}, Proposition 7.7). Now, note the following: for $(M,\delta)>0$ ($M$ to be chosen imminently), 
\[
\Lev_1(I_t) = (\Lev_{1}(I_t)\cap B_M) \cup B_M^c \subseteq [\Lev_1(I^*)\cap B_M]^{1+\delta}\cup B_M^c
\]
for all $t\geq t_0$.
Therefore, $\mathcal Y_{\alpha}^* \subset \arg\min\{ c(\xx): \xx \in [\Lev_1(I^*)\cap B_M]^{1+\delta}\cup B_M^c\}$. Since $c(\cdot)$ is positively homogeneous, we can drop $B_M^c$ from the above arg-min if $M>M_0$ for some $M_0$. This suggests that for all $t>t_0$,  $\mathcal Y_{\alpha}^*$ are contained the same compact set $\text{cl}([\Lev_1(I^*)\cap B_M]^{1+\delta})$.\qed

\noindent \textbf{Proof of Lemma~\ref{lem:extended_functional}: }
It is sufficient to show that
for every $\xx\in \Lev_1(I^*)$, there exists $\xx_t\to \xx$ such that for all sufficiently large $t$, $\xx_t \in \Lev_1(I_t)$. First, note that $I^*(\cdot)$ by definition is a homogeneous function, and therefore, for any $c>0$, $I(c\xx) = c^k \xx$ for the appropriate $k$. Now, consider the sequence $\xx_t = c_t \xx$, where $c_t \to 1$, such that $c^k_t <1$ for all $t$. Note that $I(\xx_t) \leq c_t^k $. 
$I_t$ converges to $I^*$  compactly (see \cite{rockafellar2009variational}, Theorem 7.14),
\[
\sup_{t \geq t_0}|I_t(\xx_t) - I(\xx_t)| \to 0.
\]
Fix a $\delta>0$. Then, whenever $t\geq t_0$, the above implies that $|I_t(\xx_t) - I(\xx_t)| \leq  \delta/2$. Further, for $t\geq t_1,$ the choice of $c_t$ implies that $I(\xx_t) \leq 1-\delta/2$. Combining these together yields that $I_t(\xx_t) \leq 1$ for all $t\geq \max\{t_0,t_1\}$, and therefore $\xx_t\to \xx$ and $\xx_t\in \Lev_1(I_t)$ for all $t$ large enough.\qed



\noindent \textbf{Proof of Lemma~\ref{lem:solution_of_s}: }Recall that $f$ is a smooth, convex function whose minimum occurs at $x=1$. Therefore, $f^\prime(0)=f(0)=0$.
Let
\[
f_1(t) = f\left(\frac{1-t}{1-t/s}\right)
\]
Expanding  the expression in Lemma~\ref{lem:solution_of_s} about $t_0=1$ using a Taylor series: 
\[
\eta = t \frac{f(s)}{s} +\frac{t^2}{2} \frac{\left(1-1/s\right)^2}{(1-t/s)}f^{\prime\prime}_1(\tilde t), \text{ where } \tilde t\in\left[\frac{1-t}{1-t/s},1\right].
\]
Simplifying the above expression, $s(t)$ solves the equation
\begin{equation}\label{eqn:s_solution}
    w_t(s) = \frac{\eta}{t} - g^{\leftarrow}(s) - \frac{t}{2} \frac{\left(1-1/s\right)^2}{(1-t/s)}f^{\prime\prime}_1(\tilde t) = 0.
\end{equation}
To proceed, note first that the function $f_{1}$ is twice continuously differentiable in $t$, and therefore, uniformly over $s\in [1,\infty)$, $f_1^{\prime\prime}$ is bounded. 

Consequently, for $\varepsilon$, there exists a $t_0$ such that for all $t>t_0$, at $s_1 = g(\eta/t+\varepsilon)$, $w_t(s_1) < 0$. Next, fix an $\varepsilon>0$ and let $s_2 = g(\eta/t -\varepsilon)$ and observe that there exists a $t_1$, such that for all $t<t_1$, $w_t(s_2) >0$. Using the intermediate value theorem, for all $t>t_0$, there then exists an $\varepsilon_0\in[-\varepsilon,\varepsilon]$ such that $w_t(g(\eta/t -\varepsilon_0)) = 0$. This implies that $s(t) = g(\eta/t +\kappa(t))$ where $\kappa(t) \to 0$.\qed

\noindent \textbf{Proof of Lemma~\ref{lem:var_asymptotic}: } 
For ease of notation, define $g^*(\xx,\zz) = \max_k \eta_k g_k^*(\xx,\zz)$.
A consequence from the proof of Propositions~\ref{prop:prob-convergence_LT}-\ref{prop:prob-convergence_HT} is the following: for any $u>0$, the following approximations hold uniformly over $\xx$ in  compact sets:
\begin{align*}
    \log \mathbb P\left[\max_{k}g_{k,t}(\xx,\xxi_t) \geq u\right] &\sim  \lambda(t) I(\xx,u)\quad \quad\text{ when Assumption ($\mathscr{L}$) holds and }\\
     \mathbb P\left[\max_{k}g_{k,t}(\xx,\xxi_t) \geq u\right] &\sim  \bar F(t) I(\xx,u) \quad \quad \text{ when Assumption ($\mathscr{H}$) holds and }
\end{align*}
where 
\begin{align}\label{eqn:gen_prob_asympt}
    I(\xx,u) = \begin{cases}
        \left[\inf\{\varphi^*(\zz) : g^*(\xx,\zz) \geq u\}\right]^{-1} \quad \quad\ \text{ when Assumption ($\mathscr{L}$) holds }\\
        \int_{g^*(\xx,\zz)\geq u} \varphi^*(\zz) d\zz \quad \quad \quad \quad\quad \quad\quad \text{ when Assumption ($\mathscr{H}$) holds and }
    \end{cases}
\end{align}
With $t=s_\alpha$, we have $\lambda(t) = \log \alpha$ and $\bar F(t) = \alpha$. Additionally, with $u=v_\alpha(\xx)$, the left hand side of \eqref{eqn:gen_prob_asympt} equals $\alpha.$ Consequently, $I(\xx,v_\alpha(\xx)) \sim 1$ uniformly over $\xx$ in compact sets. Noting that $I(\xx,u)$ is continuous and decreasing in $u$ and setting $I(\xx,u^*(\xx)) = 1$, it must be the case that $v_\alpha(\xx) = (1+o(1)) u^*(\xx) $\qed

\noindent\textbf{Proof of Lemma~\ref{lem:I_value}: } Note that any $\xx\in  \mathcal Y^*$ satisfies the condition that $c(\xx^*) = v^*$ by definition and is feasible $\{c(\xx)\leq v^*\}$. Therefore, $\min\{I(\xx) : c(\xx) \leq v^*\}\leq 1$. It suffices to show that such an $\xx$ is in fact in the arg-min of the previous problem. 
 Note the following sequence of equalities: for any $t>0$
\begin{align*}
    v_t& = \min\{c(\xx) : I(\xx) \leq t\} =  \min\{c(\xx) : t^{-1}I(\xx) \leq 1\} = \min\{c(\xx): I(t^{\alpha/r} \xx) \leq 1\}\\
    &= \min\{c(t^{-\alpha/r}\pp) : I(\pp)\leq1\}\  \text{ setting } t^{\alpha/r}\xx=p \text{ above}\\
    & =t^{-\alpha/r}v_1 = t^{-\alpha/r}v^*
\end{align*}
Fix $t<1$. Then the above implies that whenever $I(\xx) = t$, $c(\xx) >v^*$. Thus, there cannot exist $\xx$ such that $I(\xx)<1$ and $c(\xx) \leq v^*$. Consequently $\xx\in\mathcal Y^*$ is actually in $\arg\min\{I(\xx) :c(\xx) \leq v^*\}$\qed

\end{document}